\documentclass[a4paper,10pt]{article}

\usepackage[ascii]{inputenc}
\usepackage{amsmath,amsfonts,amssymb,amsthm}
\usepackage{comment}

\newtheoremstyle{theorem}{5pt}{5pt}{\itshape}{}{\bfseries}{.}{.5em}{}

\theoremstyle{theorem}
\newtheorem{theorem}{Theorem}
\newtheorem{lemma}[theorem]{Lemma}
\newtheorem{corollary}[theorem]{Corollary}
\newtheorem{proposition}[theorem]{Proposition}

\usepackage{titlesec}
\titleformat{\section}
{\Large\bfseries}{\thesection}{.7em}{\large}
\titlespacing*{\section}{0pt}{3.5ex plus 1ex minus .2ex}{2.3ex plus .2ex}
\titleformat{\subsection}
{\normalfont\bfseries}{\thesubsection}{.5em}{\normalsize}
\titlespacing*{\subsection}{0pt}{3.5ex plus 1ex minus .2ex}{2.3ex plus .2ex}

\usepackage[font=small,labelfont=bf,margin=1.5em,labelsep=period]{caption}

\allowdisplaybreaks[3]

\begin{document}

\title{Exponential Sums Related to Maass Forms}
\author{Jesse J\"a\"asaari\footnote{University of Helsinki, Department of Mathematics and Statistics, P.O.Box 68 (Gustaf H\"allstr\"omin katu 2b), FI-00014 University of Helsinki, FINLAND} { }and
Esa V\!. Vesalainen\footnote{University of Jyvaskyla, Department of Mathematics and Statistics, P.O.Box 35 (MaD) FI-40014 University of Jyvaskyla, FINLAND}}
\date{}
\maketitle

\vspace*{-.5cm}
\begin{abstract}
We estimate short exponential sums weighted by the Fourier coefficients of a Maass form. This requires working out a certain transformation formula for non-linear exponential sums, which is of independent interest. We also discuss how the results depend on the growth of the Fourier coefficients in question. As a byproduct of these considerations, we can slightly extend the range of validity of a short exponential sum estimate for holomorphic cusp forms.

The short estimates allow us to reduce smoothing errors. In particular, we prove an analogue of an approximate functional equation previously proven for holomorphic cusp form coefficients.

As an application of these, we remove the logarithm from the classical upper bound for long linear sums weighted by Fourier coefficients of Maass forms, the resulting estimate being the best possible. This also involves improving the upper bounds for long linear sums with rational additive twists, the gains again allowed by the estimates for the short sums. Finally, we shall use the approximate functional equation to bound somewhat longer short exponential sums.
\end{abstract}

\section{Introduction and the main results}

\subsection{Maass forms}

Let $\psi$ be a Maass form for the full modular group, corresponding to an eigenvalue $1/4+\kappa^2$ of the hyperbolic Laplacian, and with the Fourier expansion
\[\psi\!\left(x+yi\right)=y^{1/2}\sum_{n\neq0}t(n)\,K_{i\kappa}\!\left(2\pi\left|n\right|y\right)e\!\left(nx\right),\]
where $x\in\mathbb R$ and $y\in\mathbb R_+$. We may assume without loss of generality that $\psi$ is even or odd, i.e.\ that $t(-n)=t(n)$ for all $n\in\mathbb Z_+$, or that $t(-n)=-t(n)$ for all $n\in\mathbb Z_+$. For standard references on Maass forms we refer to \cite{Iwaniec2002, Motohashi1997}.

The Fourier coefficients $t(n)$ satisfy a bound of the kind
\[t(n)\ll n^{\vartheta+\varepsilon}\]
for some $\vartheta\in\left[0,\infty\right[$. The best known exponent $\vartheta=\frac7{64}$ is due to Kim and Sarnak \cite{Kim--Sarnak2003}. The Ramanujan--Petersson conjecture for Maass forms declares that $\vartheta=0$ is admissible. On average, the Fourier coefficients are of constant size. In particular, we have a Rankin--Selberg type estimate for the Fourier coefficients. One such result is the following (see e.g. \cite{Iwaniec2002}, Chapter 8):
\begin{equation}\label{rankin--selberg}
\int\limits_0^1\left|\sum_{n\leqslant M}t(n)\,e(n\alpha)\right|^2\mathrm d\alpha
=\sum_{n\leqslant M}\left|t(n)\right|^2=A\,M+O(M^{7/8}),
\end{equation}
where $A$ is a positive real constant depending on $\psi$.

\subsection{Objects of study and motivation}

In the following we will consider linear exponential sums of the form
\[\sum_{M\leqslant n\leqslant M+\Delta}t(n)\,e(n\alpha),\]
where $M\in\left[1,\infty\right[$, $\Delta\in\left[1,M\right]$ and $\alpha\in\mathbb R$. When $\Delta=o(M)$, we call such sums short.

The reasons for considering such sums are manifold. First of all, the Fourier coefficients $t(n)$ are interesting mathematical objects which are not as well understood as one might wish. The exponential sums above contain all the information about the Fourier coefficients and thus provide an interesting window into their behaviour. The study of Maass forms and their $L$-functions also naturally leads to exponential sums weighted by the corresponding Fourier coefficients of which the above linear sums are an important special case.

When $\alpha$ is a rational number $h/k$, the problem of estimating long sums with $\Delta=M$ is very analogous to classical problems in analytic number theory, such as the problems of estimating the error terms in the circle and Dirichlet divisor problems. Furthermore, the problem of estimating such sums with $\Delta=o(M)$ provides an analogue for problems such as studying the behaviour of the afore-mentioned error terms in short intervals. For further information about these classical topics, see e.g. Chapter 13 of \cite{Ivic2003} or \cite{Tsang2010}.

Finally, good estimates for the short exponential sums above can sometimes be used to reduce smoothing error. An example of such an application is given e.g. by Theorems \ref{improved-estimate} and \ref{approximate-functional-equation}
below.

For holomorphic cusp forms, short exponential sums have been studied by Jutila \cite{Jutila1987b}, and the best known bounds are due to Ernvall-Hyt\"onen and Karppinen \cite{Ernvall-Hytonen--Karppinen2008, Ernvall-Hytonen2008}.

It is interesting to study how sensitive the arguments used for holomorphic cusp forms are to the value of $\vartheta$. In a sense, the strictly positive value of $\vartheta$ is the main difference between the holomorphic and non-holomorphic cases: in both cases one applies heavily the corresponding Voronoi summation formula, and even though the Voronoi summation formulae have a different appearance, what remains after the Bessel functions have been cashed in in terms of their asymptotics is similar.

\subsection{The results: Bounds for short exponential sums with applications}

The following is a Maass form analogue of the related estimate for holomorphic cusp forms due to Ernvall-Hyt\"onen and Karppinen, Theorem 5.5 in \cite{Ernvall-Hytonen--Karppinen2008}. The proof is based on techniques analogous to those in \cite{Ernvall-Hytonen--Karppinen2008}.
\begin{theorem}\label{short-estimate}
Let $M\in\left[1,\infty\right[$ and let $\Delta\in\left[1,M\right]$ be such that $\Delta\ll M^{2/3}$. Then
\[\sum_{M\leqslant n\leqslant M+\Delta}t(n)\,e(n\alpha)\ll\Delta^{1/6-\vartheta}\,M^{1/3+\vartheta+\varepsilon},\]
uniformly for $\alpha\in\mathbb R$. This is better than estimating via absolute values when $M^{2/(5+6\vartheta)}\ll\Delta\ll M^{2/3}$.
\end{theorem}
When $\Delta=M^{2/3}$ this gives the upper bound $\ll M^{\vartheta/3+4/9+\varepsilon}$, and so splitting a longer sum into sums of this length and estimating the subsums separately gives the following bound for longer sums.
\begin{corollary}\label{short-estimate-trivial-corollary}
Let $M\in\left[1,\infty\right[$ and let $\Delta\in\left[1,M\right]$ be such that $M^{2/3}\ll\Delta\ll M$. Then
\[\sum_{M\leqslant n\leqslant M+\Delta}t(n)\,e(n\alpha)\ll\Delta\,M^{\vartheta/3-2/9+\varepsilon}.\]
This is better than the bound $\ll M^{1/2+\varepsilon}$ when $M^{2/3}\ll\Delta\ll M^{(13-6\vartheta)/18}$.
\end{corollary}

Actually, Theorem \ref{short-estimate} is valid for a slightly larger range of $\Delta$ than Theorem 5.5 in \cite{Ernvall-Hytonen--Karppinen2008} is. In fact, with a minor modification \cite{Ernvall-Hytonen--communication}, the proof of Theorem 5.5 of \cite{Ernvall-Hytonen--Karppinen2008} can be easily modified to give the analogous result for holomorphic cusp forms:
\begin{theorem}\label{holomorphic-short-estimate}
Let us consider a fixed holomorphic cusp form of weight $\kappa\in\mathbb Z_+$ for the full modular group with the Fourier expansion
\[\sum_{n=1}^\infty a(n)\,n^{(\kappa-1)/2}\,e(nz)\]
for $z\in\mathbb C$ with $\Im z>0$. Also, let $M\in\left[1,\infty\right[$, $\Delta\in\left[1,M\right]$, and let $\alpha\in\mathbb R$. If $\Delta\ll M^{2/3}$, then
\[\sum_{M\leqslant n\leqslant M+\Delta}a(n)\,e(n\alpha)\ll\Delta^{1/6}\,M^{1/3+\varepsilon},\]
where the implicit constant depends only on the underlying cusp forms and $\varepsilon$. Similarly, if $M^{2/3}\ll\Delta$, then
\[\sum_{M\leqslant n\leqslant M+\Delta}a(n)\,e(n\alpha)\ll\Delta\,M^{-2/9+\varepsilon}.\]
\end{theorem}

The proof of Theorem \ref{short-estimate} depends on an estimate for short non-linear sums, analogous to Theorem 4.1 in \cite{Ernvall-Hytonen--Karppinen2008}. Fortunately, the proof in \cite{Ernvall-Hytonen--Karppinen2008} works almost verbatim for Maass forms and we shall indicate the differences later. On the other hand, the proof of the non-linear estimate requires a transformation formula of a certain shape for smoothed exponential sums, and this particular result does not seem to have been worked out before yet. Thus, in Section \ref{transformation-section}, we will give an analogue of the relevant Theorem 3.4 of Jutila's monograph \cite{Jutila1987a}, which considers smooth sums with holomorphic cusp form coefficients, with full details for Maass forms. An analogue of Theorem 3.2 of \cite{Jutila1987a} has been given by Meurman in \cite{Meurman1987}.

The following estimates provide a concrete example of how estimates for short sums allow one to reduce smoothing errors thereby leading to improved upper bounds.
\begin{theorem}\label{improved-estimate}
Let $M\in\left[1,\infty\right[$, $h\in\mathbb Z$, $k\in\mathbb Z_+$ and $\left(h,k\right)=1$. Also, let $\delta\in\left]0,1/2\right[$ and assume that $k\ll M^{1/2-\delta}$. Then
\[\sum_{n\leqslant M}t(n)\,e\!\left(\frac{nh}k\right)\ll_\delta k^{2/3}\,M^{1/3+\vartheta/3+\varepsilon}.\]
When $M^{3/(5+6\vartheta)-1/2+\vartheta}\ll k\ll M^{5/18+\vartheta/3}$, we have the upper bound
\[\sum_{n\leqslant M}t(n)\,e\!\left(\frac{nh}k\right)\ll k^{(1-6\vartheta)/(4-6\vartheta)}\,M^{3/(8-12\vartheta)+\varepsilon}.\]
Similarly, for $M^{5/18+\vartheta/3}\ll k\ll M^{1/2-\varepsilon}$, we have
\[\sum_{n\leqslant M}t(n)\,e\!\left(\frac{nh}k\right)\ll k^{2/3}\,M^{7/27+\vartheta/9+\varepsilon}.\]
\end{theorem}
The case $k=1$ was considered by Hafner and Ivi\'c \cite{Hafner--Ivic1989} who essentially obtained the bound $\ll M^{1/3+\vartheta/3}$.
Similar reduction for certain ranges of $k$ in the case of holomorphic cusp forms have recently been proved by Vesalainen \cite{Vesalainen2014}. The proof is analogous to the approach of \cite{Ivic2003}.

It is of interest to note here that for small enough $k$, the rationally twisted sum has on average (in the mean square sense) the order of magnitude $k^{1/2}\,M^{1/4}$. This kind of result was first proven by Cram\'er \cite{Cramer1922} for the error term in the Dirichlet divisor problem. Jutila \cite{Jutila1985} extended this to the divisor problem with rational additive twists, and in \cite{Jutila1987a} Jutila proved the analogous result for holomorphic cusp forms. We shall elaborate on this in the last section.

Theorem \ref{holomorphic-short-estimate} allows us to improve Theorem 1 from \cite{Vesalainen2014} in the range $k\gg M^{1/4}$:
\begin{corollary}\label{holomorphic-corollary}
Let $a(n)$ be the Fourier coefficients of a holomorphic cusp form as in Theorem \ref{holomorphic-short-estimate}. Then, for coprime integers $h$ and $k$ with $M^{1/10}\ll k\ll M^{5/18}$, we have
\[\sum_{n\leqslant M}a(n)\,e\!\left(\frac{nh}k\right)\ll k^{1/4}\,M^{3/8+\varepsilon},\]
and for $M^{5/18}\ll k\ll M^{1/2-\varepsilon}$, we have
\[\sum_{n\leqslant M}a(n)\,e\!\left(\frac{nh}k\right)\ll k^{2/3}\,M^{7/27+\varepsilon}.\]
\end{corollary}

\subsection{The results: an approximate functional equation and applications}

Wilton \cite{Wilton1933} proved an approximate functional equation for exponential sums involving the divisor function. Jutila \cite{Jutila1985} extended this to sums with additive twists, and in \cite{Jutila1987b} he proved an analogue for holomorphic cusp forms. In \cite{Ernvall-Hytonen2008} Ernvall-Hyt\"onen improved the error term. The following is an analogue of Ernvall-Hyt\"onen's result, and the proof is analogous to that in \cite{Ernvall-Hytonen2008}.
We write $\overline h$ for an integer such that $h\overline h\equiv 1\pmod k$. Also, to simplify the notation, we write
\[T(M,\Delta;\alpha)=\sum_{M\leqslant n\leqslant M+\Delta}t(n)\,e(n\alpha).\]
\begin{theorem}\label{approximate-functional-equation}
Let $\alpha\in\mathbb R$ have the rational approximation $\alpha=\frac hk+\eta$, where $h$ and $k$ are coprime integers with $1\leqslant k\leqslant M^{1/4}$ and $\left|\eta\right|\leqslant k^{-1}\,M^{-1/4}$. Furthermore, let $M\in\left[1,\infty\right[$ and $\Delta\in\left[1,M\right]$. If $k^2\,\eta^2\,M\gg1$, then
\[\frac{T(M,\Delta;\alpha)}{M^{1/2}}
=\frac{T(k^2\,\eta^2\,M,k^2\,\eta^2\,\Delta;\beta)}{(k^2\,\eta^2\,M)^{1/2}}
+O\bigl((k^2\,\eta^2\,M)^{\vartheta/2-1/12+\varepsilon}\bigr),\]
where $\displaystyle{\beta=-\frac{\overline h}k-\frac1{k^2\eta}}$.
\end{theorem}

Wilton \cite{Wilton1929} proved that for the normalized Fourier coefficients $a(n)$ of a fixed holomorphic cusp form,
\[\sum_{n\leqslant M}a(n)\,e(n\alpha)\ll M^{1/2}\,\log M,\]
uniformly in $\alpha\in\mathbb R$. The Rankin--Selberg bound on the mean square of Fourier coefficients implies that
\[\int\limits_0^1\left|\sum_{n\leqslant M}a(n)\,e(n\alpha)\right|^2\mathrm d\alpha
=\sum_{n\leqslant M}\left|a(n)\right|^2=A\,M+O(M^{3/5}),\]
for a certain positive real constant $A$ depending on the underlying cusp form,
and so at most the logarithm can be removed from Wilton's estimate, and this indeed was done by Jutila \cite{Jutila1987b}.
For Maass forms, the estimate analogous to Wilton's was proved by Epstein, Hafner and Sarnak \cite{Epstein--Hafner--Sarnak1985, Hafner1987}.
The following estimate is an analogue of Jutila's logarithm removal, and its proof is largely analogous to the arguments in \cite{Jutila1987b}.
\begin{theorem}\label{logarithm-removal}
We have
\[\sum_{n\leqslant M}t(n)\,e(n\alpha)\ll M^{1/2},\]
uniformly in $\alpha\in\mathbb R$.
\end{theorem}
\noindent This is sharp in view of \eqref{rankin--selberg}.

With the approximate functional equation at hand, we may prove further estimates for short sums.
\begin{theorem}\label{longer-short-estimates}
Let $M\in\left[1,\infty\right[$ and $\Delta\in\left[1,M\right]$ with $M^{2/3}\ll\Delta\ll M^{3/4}$, and let $\alpha\in\mathbb R$. Then
\[\sum_{M\leqslant n\leqslant M+\Delta}t(n)\,e(n\alpha)
\ll
M^{3/8+(3+12\vartheta)/(32+48\vartheta)+\varepsilon}
+\Delta\,M^{-1/4+3\vartheta/(32+48\vartheta)+\varepsilon}.
\]
In particular, for $\vartheta=7/64$ we have
\[\sum_{M\leqslant n\leqslant M+\Delta}t(n)\,e(n\alpha)\ll M^{585/1192+\varepsilon}+\Delta\,M^{-575/2384+\varepsilon},\]
which is better than the estimate $\ll M^{1/2+\varepsilon}$ for $\Delta\ll M^{1767/2384}$,
and for $\vartheta=0$, we have
\[\sum_{M\leqslant n\leqslant M+\Delta}t(n)\,e(n\alpha)\ll M^{15/32+\varepsilon}+\Delta\,M^{\varepsilon-1/4}.\]
\end{theorem}

\subsection{$\Omega$-results}

Finally, it is naturally interesting to consider what are the limits of estimating short sums.  In \cite{Ernvall-Hytonen2009b} Ernvall-Hyt\"onen proved that, if $d\in\mathbb Z_+$ is a fixed integer such that $t(d)\neq0$, then
\[\sum_{M\leqslant n\leqslant M+\Delta}t(n)\,e(n\alpha)\,w(n)\asymp\Delta\,M^{-1/4},\]
where $w$ is a suitable weight function, $\alpha=\sqrt{d}/\sqrt{M}$, and $M^{1/2+\varepsilon}\ll\Delta\leqslant d^{-1/2}\,M^{3/4}$. This immediately implies that for this range of lengths $\Delta$,
\begin{equation}\label{gl2-resonance}
\sum_{M\leqslant n\leqslant M+\Delta}t(n)\,e\!\left(\frac{n\sqrt d}{\sqrt M}\right)=\Omega(\Delta\,M^{-1/4}).
\end{equation}
This result also has counterparts for the divisor function and Fourier coefficients of holomorphic cusp forms in the papers of Ernvall-Hyt\"onen and Karppinen \cite{Ernvall-Hytonen--Karppinen2008} and Ernvall-Hyt\"onen \cite{Ernvall-Hytonen2008, Ernvall-Hytonen2009a}.
The above $\Omega$-result implies that the bound
\[\sum_{M\leqslant n\leqslant M+\Delta}t(n)\,e(n\alpha)\ll M^{1/2}\]
from Theorem \ref{logarithm-removal} is sharp for $M^{3/4}\ll\Delta\ll M$.

For sums of length $\Delta\ll M^{1/2}$, it turns out that square root cancellation is the best that could be hoped for. Essentially, combining the truncated Voronoi identity of Meurman \cite{Meurman1988} with the arguments of Jutila \cite{Jutila1984}, one gets the following mean square asymptotics
\[\int\limits_M^{2M}\,\left|\sum_{x\leqslant n\leqslant x+\Delta}t(n)\right|^2\mathrm dx\asymp\Delta\,M,\]
for $M^{2\vartheta+\varepsilon}\ll\Delta\ll M^{1/2-\varepsilon}$.
In fact, a sharper result could be obtained, but this is enough for the relevant $\Omega$-result. The paper \cite{Jutila1984} actually considered the behaviour of the error terms in the Dirichlet divisor problem and the second moment for the Riemann $\zeta$-function in short intervals, but the proof for the divisor function carries through fairly easily for Fourier coefficients of holomorphic cusp forms or Maass forms. In the last section, we will discuss the second moments with more details, and add here only that for holomorphic cusp forms second moments of rationally additively twisted short sums have been considered in the works \cite{Ernvall-Hytonen2011, Ernvall-Hytonen2015, Vesalainen2014}.

We would like to emphasize that there are reasons to believe that even if the best possible upper bounds conform to the above $\Omega$-results, they are likely to be very difficult to prove. For example, the conjectural upper bounds
\[\sum_{M\leqslant n\leqslant M+M^{1/2}}a(n)\,e(n\alpha)\ll M^{1/4+\varepsilon}\]
and
\[\sum_{M\leqslant n\leqslant M+M^{1/2}}t(n)\,e(n\alpha)\ll M^{1/4+\varepsilon}\]
would be analogous to the conjectural upper bound
\[\Delta(M+M^{1/2})-\Delta(M)\ll M^{1/4+\varepsilon}\]
for the error term in the Dirichlet divisor problem. Jutila \cite{Jutila1983} has proved that if the last estimate is true, then Riemann's zeta-function satisfies the bound $\zeta(1/2+it)\ll t^{3/20+\varepsilon}$ on the critical line, and this exponent $3/20$ is better than the best known exponent 53/342 due to Bourgain \cite{Bourgain2014} or Huxley's exponent 32/205 \cite{Huxley2005}.

\begin{figure}[h]
\begin{center}
\setlength{\unitlength}{1cm}
\begin{picture}(12.3,6.5)(-.8,-.5)
\qbezier(0,-.25)(0,2)(0,5.5)
\qbezier(-.25,0)(5,0)(10.5,0)
\qbezier(0,5.5)(.025,5.4)(.07,5.3)
\qbezier(0,5.5)(-.025,5.4)(-.07,5.3)
\qbezier(10.5,0)(10.4,.025)(10.3,.07)
\qbezier(10.5,0)(10.4,-.025)(10.3,-.07)
\put(-.105,5.61){$\beta$}
\put(10.6,-.08){$\gamma$}
\qbezier(-.05,5)(0,5)(.05,5)
\qbezier(-.05,2.5)(0,2.5)(.05,2.5)
\qbezier(-.05,4)(0,4)(.05,4)
\qbezier(-.05,4.44)(0,4.44)(.05,4.44)
\qbezier(-.05,4.69)(0,4.69)(.05,4.69)
\put(-.55,4.95){$\scriptstyle{1/2}$}
\put(-.82,4.64){$\scriptstyle{15/32}$}
\put(-.55,2.44){$\scriptstyle{1/4}$}
\put(-.55,3.95){$\scriptstyle{2/5}$}
\put(-.55,4.38){$\scriptstyle{4/9}$}
\qbezier(4,.05)(4,0)(4,-.05)
\qbezier(5,.05)(5,0)(5,-.05)
\qbezier(7.5,.05)(7.5,0)(7.5,-.05)
\qbezier(10,.05)(10,0)(10,-.05)
\qbezier(6.66,.05)(6.66,0)(6.66,-.05)
\qbezier(6.91,.05)(6.91,0)(6.91,-.05)
\qbezier(7.19,.05)(7.19,0)(7.19,-.05)
\put(3.88,-.4){$\textstyle{\frac25}$}
\put(4.88,-.4){$\textstyle{\frac12}$}
\put(6.54,-.4){$\textstyle{\frac23}$}
\put(7.38,-.4){$\textstyle{\frac34}$}
\put(9.93,-.27){$\scriptstyle1$}
\put(7.01,-.4){$\textstyle{\frac{23}{32}}$}
\put(6.66,.25){$\textstyle{\frac{199}{288}}$}
\linethickness{.1mm}
\multiput(4,0)(0,.125){32}{\qbezier(0,0)(0,.025)(0,.05)}
\multiput(5,0)(0,.125){20}{\qbezier(0,0)(0,.025)(0,.05)}
\multiput(6.66,0)(0,.125){36}{\qbezier(0,0)(0,.025)(0,.05)}
\multiput(7.5,0)(0,.125){40}{\qbezier(0,0)(0,.025)(0,.05)}
\multiput(7.19,0)(0,.125){38}{\qbezier(0,0)(0,.025)(0,.05)}
\multiput(6.91,.625)(0,.125){33}{\qbezier(0,0)(0,.025)(0,.05)}
\linethickness{.1mm}
\multiput(0,2.5)(.125,0){40}{\qbezier(0,0)(.025,0)(.05,0)}
\multiput(0,4)(.125,0){32}{\qbezier(0,0)(.025,0)(.05,0)}
\multiput(0,4.44)(.125,0){53}{\qbezier(0,0)(.025,0)(.05,0)}
\multiput(0,4.69)(.125,0){55}{\qbezier(0,0)(.025,0)(.05,0)}
\multiput(0,5)(.125,0){60}{\qbezier(0,0)(.025,0)(.05,0)}
\linethickness{.25mm}
\qbezier(0,0)(2,2)(4,4)
\qbezier(7.5,5)(8.75,5)(10,5)
\qbezier(4,4)(5.33,4.22)(6.66,4.44)
\qbezier(7.19,4.69)(7.345,4.845)(7.5,5)
\qbezier(6.91,4.69)(7.05,4.69)(7.19,4.69)
\qbezier(6.66,4.44)(6.785,4.565)(6.91,4.69)
\multiput(0,0)(.25,.125){20}{\qbezier(0,0)(.075,.0375)(.15,.075)}
\multiput(5,2.5)(.25,.25){10}{\qbezier(0,0)(.075,.075)(.15,.15)}
\end{picture}
\end{center}
\caption{\label{figure1}Estimates for short linear sums related to holomorphic cusp forms: a point $\left\langle\gamma,\beta\right\rangle$ on the solid thick line or on the dotted thick line signifies an estimate of the form\\[3mm]
\hspace*{.8em}$\displaystyle{\sum_{M\leqslant n\leqslant M+M^\gamma}a(n)\,e(n\alpha)\ll M^{\beta+\varepsilon},\quad\text{or}\quad\sum_{M\leqslant n\leqslant M+M^\gamma}a(n)\,e(n\alpha)=\Omega(M^\beta),}$\\[2mm]
respectively.
The first estimate holds uniformly in $\alpha\in\mathbb R$. The first upper bound segment comes from estimating by absolute values with Deligne's estimate for individual Fourier coefficients from \cite{Deligne1974}, the second and third segments from Theorem 5.5 of \cite{Ernvall-Hytonen--Karppinen2008} and Theorem \ref{holomorphic-short-estimate}, the fourth and fifth segments from Theorem 5.16 in \cite{Ernvall-Hytonen--Karppinen2008}, and the horizontal sixth segment from Wilton's and Jutila's estimates \cite{Wilton1929, Jutila1987b}. The first lower bound segment follows from the work of Jutila \cite{Jutila1984}, and the second follows from Theorem 6.1 in \cite{Ernvall-Hytonen--Karppinen2008}. The upper bounds in the range $M^{23/32+\varepsilon}\ll\Delta\ll M$ are sharp.}
\end{figure}

\begin{figure}[h]
\begin{center}
\setlength{\unitlength}{1cm}
\begin{picture}(12.3,6.5)(-.8,-.5)
\qbezier(0,-.25)(0,2)(0,5.5)
\qbezier(-.25,0)(5,0)(10.5,0)
\qbezier(0,5.5)(.025,5.4)(.07,5.3)
\qbezier(0,5.5)(-.025,5.4)(-.07,5.3)
\qbezier(10.5,0)(10.4,.025)(10.3,.07)
\qbezier(10.5,0)(10.4,-.025)(10.3,-.07)
\put(-.105,5.61){$\beta$}
\put(10.6,-.08){$\gamma$}
\qbezier(-.05,5)(0,5)(.05,5)
\qbezier(-.05,2.5)(0,2.5)(.05,2.5)
\qbezier(-.05,4.64)(0,4.64)(.05,4.64)
\qbezier(-.05,4.81)(0,4.81)(.05,4.81)
\qbezier(-.05,4.91)(0,4.91)(.05,4.91)
\put(-.55,4.99){$\scriptstyle{1/2}$}
\put(-.55,2.44){$\scriptstyle{1/4}$}
\put(-.43,4.68){\raisebox{.45mm}{$\scriptstyle\star$}$\{$}
\qbezier(-.05,1.09)(0,1.09)(.05,1.09)
\put(-.70,1.04){$\scriptstyle{7/64}$}
\qbezier(3.53,.05)(3.53,0)(3.53,-.05)
\qbezier(5,.05)(5,0)(5,-.05)
\qbezier(7.5,.05)(7.5,0)(7.5,-.05)
\qbezier(10,.05)(10,0)(10,-.05)
\qbezier(6.66,.05)(6.66,0)(6.66,-.05)
\qbezier(6.76,.05)(6.76,0)(6.76,-.05)
\qbezier(7.32,.05)(7.32,0)(7.32,-.05)
\qbezier(7.41,.05)(7.41,0)(7.41,-.05)
\put(3.27,-.4){$\textstyle{\frac{64}{181}}$}
\put(4.88,-.4){$\textstyle{\frac12}$}
\put(9.93,-.27){$\scriptstyle1$}
\qbezier(2.18,.05)(2.18,0)(2.18,-.05)
\put(1.985,-.4){$\textstyle{\frac7{32}}$}
\put(6.55,0){$\underbrace{\phantom{abcdef}}_{\star\star}$}
\linethickness{.1mm}
\multiput(2.18,0)(0,.125){9}{\qbezier(0,0)(0,.025)(0,.05)}
\multiput(3.53,0)(0,.125){37}{\qbezier(0,0)(0,.025)(0,.05)}
\multiput(5,0)(0,.125){20}{\qbezier(0,0)(0,.025)(0,.05)}
\multiput(6.66,0)(0,.125){39}{\qbezier(0,0)(0,.025)(0,.05)}
\multiput(6.76,0)(0,.125){39}{\qbezier(0,0)(0,.025)(0,.05)}
\multiput(7.32,0)(0,.125){39}{\qbezier(0,0)(0,.025)(0,.05)}
\multiput(7.41,0)(0,.125){40}{\qbezier(0,0)(0,.025)(0,.05)}
\multiput(7.5,0)(0,.125){40}{\qbezier(0,0)(0,.025)(0,.05)}
\linethickness{.1mm}
\multiput(0,2.5)(.125,0){40}{\qbezier(0,0)(.025,0)(.05,0)}
\multiput(0,4.64)(.125,0){28}{\qbezier(0,0)(.025,0)(.05,0)}
\multiput(0,4.81)(.125,0){53}{\qbezier(0,0)(.025,0)(.05,0)}
\multiput(0,4.91)(.125,0){54}{\qbezier(0,0)(.025,0)(.05,0)}
\multiput(0,5)(.125,0){60}{\qbezier(0,0)(.025,0)(.05,0)}
\multiput(0,1.09)(.125,0){18}{\qbezier(0,0)(.025,0)(.05,0)}
\qbezier(0,0)(2,2)(4,4)
\qbezier(4,4)(5.33,4.22)(6.66,4.44)
\qbezier(7.19,4.69)(7.345,4.845)(7.5,5)
\qbezier(6.91,4.69)(7.05,4.69)(7.19,4.69)
\qbezier(6.66,4.44)(6.785,4.565)(6.91,4.69)
\linethickness{.25mm}
\multiput(2.25,1.125)(1,1){1}{\qbezier(0,0)(-.035,-.0175)(-.07,-.035)}
\multiput(5,2.5)(-.25,-.125){11}{\qbezier(0,0)(-.075,-.0375)(-.15,-.075)}
\multiput(5,2.5)(.25,.25){10}{\qbezier(0,0)(.075,.075)(.15,.15)}
\qbezier(0,1.09)(1.768,2.865)(3.53,4.64)
\qbezier(3.53,4.64)(5.095,4.725)(6.66,4.81)
\qbezier(6.66,4.81)(6.71,4.86)(6.76,4.91)
\qbezier(6.76,4.91)(7.04,4.91)(7.32,4.91)
\qbezier(7.32,4.91)(7.365,4.955)(7.41,5)
\qbezier(7.41,5)(8.705,5)(10,5)
\put(8,4){$\scriptstyle{\star\left\{\!\!\!\begin{array}{l}\scriptstyle{585/1192}\\[-.1mm]\scriptstyle{277/576}\\[-.1mm]\scriptstyle{64/181}\end{array}\right.}$}
\put(7.8,1.5){$\underbrace{\textstyle{\frac23\,\,\frac{58063}{85824}\,\,\frac{1745}{3284}\,\,\frac{1767}{2384}\,\,\frac34}}_{\star\star}$}
\end{picture}
\end{center}
\caption{\label{figure2}Estimates for short linear sums related to Maass forms: a point $\left\langle\gamma,\beta\right\rangle$ on the solid thick line or on the dotted thick line signifies an estimate of the form\\[3mm]
\hspace*{.8em}$\displaystyle{\sum_{M\leqslant n\leqslant M+M^\gamma}a(n)\,e(n\alpha)\ll M^{\beta+\varepsilon},\quad\text{or}\quad\sum_{M\leqslant n\leqslant M+M^\gamma}a(n)\,e(n\alpha)=\Omega(M^\beta),}$\\[2mm]
respectively, when $\vartheta$ is taken to be $7/64$.
The first estimate holds uniformly in $\alpha\in\mathbb R$.
The first upper bound segment comes from estimating by absolute values, the second and third from Theorem \ref{short-estimate}, the fourth and fifth segments from Theorem \ref{longer-short-estimates}, and the sixth horizontal segment from \cite{Epstein--Hafner--Sarnak1985, Hafner1987} and Theorem \ref{logarithm-removal}. The first lower bound segment follows from the arguments in \cite{Jutila1984} (but see also Theorem \ref{short-mean-square} below), the second lower bound segment comes from Theorem 2 in \cite{Ernvall-Hytonen2009b}. The upper bounds in the range $M^{3/4}\ll\Delta\ll M$ are sharp.
If $\vartheta$ can be taken to be zero, then this picture reduces to the one in Figure \ref{figure1}.}
\end{figure}

\subsection{Complements: uniformity and higher rank}

We would like to say a few words about Maass forms for $\mathrm{GL}(n)$. For them, a Voronoi summation formula exists and was implemented in \cite{Miller--Schmid2006, Miller--Schmid2008, Goldfeld--Li2006, Goldfeld--Li2008}. It has been applied to exponential sums weighted by Fourier coefficients of $\mathrm{GL}(n)$ Maass forms. As examples, we mention the works \cite{Miller2006, Li--Young2012, Ernvall-Hytonen2010, Ernvall-Hytonen--Jaasaari--Vesalainen2015, Godber2013}.
In particular, \cite{Miller2006} gives an upper bound for long linear sums in $\mathrm{GL}(3)$, and \cite{Ernvall-Hytonen2010, Ernvall-Hytonen--Jaasaari--Vesalainen2015} give higher rank analogues of the above mentioned $\Omega$-result \eqref{gl2-resonance}.

We have only considered a fixed cusp form for the full modular group. The dependence of the upper bound for long linear sums on the underlying cusp form has been considered in \cite{Li--Young2012} and \cite{Godber2013}. The discussion of the Farey and similar methods for holomorphic cusp forms in \cite{Huxley1996} also considers the depence on the underlying cusp forms, and the papers \cite{Meurman1987} and \cite{Meurman1988} consider the dependence on the underlying Maass form.

\subsection{Notation}

All the implicit constants are allowed to depend on the underlying Maass form, and $\varepsilon$, which denotes an arbitrarily small fixed positive number, which is not necessarily the same on each occurrence. Implicit constants depend also on chosen positive integers $J$ and $K$, when they appear.

The symbols $\ll$, $\gg$, $\asymp$, and $O$ are used for the usual asymptotic notation: for complex valued functions $f$ and $g$ in some set $\Omega$, the notation $f\ll g$ means that $\left|f(x)\right|\leqslant C\left|g(x)\right|$ for all $x\in\Omega$ for some implicit constant $C\in\mathbb R_+$. When the implicit constant depends on some parameters $\alpha,\beta,\ldots$, we use $\ll_{\alpha,\beta,\ldots}$ instead of mere $\ll$. The notation $g\gg f$ means $f\ll g$, and $f\asymp g$ means $f\ll g\ll f$.

The symbol $\sum_{a\leqslant n\leqslant b}'$ signifies summation over the integers $n$ with $a\leqslant n\leqslant b$, with possible terms corresponding to $a$ and $b$ halved if $a$ or $b$ is an integer. The symbol
\[\sum_{\substack{L\leqslant X\\\mathrm{dyadic}}}\]
signifies summation over the values $L=X$, $L=X/2$, $L=X/4$, \dots

Finally, the characteristic function of a set $B$ is denoted by $\chi_B$, and $e(x)$ denotes $e^{2\pi i x}$ for all $x\in\mathbb R$.

\section{The Voronoi type summation formula for Maass forms}

The main tool in the following is a Voronoi type summation formula for Maass forms with rational additive twists, proved by Meurman \cite{Meurman1988}. The following result is Theorem 2 in \cite{Meurman1988}.
\begin{theorem}\label{full-voronoi-for-maass-waves}
For a function $f\in C^1\!\left(\left[a,b\right]\right)$, where $a<b$ are positive real numbers, and for a positive integer $k$ and an integer $h$ coprime to $k$, we have
\begin{align*}
&\sum_{a\leqslant n\leqslant b}'t\!\left(n\right)e\!\left(\frac{nh}k\right)f\!\left(n\right)\\
&=\frac{\pi\,i}{k\,\sinh\pi\kappa}\sum_{n=1}^\infty t\!\left(n\right)e\!\left(\frac{-n\overline h}k\right)
\int\limits_a^b\!\left(J_{2i\kappa}\!\left(\frac{4\pi\sqrt{nx}}k\right)-J_{-2i\kappa}\!\left(\frac{4\pi\sqrt{nx}}k\right)\right)f\!\left(x\right)\mathrm dx\\
&\qquad+\frac{4\,\cosh\pi\kappa}k\sum_{n=1}^\infty t\!\left(-n\right)e\!\left(\frac{n\overline h}k\right)
\int\limits_a^bK_{2i\kappa}\!\left(\frac{4\pi\sqrt{nx}}k\right)f\!\left(x\right)\mathrm dx.
\end{align*}
\end{theorem}

The following upper bound for the $K$-Bessel function will be enough for estimating all the integrals involving it:
\[K_\nu(x)\ll_\nu x^{-1/2}\,e^{-x}\ll_Ax^{-A},\]
where $A>0$ is fixed and $x\gg1$.
This follows from (5.11.9) in \cite{Lebedev}. Here $\nu$ is fixed. In particular, we may estimate
\begin{equation}\label{K-asymptotics}
K_{2i\kappa}\!\left(\frac{4\pi\sqrt{nx}}k\right)
\ll_Ak^A\,n^{-A/2}\,x^{-A/2},
\end{equation}
for $x\in\left[1,\infty\right[$ and $n,k\in\mathbb Z_+$ satisfying $nx\gg k^2$.

For the  J-Bessel functions appearing in the Voronoi summation formula, we have the following asymptotics For every $K\in\mathbb Z_+$, we have the asymptotics, again for $n\,x\gg k^2$,
\begin{align}\label{Voronoi-J-asymptotics}
&J_{2i\kappa}\!\left(\frac{4\pi\sqrt{nx}}k\right)-J_{-2i\kappa}\!\left(\frac{4\pi\sqrt{nx}}k\right)\nonumber\\
&=\frac{k^{1/2}\,\sinh\pi\kappa}{\pi\,\sqrt2}\,n^{-1/4}\,x^{-1/4}\nonumber\\
&\qquad\cdot\sum_\pm (\pm1)e\!\left(\mp\frac18\pm\frac{2\sqrt{nx}}k\right)
\left(1+\sum_{\ell=1}^Kc_\ell^\pm\,k^\ell\,n^{-\ell/2}\,x^{-\ell/2}\right)\nonumber\\
&\qquad\qquad+O_K\!\left(k^{1/2+(K+1)}\,n^{-1/4-(K+1)/2}\,x^{-1/4-(K+1)/2}\right).
\end{align}
This follows form (5.11.6) of \cite{Lebedev}.

\section{Theorems on exponential integrals}\label{integral-section}

The use of the Voronoi summation formula leads to many exponential integrals, and so we will introduce several facts about such integrals.

Let us consider an interval $\left[M_1,M_2\right]\subseteq\mathbb R_+$, and let $U\in\mathbb R_+$ and $J\in\mathbb Z_+$ be such that $2JU<M_2-M_1$. Following \cite{Jutila1987a}, we introduce weight function $\eta_J$ by requiring that
\begin{align}\label{Jutilaweightfunction}
\int\limits_{M_1}^{M_2}\eta_J(x)\,h(x)\,\mathrm dx
=U^{-J}\int\limits_0^U\int\limits_0^U\cdots\int\limits_0^U
\int\limits_{M_1+u_1+\ldots+u_J}^{M_2-u_1-\ldots-u_J}h(x)\,\mathrm dx\,\mathrm du_J\cdots\mathrm du_2\,\mathrm du_1
\end{align}
for any integrable function $h$ on $\mathbb R$. It is not too difficult to see that actually $\eta_J$ is given by the convolution
\[\eta_J=\frac1U\chi_{\left[0,U\right]}\ast
\frac1U\chi_{\left[0,U\right]}\ast\ldots\ast
\frac1U\chi_{\left[0,U\right]}
\ast\chi_{\left[M_1,M_2-JU\right]},\]
with $U^{-1}\chi_{\left[0,U\right]}$ appearing $J$ times. In particular, $\eta_J$ is $J-1$ times continuously differentiable on $\mathbb R$, and supported in $\left[M_1,M_2\right]$.

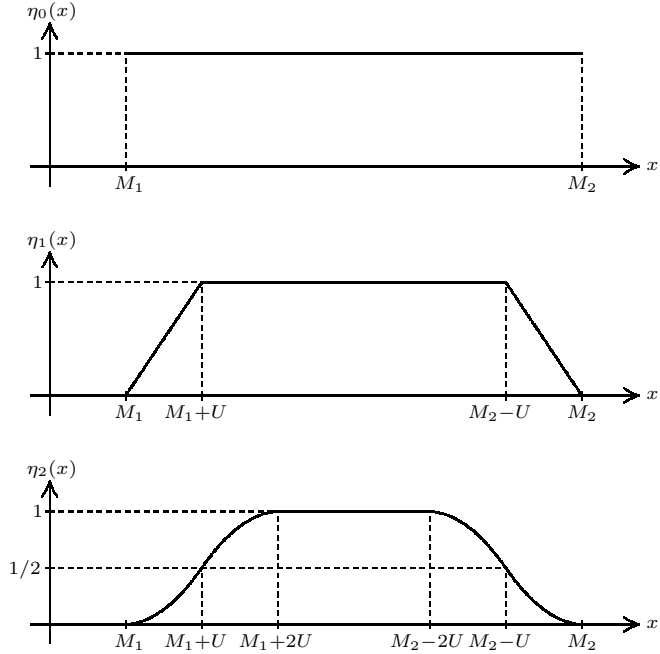
\begin{figure}[h]
\begin{center}
\setlength{\unitlength}{.5cm}
\begin{picture}(16,5.3)(-1,-.8)
\qbezier(-.5,0)(6.5,0)(15.5,0)
\qbezier(0,-.5)(0,1.75)(0,3.8)
\qbezier(2,-.1)(2,0)(2,.1)
\qbezier(14,-.1)(14,0)(14,.1)
\put(1.7,-.6){$\scriptstyle M_1$}
\put(13.6,-.6){$\scriptstyle{M_2}$}
\qbezier(15.5,0)(15.3,.05)(15.1,.2)
\qbezier(15.5,0)(15.3,-.05)(15.1,-.2)
\qbezier(0,3.8)(.05,3.6)(.2,3.4)
\qbezier(0,3.8)(-.05,3.6)(-.2,3.4)
\qbezier(-.1,3)(0,3)(.1,3)
\put(-.45,2.85){$\scriptstyle1$}
\put(15.7,-.1){$\scriptstyle x$}
\put(-.6,4){$\scriptstyle{\eta_0(x)}$}
\multiput(0,3)(.25,0){18}{\qbezier(0,0)(.05,0)(.1,0)}
\multiput(2,0)(0,.25){12}{\qbezier(0,0)(0,.05)(0,.1)}
\multiput(14,0)(0,.25){12}{\qbezier(0,0)(0,.05)(0,.1)}
\linethickness{.2mm}
\linethickness{.25mm}
\qbezier(2,3)(8,3)(14,3)
\end{picture}
\end{center}
\begin{center}
\setlength{\unitlength}{.5cm}
\begin{picture}(16,5.3)(-1,-.8)
\qbezier(-.5,0)(6.5,0)(15.5,0)
\qbezier(0,-.5)(0,1.75)(0,3.8)
\qbezier(2,-.1)(2,0)(2,.1)
\qbezier(4,-.1)(4,0)(4,.1)
\qbezier(12,-.1)(12,0)(12,.1)
\qbezier(14,-.1)(14,0)(14,.1)
\put(1.7,-.6){$\scriptstyle M_1$}
\put(3,-.6){$\scriptstyle{M_1+U}$}
\put(11,-.6){$\scriptstyle{M_2-U}$}
\put(13.6,-.6){$\scriptstyle{M_2}$}
\qbezier(15.5,0)(15.3,.05)(15.1,.2)
\qbezier(15.5,0)(15.3,-.05)(15.1,-.2)
\qbezier(0,3.8)(.05,3.6)(.2,3.4)
\qbezier(0,3.8)(-.05,3.6)(-.2,3.4)
\qbezier(-.1,3)(0,3)(.1,3)
\put(-.45,2.85){$\scriptstyle1$}
\put(15.7,-.1){$\scriptstyle x$}
\put(-.6,4){$\scriptstyle{\eta_1(x)}$}
\multiput(0,3)(.25,0){18}{\qbezier(0,0)(.05,0)(.1,0)}
\multiput(4,0)(0,.25){12}{\qbezier(0,0)(0,.05)(0,.1)}
\multiput(12,0)(0,.25){12}{\qbezier(0,0)(0,.05)(0,.1)}
\linethickness{.2mm}
\qbezier(2,0)(3,1.5)(4,3)
\qbezier(12,3)(13,1.5)(14,0)
\linethickness{.25mm}
\qbezier(4,3)(6,3)(12,3)
\end{picture}
\end{center}
\begin{center}
\setlength{\unitlength}{.5cm}
\begin{picture}(16,5.3)(-1,-.8)
\qbezier(-.5,0)(6.5,0)(15.5,0)
\qbezier(0,-.5)(0,1.75)(0,3.8)
\qbezier(2,-.1)(2,0)(2,.1)
\qbezier(4,-.1)(4,0)(4,.1)
\qbezier(6,-.1)(6,0)(6,.1)
\qbezier(10,-.1)(10,0)(10,.1)
\qbezier(12,-.1)(12,0)(12,.1)
\qbezier(14,-.1)(14,0)(14,.1)
\put(1.7,-.6){$\scriptstyle M_1$}
\put(3,-.6){$\scriptstyle{M_1+U}$}
\put(4.97,-.6){$\scriptstyle{M_1+2U}$}
\put(8.95,-.6){$\scriptstyle{M_2-2U}$}
\put(11,-.6){$\scriptstyle{M_2-U}$}
\put(13.6,-.6){$\scriptstyle{M_2}$}
\qbezier(15.5,0)(15.3,.05)(15.1,.2)
\qbezier(15.5,0)(15.3,-.05)(15.1,-.2)
\qbezier(0,3.8)(.05,3.6)(.2,3.4)
\qbezier(0,3.8)(-.05,3.6)(-.2,3.4)
\qbezier(-.1,3)(0,3)(.1,3)
\put(-.45,2.85){$\scriptstyle1$}
\qbezier(-.1,1.5)(0,1.5)(.1,1.5)
\put(-1.08,1.35){$\scriptstyle{1/2}$}
\put(15.7,-.1){$\scriptstyle x$}
\put(-.6,4){$\scriptstyle{\eta_2(x)}$}
\multiput(0,3)(.25,0){24}{\qbezier(0,0)(.05,0)(.1,0)}
\multiput(0,1.5)(.25,0){48}{\qbezier(0,0)(.05,0)(.1,0)}
\multiput(4,0)(0,.25){6}{\qbezier(0,0)(0,.05)(0,.1)}
\multiput(6,0)(0,.25){12}{\qbezier(0,0)(0,.05)(0,.1)}
\multiput(10,0)(0,.25){12}{\qbezier(0,0)(0,.05)(0,.1)}
\multiput(12,0)(0,.25){6}{\qbezier(0,0)(0,.05)(0,.1)}
\linethickness{.2mm}
\qbezier(2,0)(3,.05)(4,1.5)
\qbezier(4,1.5)(5,2.95)(6,3)
\qbezier(10,3)(11,2.95)(12,1.5)
\qbezier(12,1.5)(13,.05)(14,0)
\linethickness{.25mm}
\qbezier(6,3)(8,3)(10,3)
\end{picture}
\end{center}
\caption{\label{weight-functions-of-jutila}
A sketch of the weight functions $\eta_0$, $\eta_1$ and $\eta_2$. Please note that for $J\geqslant1$ the weight function $\eta_J$ is $C^\infty$-smooth except for the points $M_1+\ell U$ and $M_2-\ell U$, where $\ell\in\left\{0,1,\ldots,J\right\}$, where it belongs only to $C^{J-1}$.}
\end{figure}

The following saddle point theorem is a special case of Theorem 2.2 in \cite{Jutila1987a}.
\begin{theorem}\label{saddle-point-lemma}
Let us consider an interval $\left[M_1,M_2\right]\subseteq\mathbb R_+$, let $\mu\in\mathbb R_+$, and let $D$ stand for the domain
\[D=\bigl\{z\in\mathbb C\mathbin:\text{$\left|z-x\right|<\mu$ for some $x\in\left[M_1,M_2\right]$}\bigr\}.\]
Let $f,g\colon D\longrightarrow\mathbb C$ be holomorphic, let $F,G\in\mathbb R_+$, and assume that
\[f(x)\in\mathbb R,\quad f''(x)>0\quad\text{and}\quad f''(x)\gg F\,\mu^{-2},\]
for $x\in\left[M_1,M_2\right]$, and that
\[f'(z)\ll F\,\mu^{-1}\quad\text{and}\quad g(z)\ll G\]
for $z\in D$.

Next, let $U\in\mathbb R_+$ and $J\in\mathbb Z_+$ be such that $2\,J\,U<M_2-M_1$, and let $\eta_J$ denote the weight function defined as above, namely the convolution
\[\eta_J=\frac1U\chi_{\left[0,U\right]}\ast
\frac1U\chi_{\left[0,U\right]}\ast\ldots\ast
\frac1U\chi_{\left[0,U\right]}
\ast\chi_{\left[M_1,M_2-JU\right]},\]
with $U^{-1}\chi_{\left[0,U\right]}$ appearing $J$ times.

Finally, let $\alpha\in\mathbb R$, and let $x_0\in\left]M_1,M_2\right[$ be such that $f'(x_0)+\alpha=0$.
Then
\begin{multline*}
\int\limits_{M_1}^{M_2}g(x)\,e(f(x)+\alpha x)\,\eta_J(x)\,\mathrm dx\\
=\xi_J(x_0)\,g(x_0)\,f''(x_0)^{-1/2}\,e\bigl(f(x_0)+\alpha x_0+1/8\bigr)+\text{error},
\end{multline*}
where the error is
\begin{align*}
&\ll\left(M_2-M_1\right)\left(1+\mu^J\,U^{-J}\right)\,G\,e^{-A\left|\alpha\right|\mu-A\,F}\\
&\qquad+\left(1+\chi(x_0)\,F^{1/2}\right)G\,\mu\,F^{-3/2}\\
&\qquad+U^{-J}\sum_{j=0}^J\left(E_J(M_1+jU)+E_J(M_2-jU)\right).
\end{align*}
Here $A$ is some positive real constant independent of $f$, $g$, $\alpha$, and $\left[M_1,M_2\right]$, $\chi$ denotes the characteristic function of the set $\left]M_1,M_1+JU\right[\cup\left]M_2-JU,M_2\right[$, the symbol $E_J(x)$ stands for
\[E_J(x)=\frac G{\left(\left|f'(x)+\alpha\right|+f''(x)^{1/2}\right)^{J+1}},\]
and the factor $\xi_J(x_0)$ is as follows:
\begin{enumerate}\setlength{\itemsep}{0pt}
\item\label{property-one}
If $M_1+JU<x_0<M_2-JU$, then
\[\xi_J(x_0)=1.\]
\item\label{property-two}
If $M_1<x_0\leqslant M_1+JU$, then
\[\xi_J(x_0)=(J!\,U^J)^{-1}\sum_{j=0}^{j_1}\binom Jj(-1)^j\sum_{0\leqslant\nu\leqslant J/2}c_\nu\,f''(x_0)^{-\nu}\left(x_0-M_1-jU\right)^{J-2\nu},\]
where $j_1$ is the largest integer with $M_1+j_1U<x_0$.
\item\label{property-three}
If $M_2-JU\leqslant x_0<M_2$, then
\[\xi_J(x_0)=(J!\,U^J)^{-1}\sum_{j=0}^{j_2}\binom Jj(-1)^j\sum_{0\leqslant\nu\leqslant J/2}c_\nu\,f''(x_0)^{-\nu}\left(M_2-jU-x_0\right)^{J-2\nu},\]
where $j_2$ is the largest integer with $M_2-j_2U>x_0$.
\end{enumerate}
The coefficients $c_\nu$ are fixed numerical constants only depending on $J$. Furthermore, $\xi_J(x)$ is continuously differentiable in the intervals $\left]M_1,M_1+JU\right[$ and $\left]M_2-JU,M_2\right[$ except for the points $M_1+jU$ and $M_2-jU$ appearing in the above sums, and the derivative satisfies in these intervals, where it exists, the estimate $\xi_J'(x)\ll U^{-1}$.
\end{theorem}
\noindent
Strictly speaking, the last statement about $\xi_J'(x)$ does not appear in the statement of Theorem 2.2 in \cite{Jutila1987a}, but it follows easily by inspecting the above sums for $\xi_J(x_0)$.

Some of the exponential integrals we will meet will not have saddle points. They can be handled with the following theorem, which is a special case of Theorem 2.3 in \cite{Jutila1987a}.
\begin{theorem}\label{saddle-point-lemma-without-saddle-points}
Let us consider an interval $\left[M_1,M_2\right]\subseteq\mathbb R_+$, let $\mu\in\mathbb R_+$, and let $D$ stand for the domain
\[D=\bigl\{z\in\mathbb C\mathbin:\text{$\left|z-x\right|<\mu$ for some $x\in\left[M_1,M_2\right]$}\bigr\}.\]
Let $f,g\colon D\longrightarrow\mathbb C$ be holomorphic, let $F,G\in\mathbb R_+$, and assume that
\[f(x)\in\mathbb R,\quad\text{and}\quad f'(x)\asymp F\,\mu^{-1},\]
for $x\in\left[M_1,M_2\right]$, and that
\[f'(z)\ll F\,\mu^{-1}\quad\text{and}\quad g(z)\ll G\]
for $z\in D$.

Next, let $U\in\mathbb R_+$ and $J\in\mathbb Z_+$ be such that $2\,J\,U<M_2-M_1$, and let $\eta_J$ denote the weight function defined as above, namely the convolution
\[\eta_J=\frac1U\chi_{\left[0,U\right]}\ast
\frac1U\chi_{\left[0,U\right]}\ast\ldots\ast
\frac1U\chi_{\left[0,U\right]}
\ast\chi_{\left[M_1,M_2-JU\right]},\]
with $U^{-1}\chi_{\left[0,U\right]}$ appearing $J$ times.

Finally, let $\alpha\in\mathbb R$. Then
\begin{multline*}
\int\limits_{M_1}^{M_2}g(x)\,e(f(x)+\alpha x)\,\eta_J(x)\,\mathrm dx\\
\ll U^{-J}\,G\,\mu^{J+1}\,F^{-J-1}+\left(\mu^J\,U^{1-J}+M_2-M_1\right)G\,e^{-A\,F}.
\end{multline*}
Here $A$ is some positive real constant independent of $f$, $g$, $\alpha$, and $\left[M_1,M_2\right]$.
\end{theorem}

We will also use the following lemma for estimating exponential integrals. It is Lemma 6 in \cite{Jutila--Motohashi}.
\begin{lemma}\label{jutila-motohashi-lemma}
Let $M_1,M_2\in\mathbb R_+$ and $M_1<M_2$, let $J\in\mathbb Z_+$, and let $g\in C_{\mathrm c}^J(\mathbb R_+)$ with $\mathrm{supp}\,g\subseteq\left[a,b\right]$, and let $G_0$ and $G_1$ be such that
\[g^{(\nu)}(x)\ll_\nu G_0\,G_1^{-\nu}\]
for all $x\in\mathbb R_+$ for each $\nu\in\left\{0,1,2,\ldots,J\right\}$. Also, let $f$ be a holomorphic function defined in $D\subseteq\mathbb C$, which consists all points in the complex plane with distance smaller than $\mu\in\mathbb R_+$ from the interval $\left[M_1,M_2\right]$ of the real axis. Assume that $f$ is real-valued on $\left[M_1,M_2\right]$ and let $F_1\in\mathbb R_+$ be such that
\[F_1\ll\left|f'(z)\right|\]
for all $z\in D$. Then, for all all $P\in\left\{1,2,\ldots,J\right\}$,
\[\int\limits_{M_1}^{M_2}g(x)\,e(f(x))\,\mathrm dx\ll_PG_0\left(G_1\,F_1\right)^{-P}\left(1+\frac{G_1}\mu\right)^P\left(M_2-M_1\right).\]
\end{lemma}

Finally, for completeness, we state the following two classical tools, known as the first derivative test and the second derivative test, respectively. For a discussion of these, see e.g.\ Section 5.1 in \cite{Huxley1996}.

\begin{lemma}\label{first-derivative-test}
Let $M_1,M_2\in\mathbb R$ with $M_1<M_2$, let $\lambda\in\mathbb R_+$, and let $f$ be a real-valued continuously differentiable function on $\left]M_1,M_2\right[$ such that $\left|f'(x)\right|\geqslant\lambda$ for $x\in\left]M_1,M_2\right[$. Also, let $g$ be a complex-valued continuously differentiable function on the interval $\left[M_1,M_2\right]$, and let $G\in\mathbb R_+$ be such that
$g(x)\ll G$ for $x\in\left[M_1,M_2\right]$.
Then
\[\int\limits_{M_1}^{M_2}g(x)\,e(f(x))\,\mathrm dx\ll\frac G\lambda+\frac1\lambda\int\limits_{M_1}^{M_2}\left|g'(x)\right|\mathrm dx.\]
\end{lemma}

\begin{lemma}\label{second-derivative-test}
Let $M_1,M_2\in\mathbb R$ with $M_1<M_2$, let $\lambda\in\mathbb R_+$, and let $f$ be a real-valued twice continuously differentiable function on $\left]M_1,M_2\right[$ such that $\left|f''(x)\right|\geqslant\lambda$ for $x\in\left]M_1,M_2\right[$. Also, let $g$ be a complex-valued continuously differentiable function on the interval $\left[M_1,M_2\right]$, and let $G\in\mathbb R_+$ be such that $g(x)\ll G$ for $x\in\left[M_1,M_2\right]$.
Then
\[\int\limits_{M_1}^{M_2}g(x)\,e(f(x))\,\mathrm dx\ll\frac G{\sqrt\lambda}+\frac1{\sqrt\lambda}\int\limits_{M_1}^{M_2}\left|g'(x)\right|\mathrm dx.\]
\end{lemma}

\section{A transformation formula for smoothed exponential sums}\label{transformation-section}

\subsection{Statement of the transformation formula}

In the following theorem $\delta_1$, $\delta_2$, \dots denote positive constants which may be supposed to be arbitrarily small. Further, we write $L$ for $\log M_1$.

\begin{theorem}\label{general-transformation-formula-for-cusp-forms}
Let $2\leqslant M_1<M_2\leqslant2M_1$. We assume that $M_1$ is sufficiently large, the notion of sufficiently large depending on the implicit constants in the assumptions below and on $\delta_1$.
Let $f$ and $g$ be holomorphic functions in the domain
\[D=\bigl\{z\in\mathbb C\mathbin:\text{$\left|z-x\right|<c\,M_1$ for some $x\in\left[M_1,M_2\right]$}\bigr\},\]
where $c$ is a positive constant. Suppose that $f(x)$ is real for $x\in\left[M_1,M_2\right]$. Suppose also that, for some positive numbers $F$ and $G$,
\[g(z)\ll G,\qquad f'(z)\ll F\,M_1^{-1}\]
for $z\in D$, and that
\[f''(x)>0\quad\text{and}\quad f''(x)\gg F\,M_1^{-2}\quad\text{for}\quad x\in\left[M_1,M_2\right].\]
Let $r=h/k$ be a rational number such that $\left(h,k\right)=1$,
\[1\leqslant k\ll M_1^{1/2-\delta_1},\]
\[r\asymp F\,M_1^{-1}\]
and
\[f'(M(r))=r\]
for a certain number $M(r)\in\left]M_1,M_2\right[$.
Write
\[M_j=M(r)+(-1)^j\,m_j,\quad j=1,2.\]
Suppose that $m_1\asymp m_2$, and that
\[M_1^{\delta_2}\max\left\{M_1\,F^{-1/2},\left|hk\right|\right\}\ll m_1\ll M_1^{1-\delta_3}.\]
Define for $j\in\left\{1,2\right\}$
\[p_{j,n}(x)=f(x)-rx+(-1)^{j-1}\left(\frac{2\sqrt{nx}}k-\frac18\right),\]
\[n_j=(r-f'(M_j))^2\,k^2\,M_j,\]
and for $n<n_j$ let $x_{j,n}$ be the (unique) zero of $p_{j,n}'(x)$ in the interval $\left]M_1,M_2\right[$.
Also, let $J$ be a fixed positive integer and sufficiently large depending on $\delta_2$ and $\delta_4$.
Let
\[U\gg F^{-1/2}\,M_1^{1+\delta_4}\asymp F^{1/2}\,r^{-1}\,M_1^{\delta_4},\]
where $\delta_4>\delta_2$, and assume also that
\[JU<\frac{M_2-M_1}2.\]
Write for $j\in\left\{1,2\right\}$
\[M_j'=M_j+(-1)^{j-1}\,J\,U=M(r)+(-1)^j\,m_j',\]
and suppose that $m_j'\asymp m_j$. Define
\[n_j'=(r-f'(M_j'))^2\,k^2\,M_j'.\]
Then we have
\begin{align}\label{transformation-formula}
&\sum_{M_1\leqslant m\leqslant M_2}\eta_J(m)\,t(m)\,g(m)\,e(f(m))\\
&=i\,2^{-1/2}\,k^{-1/2}\sum_{j=1}^2(-1)^{j-1}\sum_{n<n_j}w_j(n)\,t(n)\,e\!\left(\frac{-n\overline h}k\right)
n^{-1/4}\,x_{j,n}^{-1/4}\nonumber\\
&\qquad\cdot g(x_{j.n})\,p_{j,n}''(x_{j,n})^{-1/2}\,e\!\left(p_{j,n}(x_{j,n})+\frac18\right)\nonumber\\
&\qquad\qquad+O\!\left(k^{-1/2}\,m_1^{1/2}\,F^{-1}\,G\,\left|h\right|^{3/2}\,U\,L\,F^{1/4}\,M_1^{\vartheta}\right),\nonumber
\end{align}
where
\[w_j(n)=1\quad\text{for}\quad n<n_j',\]
\[w_j(n)\ll1\quad\text{for}\quad n<n_j,\]
$w_j(y)$ and $w_j'(y)$ are piecewise continuous functions in the interval $\left]n_j',n_j\right[$ with at most $J-1$ discontinuities, and
\[w_j'(y)\ll\left(n_j-n_j'\right)^{-1}\quad\text{for}\quad n_j'<y<n_j\]
whenever $w_j'(y)$ exists.
\end{theorem}

\subsection{The proof}

A word on the notation: In the following $j\in\left\{1,2\right\}$. This parameter comes about as follows: After applying the Voronoi summation formula and replacing the $J$-Bessel function by a simpler asymptotic expression, the cosine is replaced by the sum of two exponentials with phase factors of opposite signs. The value $j=1$ corresponds to the $+$-sign, and $j=2$ corresponds to the $-$-sign. For simplicity, we consider the various errors with fixed $j$; i.e. we omit the summation symbol $\sum_{j=1}^2$.

\subsubsection{Sizes of the parameters}

Suppose, to be specific, that $r>0$, and thus $h>0$. The proof is similar for $r<0$.

The assertion \eqref{transformation-formula} should be understood as an asymptotic result, in which $M_1$ and $M_2$ are large.

We observe that since,
\[h\,k\,M_1^{\delta_2}\ll M_1^{1-\delta_3},\]
we have $h\ll M_1^{1-\delta_2-\delta_3}$.

\paragraph{On the size of $F$.}
The number $F$ wil be large. In fact,
\[F\gg M_1\,r\geqslant k^{-1}\,M_1\gg M_1^{1/2+\delta_1}.\]
Before proving the latter, we observe that, in $\left[M_1,M_2\right]$,
\[f''(x)=\int\limits_{\partial B(x,cM_1/2)}\frac{f'(z)\,\mathrm dz}{(z-x)^2}\ll M_1\cdot\frac{F\,M_1^{-1}}{M_1^2}=F\,M_1^{-2},\]
so that $f''(x)\asymp F\,M_1^{-2}$. Here $c$ is the positive constant from the definition of~$D$. The same argument shows in fact more: we have $f''(z)\ll F\,M_1^{-2}$ for $z$ in, say, $D(M_1,M_2,cM_1/4)$.

We should also point out that $F$ is not very large either: Since
\[hk\ll M_1^{1-\delta_3-\delta_2}\quad\text{and}\quad\frac F{M_1}\asymp\frac hk,\]
we have
\[F\ll\frac{M_1\,h}k\ll\frac{M_1^{2-\delta_3-\delta_2}}{k^2}\ll M_1^{2-\delta_3-\delta_2}.\]
Very crudely, but more simply, $F$ is bounded from above and below by powers of $M_1$:
\[M_1^{1/2}\ll F\ll M_1^2.\]

\paragraph{On the sizes of $n_j$ and $n_j'$.}
The number $n_j$ will also be large:
\[n_j\gg h\,k\,M_1^{2\delta_2}.\]
Now
\[\frac hk-f'(M_j)\asymp\int\limits_{M_j}^{M(r)}f''(x)\,\mathrm dx\asymp m_j\,F\,M_1^{-2},\]
so that
\begin{align*}
n_j&=\left(\frac hk-f'(M_j)\right)^2k^2\,M_1
\asymp m_j^2\,F^2\,M_1^{-4}\,k^2\,M_1
\asymp F^{-1}\,h^3\,k^{-1}\,m_j^2\\
&\gg F^{-1}\,h^3\,k^{-1}\,M_1^{2+2\delta_2}\,F^{-1}
=F^{-2}\,h^3\,k^{-1}\,M_1^{2+2\delta_2}\\
&\gg k^2\,h^{-2}\,M_1^{-2}\,h^3\,k^{-1}\,M_1^{2+2\delta_2}
=h\,k\,M_1^{2\delta_2}.
\end{align*}

By using the estimate derived in the previous calculation we get
\begin{align*}
n_j&\asymp F^{-1}\,h^3\,k^{-1}\,m_j^2\asymp k\,h^{-1}\,M_1^{-1}\,h^3\,k^{-1}\,m_j^2\\
&\ll M_1^{-1}\,h^2\,m_j^2\ll M_1^{3-2\delta_2}, 
\end{align*}
since $h\ll M_1$. In particular,
\[\log n_j\ll\log M_1.\]
We also have a simple estimate
\begin{align*}
n_j\asymp F^{-1}h^3k^{-1}m_j^2\asymp F^{-1}h^3k^{-1}m_j'^2\asymp n_j'
\end{align*}
due to the fact $m_j\asymp m_j'$. 

We will also need to know the size of $n_j-n_j'$. To this end, let us write
\begin{align*}
n_j-n_j'&=\left(r-f'(M_j)\right)^2\,k^2\,M_j-(r-f'(M_j'))^2\,k^2\,M_j'\\
&=k^2\int\limits_{M_j'}^{M_j}\frac{\mathrm d}{\mathrm dt}\left(\left(r-f'(t)\right)^2t\right)\mathrm dt\\
&=k^2\int\limits_{M_j'}^{M_j}\left(\left(r-f'(t)\right)^2-2t\left(r-f'(t)\right)f''(t)\right)\mathrm dt.
\end{align*}
In the integrand the first term is
\[\asymp\left(\int\limits_t^{M(r)}f''(t)\,\mathrm dt\right)^{\!\!2}
\asymp\left(m_j\,F\,M_1^{-2}\right)^2\asymp m_j^2\,F^2\,M_1^{-4},\]
and the second term is similarly
\[\asymp M_1\,m_j\,F^2\,M_1^{-4}\asymp m_j\,F^2\,M_1^{-3}.\]
The first term is smaller than the second by the extra factor $m_j\,M_1^{-1}\ll M_1^{-\delta_3}$, and so the second term dominates in the integral, and we get
\[n_j-n_j'\asymp k^2\,U\,m_j\,F^2\,M_1^{-3}.\]

\subsubsection{The behaviour of $p_{j,n}$}

After the application of Voronoi's summation formula and replacing the $J$-Bessel function by its asymptotics, the phase functions in the individual integrals will be given by
\[p_{j,n}(x)=f(x)-\frac{hx}k+(-1)^{j-1}\left(\frac{2\sqrt{nx}}k+\frac18\right),\]
where $x$ ranges over $\left[M_1,M_2\right]$.

The parameter $n_j$ (which, despite the notation, is not necessarily an integer) is chosen so that
\[p_{j,n_j}'(M_j)=0.\]
As the derivative is
\[p_{j,n}'(x)=f'(x)-\frac hk+(-1)^{j-1}\frac{\sqrt n}{k\sqrt x},\]
this simplifies to
\[n_j=\left(\frac hk-f'(M_j)\right)^2k^2\,M_j.\]

Now the first salient feature of the function $p_{j,n}$ is that $p_{j,n}'$ has a unique zero $x_{j,n}$ in the interval $\left]M_1,M_2\right[$ for $n<n_j$.
The second feature is that $p_{j,n}'(x)$ has no zero in $\left]M_1,M_2\right[$ when $n\geqslant n_j$.

The existence of a zero in $\left]M_1,M_2\right[$ when $n<n_j$ is easily seen from the inequalities
\[(-1)^j\,p_{j,n}'(M(r))<0\]
and
\[(-1)^j\,p_{j,n}'(M_j)>(-1)^j\,p_{j,n_j}'(M_j)=0,\]
where the latter follows from the fact that $p_{j,n}'$ behaves monotonically with respect to~$n$.
Furthermore, the zero $x_{j,n}$, whose existence is guaranteed when $n<n_j$, lies on $\left]M_1,M(r)\right[$ when $j=1$, and on $\left]M(r),M_2\right[$ when $j=2$.

When $j=2$, the derivative $p_{2,n}'(x)$ is monotonically increasing and therefore it is clear that $x_{2,n}$ is unique for $n<n_2$, and that there is no zero when $n\geqslant n_j$ as \[p_{2,n}'(M_2)<p_{2,n_2}'(M_2)=0.\]
We mention in passing that, in fact,
\[p_{2,n}''(x)=f''(x)+\frac{\sqrt{n}}{2k\,x^{3/2}}\asymp F\,M_1^{-2}\]
on the interval $\left[M_1,M_2\right]$.

The case $j=1$ is slightly less obvious. The main point is that by inspecting $p_{1,n}''(x)$, we will see that $p_{1,n}'$ is strictly increasing in $\left[M_1,M_2\right]$ when $n\leqslant n_j$, which will guarantee the uniqueness of $x_{1,n}$ and the non-existence of zero $n=n_1$ (if $n_1$ happens to be an integer) as
\[p_{1,n_1}'(M_1)=0.\]
In particular, $p_{1,n_1}'$ takes only non-negative values in $\left[M_1,M_2\right]$ and we get, for $n>n_1$, that
\[p_{1,n}'(x)>p_{1,n_1}'(x)\geqslant0,\]
thereby excluding the possibility of zeros.

Now it only remains to show that $p_{1,n}''(x)$ is positive for $n\leqslant n_1$, when $M_1$ is supposed to be sufficiently large.
Since
\[p_{1,n}''(x)=f''(x)-\frac12\,k^{-1}\,n^{1/2}\,x^{-3/2},\]
and
\begin{align*}
&\frac12\,k^{-1}\,n^{1/2}\,x^{-3/2}\ll k^{-1}\,n_1^{1/2}\,M_1^{-3/2}\\
&\asymp k^{-1}\,m_1\,F\,M_1^{-2}\,k\,M_1^{1/2}\,M_1^{-3/2}
=m_1\,F\,M_1^{-3},
\end{align*}
as well as $m_1\ll M^{1-\delta_1}$,
we indeed have $p_{1,n}''(x)\asymp f''(x)\asymp F\,M_1^{-2}$ if only $M_1$ is sufficiently large, depending (at most) on the implicit constants in the assumptions of the theorem and $\delta_1$.

The reason for introducing the numbers $n_j'$ is the following: when $n<n_j'$, the corresponding saddle-points $x_{j,n}$ lie on the interval $\left]M_1,M_2\right[$. This is not hard to see: for $j=1$ the saddle-point $x_{1,n}$ decreases strictly monotonically as $n$ increases, and the value $n_1'$ corresponds to the situation where $x_{1,n}$ lies precisely at $M_1$. For $j=2$ things work similarly, except that $x_{2,n}$ increase monotonically as $n$ increases. The monotonicity of $x_{j,n}$ with respect to $n$ follows from the fact that the expression for $p_{j,n}'(x)$ depends strictly monotonically on~$n$.

\subsubsection{The derivative of the saddle-point}

Let us define more generally for $y\in\left]n_j',n_j\right[$ the saddle-point $x_{j,y}$ so that it is the unique zero of
\[p_{j,y}'(x)=f'(x)-r+(-1)^{j-1}\frac{\sqrt y}{k\sqrt x},\]
so that
\[f'(x_{j,y})-r+(-1)^{j-1}\frac{\sqrt y}{k\sqrt{x_{j,y}}}=0.\]
Differentiating this with respect to $y$ gives first
\[\frac{\mathrm dx_{j,y}}{\mathrm dy}\,f''(x_{j,y})+(-1)^{j-1}\frac1{2\sqrt y\,k\sqrt{x_{j,y}}}+(-1)^{j-1}\frac{\sqrt y}k\left(-\frac12\right)x_{j,y}^{-3/2}\,\frac{\mathrm dx_{j,y}}{\mathrm dy}=0,\]
which simplifies to
\[\frac{\mathrm dx_{j,y}}{\mathrm dy}\,p''(x_{j,y})=\frac{(-1)^j}{2\sqrt y\,k\,\sqrt{x_{j,y}}}.\]
Since $p''(x_{j,y})\asymp F\,M_1^{-2}$ and $y\asymp n_j$, we have
\[\frac{\mathrm dx_{n,y}}{\mathrm dy}\asymp F^{-1}\,M_1^2\,k^{-1}\,M_1^{-1/2}\,n_1^{-1/2}\asymp\frac{m_1}{n_1}.\]

\subsubsection{Voronoi summation and Bessel asymptotics}

We begin the transformation of the exponential sum by applying the Voronoi-type summation formula for Maass forms:      
\begin{align*}
&\sum_{M_1\leqslant m\leqslant M_2}\eta_J(m)\,t(m)\,g(m)\,e(f(m))\\
&=\sum_{M_1\leqslant m\leqslant M_2}\eta_J(m)\,t(m)\,g(m)\,e\!\left(\frac{mh}k\right)e\!\left(f(m)-\frac{mh}k\right)\\
&=\frac{\pi i}{k \sinh \pi\kappa}\,\sum_{n=1}^\infty t(n)\,e\!\left(\frac{-n\overline h}k\right)\\
&\qquad\cdot\int\limits_{M_1}^{M_2}\left(J_{2i\kappa}\!\left(\frac{4\pi\sqrt{nx}}k\right)-J_{-2i\kappa}\!\left(\frac{4\pi\sqrt{nx}}k\right)\right)
\eta_J(x)\,g(x)\,e\!\left(f(x)-\frac{hx}k\right)\mathrm dx\\
&\quad+\frac{4\cosh \pi\kappa}k\,\sum_{n=1}^\infty t(-n)\,e\!\left(\frac{n\overline h}k\right)\\
&\qquad\cdot\int\limits_{M_1}^{M_2}K_{2i\kappa}\!\left(\frac{4\pi\sqrt{nx}}k\right)\eta_J(x)\,g(x)\,e\!\left(f(x)-\frac{hx}k\right)\mathrm dx.
\end{align*}
Using the asymptotics of $J$- and $K$-Bessel functions we combine above calculations to
\begin{align}\label{voronoi-sum}
&\sum_{M_1\leqslant m\leqslant M_2}\eta_J(m)\,t(m)\,g(m)\,e(f(m))\\
&=i\,2^{-1/2}\,k^{-1/2}\sum_{j=1}^2(-1)^{j-1}\sum_{n=1}^\infty t(n)\,n^{-1/4}\,e\!\left(\frac{-n\overline h}k\right)\nonumber\\
&\qquad\cdot\int\limits_{M_1}^{M_2}\,x^{-1/4}\,e\!\left(p_{j,n}(x)\right)\left(1+\sum_{\ell=1}^Kc_\ell^{(j)}k^\ell\,n^{-\ell/2}\,x^{-\ell/2}\right)\,\eta_J(x)\,g(x)\,\mathrm dx\nonumber\\
&+O\Bigg(\frac1k\sum_{n=1}^\infty \left|t(n)\right|k^{1/2+K+1}\,n^{-1/4-(K+1)/2}
\int\limits_{M_1}^{M_2}x^{-1/4-(K+1)/2}\eta_J(x)\,\left|g(x)\right|\,\mathrm dx\Bigg)\nonumber\\
&+O\Bigg(\frac1k\sum_{n=1}^{\infty}\left|t(-n)\right|\int\limits_{M_1}^{M_2}k^A\,n^{-A/2}\,x^{-A/2}\,\eta_J(x)\left|g(x)\right|\mathrm dx\Bigg)\nonumber
\end{align}
for any $K\in\mathbb Z_+$ and fixed $A>0$. Fixing large enough $K$ depending on the Maass form in question, the first error term on the right-hand side can be absorbed in the error term on the right-hand side of (\ref{transformation-formula}). Also, clearly the second error term is negligible in view of the error term by choosing large enough $A$.

Next, we will estimate the integral
\begin{align}\label{int1}
\int\limits_{M_1}^{M_2}\,x^{-1/4}\,e\!\left(p_{j,n}(x)\right)\left(1+\sum_{\ell=1}^Kc_\ell^{(j)}k^\ell n^{-\ell/2}x^{-\ell/2}\right)\,\eta_J(x)\,g(x)\,\mathrm dx.
\end{align}

\subsubsection{Large frequencies}

When $n>2n_j$, the integrals are estimated using Theorem \ref{saddle-point-lemma-without-saddle-points} with
\[\mu\asymp m_j,\qquad F\,\mu^{-1}:=k^{-1}\,M_1^{-1/2}\,n^{1/2}\gg m_j\,F\,M_1^{-2},\]
and
\[G:=M_1^{-1/4}\,G.\]
Of the conditions of the theorem, only the ones related to the size of $p_{j,n}'(x)$ are not immediately checked. Also, the parameter $\mu$ is $\asymp m_1$ instead of, say, $\asymp M_1$, in order for $p_{j,n}'(x)$ to be satisfy these conditions.

Since $f''(z)\ll F\,M_1^{-2}$ in $D(M_1,M_2,cM_1/2)$, we have
\[f'(z)-\frac hk\asymp\int\limits_{M(r)}^zf''(w)\,\mathrm dw\ll m_1\,F\,M_1^{-2}\]
for $z\in D(M_1,M_2,\mu)$, where it is best to integrate along the straight line segment connecting $M(r)$ and $z$. Thus we have
\begin{multline*}
p_{j,n}'(z)=f'(z)-\frac hk+(-1)^{j-1}\frac{\sqrt n}{k\sqrt z}\\\ll
m_1\,F\,M_1^{-2}+k^{-1}\,M_1^{-1/2}\,n^{1/2}\ll k^{-1}\,M_1^{-1/2}\,n^{1/2}.
\end{multline*}

The conclusion that $p_{j,n}'(x)\asymp M$ on the interval $\left[M_1,M_2\right]$, when $n>2n_j$, can be obtained by comparing $p_{j,n}'(x)$ with $p_{j,n_j}'(x)$. More precisely, when $j=1$, the function $p_{j,n_j}'(x)$ is non-negative, bounded from above by $\ll M$ by estimates similar to the ones above, and the difference $p_{j,n}'(x)-p_{j,n_j}'(x)$ is
\begin{align*}
&=\frac{\sqrt{n\vphantom{kn_j}}-\sqrt{n_j\vphantom{kn}}}{k\sqrt x}\asymp\frac{\sqrt n}{k\sqrt{M_1}}.
\end{align*}

When $j=2$, the conclusion is obtained in the same way, except that now $p_{j,n_j}'(x)$ is non-positive, and the difference $p_{j,n}'(x)-p_{j,n_j}'(x)$ has the opposite sign.

Now that the assumptions of the theorem certainly hold, the estimate will be
\begin{align*}
&\int\limits_{M_1}^{M_2}\,x^{-1/4}\,e\!\left(p_{j,n}(x)\right)\left(1+\sum_{\ell=1}^Kc_\ell^{(j)}k^\ell n^{-\ell/2}x^{-\ell/2}\right)\,\eta_J(x)\,g(x)\,\mathrm dx\\
&\qquad\ll U^{-J}\,M_1^{-1/4}\,G\,\left(k^{-1}\,M_1^{-1/2}\,n^{1/2}\right)^{-J-1}\\
&\qquad\qquad+\left(m_1^J\,U^{1-J}+m_1\right)M_1^{-1/4}\,G\,\exp\!\left(-A\left(k^{-1}M_1^{-1/2}n^{1/2}\right)m_1\right)\\
&\qquad\ll k^{J+1}\,U^{-J}\,G\,M_1^{J/2+1/4}\,n^{-J/2-1/2}\\
&\qquad\qquad+\left(m_1^J\,U^{1-J}+m_1\right)M_1^{-1/4}\,G\,\exp\!\left(-A\left(k^{-1}M_1^{-1/2}n^{1/2}\right)m_1\right)
\end{align*}

The terms with $n>2n_j$ contribute
\begin{align*}
&\ll k^{-1/2}\,\sum_{n>2n_j}\left|t(n)\right|n^{-1/4}
\Bigg(k^{J+1}\,U^{-J}\,G\,M_1^{J/2+1/4}\,n^{-J/2-1/2}\\
&\qquad\qquad+\left(m_1^J\,U^{1-J}+m_1\right)M_1^{-1/4}\,G\,\exp\!\left(-A\left(k^{-1}M_1^{-1/2}n^{1/2}\right)m_1\right)\Bigg).
\end{align*}

\paragraph{The first error term.}
The error from the first error term is
\begin{align*}
&\ll k^{J+1/2}\,U^{-J}\,G\,M_1^{J/2+1/4}\sum_{n>2n_j}\left|t(n)\right|n^{-J/2-3/4}\\
&\ll k^{J+1/2}\,U^{-J}\,G\,M_1^{J/2+1/4}\,n_j^{-J/2+1/4}\\
&\ll k^{J+1/2}\,U^{-J}\,G\,M_1^{J/2+1/4}\left(m_1^2\,F^2\,M_1^{-4}\,k^2\,M_1\right)^{-J/2+1/4}\\
&\ll k\,U^{-J}\,G\,F^{-J+1/2}\,M_1^{2J-1/2}\,m_1^{-J+1/2}\\
&\ll k\,F^{J/2}\,M_1^{-\delta_4J-J}\,G\,F^{-J+1/2}\,M_1^{2J-1/2}\left(F^{-1/2}\,M_1^{1+\delta_2}\right)^{-J+1/2}\\
&\ll k\,F^{J/2}\,M_1^{-\delta_4J-J}\,G\,F^{-J+1/2}\,M_1^{2J-1/2}\left(F^{-1/2}\,M_1^{1+\delta_2}\right)^{-J+1/2}\\
&\ll k\,G\,F^{1/2}\,M_1^{-\delta_4J-\delta_2J-\delta_2/2},
\end{align*}
and this is, provided that $J$ is sufficiently large, depending on $\delta_2$ or $\delta_4$,
\[\ll F^{-1}\,G\,h^{3/2}\,k^{-1/2}\,m_1^{1/2}\,U\,L,\]
which is small enough.

\paragraph{The second error term.}
Since
\begin{align*}
k^{-1}\,M_1^{-1/2}\,n^{1/2}\,m_1
&\gg k^{-1}\,M_1^{-1/2}\,\left(m_1^2\,F^2\,M_1^{-4}\,k^2\,M_1\right)^{1/2}\,m_1\\
&\asymp F\,M_1^{-2}\,m_1^2\\
&\gg F\,M_1^{-2}\left(M_1^{1+\delta_2}\,F^{-1/2}\right)^2\\
&\asymp M_1^{2\delta_2}\gg1,
\end{align*}
we may estimate
\[\exp\!\left(-A\,k^{-1}\,M_1^{-1/2}\,n^{1/2}\,m_1\right)\ll_B\left(k^{-1}\,M_1^{1/2}\,n^{1/2}\,m_1\right)^{-B}\]
for any positive integer $B$.

The error from the ``middle terms'' (i.e. the terms involving $U^{1-J}$) is, provided that $k^{-1}\,M_1^{-1/2}\,n^{1/2}\,m_1\gg1$ for $n>2n_j$,
\begin{align*}
&\ll_B k^{-1/2}\,G\,M_1^{-1/4}\,m_1^J\,U^{1-J}\\
&\qquad\cdot\sum_{n>2n_j}\left|t(n)\right|n^{-1/4}\left(k^{-1}\,M_1^{-1/2}\,n^{1/2}\,m_1\right)^{-2B}\\
&\ll F^{-1}\,G\,k^{-1/2}\,m_1^{1/2}\,U\cdot
F\,k^{2B}\,U^{-J}\,M_1^{B-1/4}\,m_1^{J-1/2-2B}\,n_j^{-B+3/4}\\
&\ll F^{-1}\,G\,k^{-1/2}\,m_1^{1/2}\,U\cdot
F\,k^{2B}\left(M_1^{1+\delta_4}\,F^{-1/2}\right)^{-J}\\
&\qquad\qquad\cdot M_1^{B-1/4}\,m_1^{J-1/2-2B}\,\left(m_1^2\,F^2\,M_1^{-4}\,k^2\,M_1\right)^{-B+3/4}\\
&\ll F^{-1}\,G\,k^{-1/2}\,m_1^{1/2}\,U\cdot
k^{3/2}\,F^{1+J/2-2B+3/2}\,M_1^{-J-\delta_4J+B-1/4+3B-3}\,
m_1^{1+J-2B}\\
&\ll F^{-1}\,G\,k^{-1/2}\,m_1^{1/2}\,U\cdot
k^{3/2}\,F^{1+J/2-2B+3/2}\,M_1^{-J-\delta_4J+B-1/4+3B-3}\\
&\qquad\qquad\cdot\left(M_1^{1+\delta_2}\,F^{-1/2}\right)^{1+J-2B}\\
&\ll F^{-1}\,G\,k^{-1/2}\,m_1^{1/2}\,U\cdot
k^{3/2}\,F^{2-B}\,M_1^{\delta_2J-\delta_4J+2B-2\delta_2B-1/4-2+\delta_2}.
\end{align*}
Choosing here $B=2$ (which is sufficiently large to make everything finite) gives
\[\ll F^{-1}\,G\,k^{-1/2}\,m_1^{1/2}\,U\cdot
M_1^{-3\delta_1/2+\delta_2J-\delta_4J+5/4-3\delta_2},\]
and this is $\ll k^{-1/2}\,m_1^{1/2}\,F^{-1}M_1^{\vartheta}\,G\,h^{3/2}\,U$, provided that $J$ is so large that the exponent of $M_1$ is not positive. Thus, the lower bound for $J$ depends on $\delta_1$, $\delta_2$ and $\delta_4$, and we must have $\delta_4>\delta_2$.

The error from the ``last term'' (not involving $U$ at all) is
\begin{align*}
&\ll_B k^{-1/2}\,G\,M_1^{-1/4}\,m_1\sum_{n>2n_j}\left|t(n)\right|n^{-1/4}\left(k^{-1}\,M_1^{-1/2}\,n^{1/2}\,m_1\right)^{-2B}\\
&\ll F^{-1}\,G\,k^{-1/2}\,m_1^{1/2}\,U\cdot
k^{2B}\,F\,M_1^{B-1/4}\,m_1^{1/2-2B}\,n_j^{-B+3/4}\\
&\ll F^{-1}\,G\,k^{-1/2}\,m_1^{1/2}\,U\cdot
k^{2B}\,F\,M_1^{B-1/4}\,m_1^{1/2-2B}\left(m_1^2\,F^2\,M_1^{-4}\,k^2\,M_1\right)^{-B+3/4}\\
&\ll F^{-1}\,G\,k^{-1/2}\,m_1^{1/2}\,U\cdot
k^{3/2}\,F^{5/2-2B}\,M_1^{4B-5/2}\left(M_1^{1+\delta_2}\,F^{-1/2}\right)^{2-4B}\\
&\ll F^{-1}\,G\,k^{-1/2}\,m_1^{1/2}\,U\cdot
k^{3/2}\,F^{3/2}\,M_1^{-1/2+2\delta_2-4\delta_2B}\\
&\asymp F^{-1}\,G\,k^{-1/2}\,m_1^{1/2}\,U\cdot
h^{3/2}\,M_1^{1+2\delta_2-4\delta_2B}\\
&\ll F^{-1}\,G\,h^{3/2}\,k^{-1/2}\,m_1^{1/2}\,U,
\end{align*}
provided that $B$ is sufficiently large, depending on $\delta_2$.

\subsubsection{Applying the saddle-point theorem: the error terms}

In this section we shall treat the error terms coming from the saddle point Theorem \ref{saddle-point-lemma}. It is applied with the parameters
\[G:=G\,M_1^{-1/4},\qquad
F:=F,\qquad
\mu:=\frac12cM_1.\]

\paragraph{The first error term.}
For a single integral, the first error term arising from the saddle point theorem is, in view of the estimates
\[1\ll\frac{M_1}{m_1}\ll\frac{M_1}U\ll F^{1/2}\,M_1^{-\delta_4},\]
at most
\begin{align*}
&\ll m_1\left(1+M_1^J\,U^{-J}\right)G\,M_1^{-1/4}\,e^{-ArM_1-AF}\\
&\ll m_1\,F^{J/2}\,M_1^{-1/4-\delta_4J}\,G\,e^{-AF}.
\end{align*}
The total error is then
\begin{align*}
&\ll k^{-1/2}\sum_{n\leqslant2n_j}\left|t(n)\right|n^{-1/4}\,m_1\,F^{J/2}\,M_1^{-1/4-\delta_4J}\,G\,e^{-AF}\\
&\ll k^{-1/2}\,F^{-1}\,G\,m_1^{1/2}\cdot n_j^{3/4}\,m_1^{1/2}\,F^{1+J/2}\,M_1^{-1/4-\delta_4J}\,e^{-AF}\\
&\ll k^{-1/2}\,F^{-1}\,G\,m_1^{1/2}\cdot m_1^{3/2}\,F^{3/2}\,k^{3/2}\,M_1^{-9/4}\,m_1^{1/2}\,F^{1+J/2}\,M_1^{-1/4-\delta_4J}\,e^{-AF}\\
&\ll k^{-1/2}\,F^{-1}\,G\,m_1^{1/2}\,h^{3/2}\cdot
m_1^2\,M_1^{-1-\delta_4J}\,F^{1+J/2}\,e^{-AF},
\end{align*}
and since
\begin{align*}
m_1^2&\,M_1^{-1-\delta_4J}\,F^{1+J/2}\,e^{-AF}\ll_B M_1^{1+\delta_4J}\,F^{1+J/2}\,F^{-B},
\end{align*}
we are done once we choose $B$ to be sufficiently large depending on $\delta_4$ and~$J$.

\paragraph{The second error term.}
Let us recall that
\[n_j-n_j'\asymp k^2\,m_j\,F^2\,U\,M_1^{-3}.\]
When $n<n_j'$, the saddle point $x_{j,n}$ lies inside the interval $\left]M_1',M_2'\right[$ and a single second error term coming from the saddle point theorem is
\[\ll G\,M_1^{-1/4}\,M_1\,F^{-3/2}.\]
In total these contribute
\begin{align*}
&\ll k^{-1/2}\sum_{n\leqslant n_j'}\left|t(n)\right|n^{-1/4}\,G\,F^{-3/2}\,M_1^{3/4}\\
&\ll k^{-1/2}\,G\,F^{-3/2}\,M_1^{3/4}\,n_j'^{3/4}\\
&\ll k^{-1/2}\,G\,F^{-3/2}\,M_1^{3/4}\,m_j^{3/2}\,F^{3/2}\,M_1^{-9/4}\,k^{3/2}\\
&\ll k^{-1/2}\,G\,F^{-1}\,m_1^{1/2}\cdot F\,m_1\,M_1^{-3/2}\,k^{3/2}\\
&\ll k^{-1/2}\,G\,F^{-1}\,m_1^{1/2}\cdot F\,m_1\,h^{3/2}\,F^{-3/2},
\end{align*}
and since $m_1\,F^{-1/2}\ll M_1\,F^{-1/2}\ll U\ll U\,L$, this is
\[\ll k^{-1/2}\,G\,F^{-1}\,m_1^{1/2}\,h^{3/2}\,U\,L\,\]
as required.

When $n_j'\leqslant n<n_j$, the saddle point $x_{j,n}$ is in the range $\left]M_1,M_1'\right]\cup\left[M_2',M_2\right[$, and a single second error term is
\[\ll G\,M_1^{3/4}\,F^{-1}.\]
Since $t(n)\ll n^{\vartheta+\varepsilon}$ for all $n$, the total contribution is
\begin{align*}
&\ll k^{-1/2}\sum_{n_j'\leqslant n<n_j}\left|t(n)\right|n^{-1/4}\,G\,M_1^{3/4}\,F^{-1}\\
&\ll k^{-1/2}\,G\,F^{-1}\,M_1^{3/4}\,n_j^{\varepsilon-1/4+\vartheta}\left(n_j-n_j'\right)\\
&\ll k^{-1/2}\,G\,F^{-1}\,L\,M_1^{3/4}\,M_1^{\varepsilon+\delta_2(\vartheta-1/4)}\left|h\right|^{\vartheta-1/4}k^{\vartheta-1/4}
\,k^2\,m_j\,F^2\,U\,M_1^{-3}\\
&\ll k^{-1/2}\,G\,F^{-1}\,L\,U\,m_1^{1/2}\left|h\right|^{3/2}
\cdot\left|h\right|^{-7/4}\,k^{7/4}\,F^{2}\,M_1^{3/4-3}\,m_1^{1/2}\\
&\qquad\cdot M_1^{\varepsilon+\delta_2(\vartheta-1/4)}\,\left|h\right|^\vartheta k^\vartheta\\
&\ll k^{-1/2}\,G\,F^{-1}\,L\,U\,m_1^{1/2}\left|h\right|^{3/2}\cdot
F^{1/4}\,M_1^{-1/2}\,m_1^{1/2}\cdot\left|h\right|^\vartheta\,k^\vartheta\\
&\ll k^{-1/2}\,G\,F^{-1}\,L\,U\,m_1^{1/2}\left|h\right|^{3/2}\cdot
F^{1/4}\,M_1^{\vartheta}
\end{align*}
which is the desired error term. 

\paragraph{The third error term.}
A single last error term is
\[\ll U^{-J}\sum_{\ell=0}^J\left(E_J(M_1+\ell\,U)+E_J(M_2-\ell\,U)\right),\]
where
\[E_J(x)=G\,M_1^{-1/4}\left(\left|p_{j,n}'(x)\right|+p_{j,n}''(x)^{1/2}\right)^{-J-1}.\]
Since $p_{j,n}''(x)\asymp F\,M_1^{-2}$, we have
\[E_J(x)\asymp G\,M_1^{-1/4}\left(F^{1/2}\,M_1^{-1}\right)^{-J-1}.\]
The total error from these error terms is
\begin{align*}
&\ll k^{-1/2}\sum_{n<2n_j}\left|t(n)\right|n^{-1/4}\,U^{-J}\,G\,M_1^{-1/4}\left(F^{1/2}\,M_1^{-1}\right)^{-J-1}\\
&\ll k^{-1/2}\,G\,U\,n_j^{3/4}\,M_1^{-1/4}\left(F^{1/2}\,U\,M_1^{-1}\right)^{-J-1}\\
&\ll k^{-1/2}\,G\,U\,m_1^{3/2}\,F^{3/2}\,k^{3/2}\,M_1^{-3/4}\,M_1^{-1/4}\left(M_1^{\delta_4}\right)^{-J-1}\\
&\ll k^{-1/2}\,G\,U\,m_1^{1/2}\,h^{3/2}\,F^{-1}\cdot m_1\,F\,M_1^{1/2}\left(M_1^{\delta_4}\right)^{-J-1}\\
&\ll k^{-1/2}\,G\,U\,m_1^{1/2}\,h^{3/2}\,F^{-1}\cdot M_1^3\left(M_1^{\delta_4}\right)^{-J-1}.
\end{align*}
Now, if $J$ is sufficiently large with respect to $\delta_4$, then this is
\[\ll k^{-1/2}\,G\,U\,m_1^{1/2}\,h^{3/2}\,F^{-1}\]
and we are done.

\subsubsection{Applying the saddle-point theorem: the main terms}

\paragraph{Obtaining the main terms.}
For each $n<n_j$ in the integral (\ref{int1}) we get a saddle-point term
\[\xi_J(x_{j,n})\,x_{j,n}^{-1/4}\left(1+\sum_{\ell=1}^Kc_\ell^{(j)}k^\ell n^{-\ell/2}x_{j,n}^{-\ell/2}\right)g(x_{j,n})\,p_{j,n}''(x_{j,n})^{-1/2}\,e\!\left(p_{j,n}(x_{j,n})+\frac18\right).\]
Substituting this back to (\ref{voronoi-sum}) gives
\begin{multline*}
i\,2^{-1/2}\,k^{-1/2}\sum_{j=1}^2(-1)^{j-1}\sum_{n<n_j}\xi_J(x_{j,n})\,t(n)\,n^{-1/4}e\left(-\frac{n\overline h}k\right)\,x_{j,n}^{-1/4}\,g(x_{j,n})\\
\cdot\left(1+\sum_{\ell=1}^Kc_\ell^{(j)}k^\ell n^{-\ell/2}x_{j,n}^{-\ell/2}\right)p_{j,n}''(x_{j,n})^{-1/2}e\!\left(p_{j,n}(x_{j,n})+\frac18\right).
\end{multline*}
This is exactly what it should be except for the term in brackets involving a sum over $\ell$, the removal of which gives an error (for each $j$ and $\ell$)
\begin{align*}
&\ll k^{-1/2}\sum_{n<n_j}\left|t(n)\right|n^{-1/4}\,M_1^{-1/4}\,G\,n^{-\ell/2}\,M_1^{-\ell/2}\,k^\ell\,F^{-1/2}\,M_1\\
&\ll k^{-1/2}\,n_j^{1/4}\,G\,F^{-1/2}\,k\,M_1^{1/4}\\
&\ll k^{-1/2}\,m_1^{1/2}\,F^{1/2}\,k^{1/2}\,M_1^{-3/4}\,G\,F^{-1/2}\,k\,M_1^{1/4}\\
&\ll k^{-1/2}\,m_1^{1/2}\,F^{-1}\,G\cdot k^{3/2}\,F^{3/2}\, F^{-1/2}\,M_1^{-1/2}\\
&\ll k^{-1/2}\,m_1^{1/2}\,F^{-1}\,G\,h^{3/2}\,M_1\,F^{-1/2}\\
&\ll k^{-1/2}\,m_1^{1/2}\,F^{-1}\,G\,h^{3/2}\,U,
\end{align*}
which is small enough.

\paragraph{The new weight functions $w_j(n)$.}
In this subsection we show that the function $w_j(n)=\xi_J(x_{j,n})$, where $n<n_j$, has the claimed properties.

The identity $w_j(n)=1$ for $n<n_j'$ follows at once from property $1$ of the function $\xi_J(x)$ on p.\ \pageref{property-one} for $M_1'<x_{j,n}<M_2'$. To prove the estimate $w_j(n)\ll1$ for $n<n_j$, we have three cases to consider. If $M_1'<x_{j,n}<M_2'$ the claim is trivial by the property $1$ of $\xi_J$.

On the other hand, if $M_1<x_{j,n}\leqslant M_1'$ the claim follows from property $2$ of $\xi_J(x)$ using the estimates $p_{j,n}''\asymp F\,M_1^{-2}$ and $U\gg F^{-\frac{1}{2}}\,M_1^{(1+\delta_4)}$:
\begin{align*}
w_j(n)
&=(J!\,U^{J})^{-1}\sum_{j=0}^{j_1}\binom{J}{j}(-1)^j\\
&\qquad\qquad\cdot\sum_{0\leqslant v
\leqslant\frac{J}{2}}c_v\,p_{j,n}''(x_{j,n})^{-v}\,(x_{j,n}-M_1')^{J-2v}\\
&\ll F^{\frac{J}{2}}M_1^{-J(1+\delta_4)}\,F^{-\frac{J}{2}}\,M_1^{-J}
\ll M_1^{-J\delta_4}\ll 1,
\end{align*}
where $j_{1}$ is the largest integer such that $M_1+j_1U<x_{j,n}$

The third case $M_2'\leqslant x_{j,n}<M_2$ is similar; we have by property $3$ on p.\ \pageref{property-three} that
\begin{align*}
w_j(n)
&=(J!\,U^{J})^{-1}\sum_{j=0}^{j_2}\binom{J}{j}(-1)^j\\
&\qquad\qquad\cdot\sum_{0\leqslant v\leqslant\frac{J}{2}}c_v\,p_{j,n}''(x_{j,n})^{-v}\,(M_2'-x_{j,n})^{J-2v}\\
&\ll F^{\frac{J}{2}}\,M_1^{J(1+\delta_4)}\,F^{-\frac{J}{2}}\,M_1^{-J}
\ll M_1^{-J\delta_4}\ll 1,
\end{align*}
where $j_2$ is the largest integer such that $M_2-j_2U>x_{j,n}$.

To check the upper bound for $w_j'(y)$,
we recall that $\xi_J'(x_0)\ll U^{-1}$ in the saddle point theorem~\ref{saddle-point-lemma}, whenever the derivative exists. Furthermore, since
\[n_j-n_j'\asymp k^2\,U\,m_j\,F^2\,M_1^{-3}\asymp\frac{n_j\,U}{m_j}\]
and
\[\frac{\mathrm dx_{j,y}}{\mathrm dy}\asymp\frac{m_j}{n_j},\]
we conclude that for $y\in\left]n_j',n_j\right[$, where the derivative exists,
\[w_j'(y)=\frac{\mathrm d}{\mathrm dy}\,\xi_J(x_{j,y})
=\xi_J'(x_{j,y})\cdot\frac{\mathrm dx_{j,y}}{\mathrm dy}
\ll\frac1U\cdot\frac{m_j}{n_j}\asymp\frac1{n_j-n_j'}.\]

\section{Estimates for non-linear sums}

The savings in the estimates for short sums depend on an estimate for the kind of nonlinear sums that appear after the application of the Voronoi type summation formula. In the following theorem, it is essential that the estimate is better when shorter sums are considered.
\begin{theorem}\label{nonlinear-estimate}
Let $M\in\left[1,\infty\right[$, $\eta\in\mathbb R$, $B\in\mathbb R$, and $\Delta\in\left[1,M\right]$. Denote $F=\left|B\right|M^{1/2}$, and assume that
\[M^2\ll\Delta\,F.\]
Let $g$ be a $C^1$-function on the interval $\left[M,M+\Delta\right]$ satisfying bounds
\[g(x)\ll G\qquad\text{and}\qquad g'(x)\ll G'\]
on $\left[M,M+\Delta\right]$ for some positive real numbers $G$ and $G'$. Then
\[\sum_{M\leqslant m\leqslant M+\Delta}t(m)\,g(m)\,e\!\left(\eta m+Bm^{1/2}\right)\ll
\Delta^{5/6}\left(G+\Delta\,G'\right)M^{\vartheta-1/3}\,F^{1/3+\varepsilon}.\]
\end{theorem}

\paragraph{Proof.} This is analogous to Theorem 4.1 in \cite{Ernvall-Hytonen--Karppinen2008}, and in fact, the proof given in \cite{Ernvall-Hytonen--Karppinen2008} works almost verbatim in our case, except that now we use Theorem \ref{general-transformation-formula-for-cusp-forms} instead of the corresponding result for holomorphic cusp form, i.e.\ Theorem 3.4 in \cite{Jutila1987a}, and naturally, when smoothing error is to be estimated, an extra $M^\vartheta$ appears in a few places. There is only one point which requires extra clarification, the error term in Theorem \ref{general-transformation-formula-for-cusp-forms} has the extra factor $F^{1/4}\,M^\vartheta$; this time the total error from the error terms coming from using the transformation formula contributes
\[\ll\frac\Delta M\,M^{1/2+\vartheta}\,F^{1/4+\varepsilon}
\ll \Delta^{5/6}\,M^{\vartheta-1/3}\,F^{1/4+\varepsilon},\]
which is smaller than the desired upper bound.

\bigbreak\noindent
It turns out that for long sums, the $\vartheta$ in the upper bound may be erased. This was proved by Karppinen in \cite{Karppinen1998} by considering the mean value of the relevant exponential sums. Earlier, Jutila \cite{Jutila1989, Jutila1990} had considered similar mean values for holomorphic cusp forms and the divisor function. The following estimate is Theorem 8.2 in \cite{Karppinen1998}.
\begin{theorem}\label{long-nonlinear-estimate}
Let $M\in\left[1,\infty\right[$, $\eta\in\mathbb R$, $B\in\mathbb R$, and $\Delta\in\left[1,M\right]$. Denote $F=\left|B\right|M^{1/2}$, and assume that $M\ll F$.
Let $g$ be a $C^1$-function on the interval $\left[M,M+\Delta\right]$ satisfying the bounds
\[g(x)\ll G\qquad\text{and}\qquad g'(x)\ll G'\]
on $\left[M,M+\Delta\right]$ for some positive real numbers $G$ and $G'$. Then
\[\sum_{M\leqslant m\leqslant M+\Delta}t(m)\,g(m)\,e\!\left(\eta m+Bm^{1/2}\right)\ll
\left(G+\Delta\,G'\right)M^{1/2}\,F^{1/3+\varepsilon}.\]
\end{theorem}
In fact, the proofs of our main theorems do not require Theorem \ref{long-nonlinear-estimate} as we could use Theorem \ref{nonlinear-estimate} instead. However, one of the upper bounds in Theorem \ref{estimate-for-smooth-short-sums} is better if Theorem \ref{long-nonlinear-estimate} is applied instead of Theorem \ref{nonlinear-estimate}.

\section{Proof of Theorem \ref{short-estimate}}

We shall prove Theorem 1 by first proving estimates for smooth short exponential sums.  For this purpose, we shall use a wide weight function $w\in C_{\mathrm c}^\infty(\mathbb R_+)$ taking only values from $\left[0,1\right]$, supported in $\left[M,M+\Delta\right]$,
and for which
\[w^{(\nu)}(x)\ll_\nu\Delta^{-\nu},\]
for every nonnegative integer~$\nu$. The following estimates give an analogue of Theorem 5.1 of \cite{Ernvall-Hytonen--Karppinen2008}.
\begin{theorem}\label{estimate-for-smooth-short-sums}
Let $M\in\left[1,\infty\right[$, and let $\Delta\in\left[1,M\right]$ with $\Delta\gg M^\beta$ for some arbitrarily small fixed $\beta\in\mathbb R_+$. Furthermore, let $\alpha\in\mathbb R$, and let $h\in\mathbb Z$, $k\in\mathbb Z_+$ and $\eta\in\mathbb R$ be such that
\[\alpha=\frac hk+\eta,\quad(h,k)=1,\quad k\leqslant K,\quad\left|\eta\right|\leqslant\frac1{kK},\]
where $K=\Delta^{1/2-\delta}$ for an arbitrarily small fixed $\delta\in\mathbb R_+$.
\begin{enumerate}
\item
If $\eta\ll\Delta^{-1+\delta}$, then
\[\sum_{M\leqslant n\leqslant M+\Delta}t(n)\,e(n\alpha)\,w(n)\ll_{\beta,\delta}\Delta^{1/6}\,M^{1/3+\varepsilon}.\]
\item
If $\Delta^{-1+\delta}\ll\eta$ and $k^2\,\eta^2\,M<1/2$, then
\[\sum_{M\leqslant n\leqslant M+\Delta}t(n)\,e(n\alpha)\,w(n)\ll_{\beta,\delta}1.\]
\item
If $\Delta^{-1+\delta}\ll\eta\ll M\,\Delta^{-2}$, $k^2\,\eta^2\,M\gg1$ and $k^2\,\eta\,M\,\Delta^{-1+\delta}\ll1$, then
\[\sum_{M\leqslant n\leqslant M+\Delta}t(n)\,e(n\alpha)\,w(n)\ll_{\beta,\delta}1+k^{-1/2}\,\Delta\,M^{-1/4}\left(k^2\,\eta^2\,M\right)^{\vartheta-1/4+\varepsilon}.\]
\item
If $\Delta^{-1+\delta}\ll\eta\ll M\,\Delta^{-2}$, $k^2\,\eta^2\,M\gg1$ and $k^2\,\eta\,M\,\Delta^{-1+\delta}\gg1$, then
\[\sum_{M\leqslant n\leqslant M+\Delta}t(n)\,e(n\alpha)\,w(n)\ll_{\beta,\delta}\left(k^2\,\eta^2\,M\right)^\vartheta\Delta^{1/6}\,M^{1/3+\varepsilon}.\]
\end{enumerate}
\end{theorem}

\paragraph{Remark.} The proof of the estimate 1 employs Theorem \ref{long-nonlinear-estimate}. We could use Theorem \ref{nonlinear-estimate} instead to obtain the upper bound $\ll\Delta^{1/6-\vartheta}\,M^{1/3+\vartheta+\varepsilon}$, which would be good enough to obtain the main theorems.

\paragraph{Proof.} We begin by applying the Voronoi summation formula to the sum under study. The proof will soon split into two cases depending on whether $\eta$ is smaller than larger than $\Delta^{-1+\delta}$. Voronoi summation yields
\begin{align*}
&\sum_{M\leqslant n\leqslant M+\Delta}t(n)\,e(n\alpha)\,w(n)\\
&=\frac{\pi i}{k\,\sinh\pi\kappa}\sum_{n=1}^\infty t(n)\,e\!\left(\frac{-n\overline h}k\right)\int\limits_M^{M+\Delta}\left(J_{2i\kappa}-J_{-2i\kappa}\right)\!\left(\frac{4\pi\sqrt{nx}}k\right)e(\eta x)\,w(x)\,\mathrm dx\\
&\qquad+\frac{4\cosh\pi\kappa}k\sum_{n=1}^\infty t(-n)\,e\!\left(\frac{n\overline h}k\right)\int\limits_M^{M+\Delta}K_{2i\kappa}\!\left(\frac{4\pi\sqrt{nx}}k\right)e(\eta x)\,w(x)\,\mathrm dx.
\end{align*}

Using the asymptotics for the $K$-Bessel function, and picking some large $A\in\mathbb R_+$, the $K$-series can be estimated by
\[\ll_A\frac1k\sum_{n=1}^\infty\left|t(n)\right|\int\limits_M^{M+\Delta}k^A\,n^{-A/2}\,x^{-A/2}\,w(x)\,\mathrm dx.\]
This is $\ll_A k^{A-1}\,\Delta\,M^{-A/2}$, provided that $A>2$, and since $k\ll M^{1/2-\delta}$, it is furthermore $\ll_\delta1$, provided that $A\gg_\delta1$.
Similarly, by replacing the $J$-Bessel expression by the asymptotics given in \eqref{Voronoi-J-asymptotics} with some $K\in\mathbb Z_+$, the resulting $O$-terms contribute
\[\ll_K\frac1k\sum_{n=1}^\infty\left|t(n)\right|\int\limits_M^{M+\Delta}k^{1/2+(K+1)}\,n^{-1/4-(K+1)/2}\,x^{-1/4-(K+1)/2}\,w(x)\,\mathrm dx,\]
and this is again $\ll_\delta1$ for a fixed $K\gg_\delta1$. Thus, we are led to
\begin{align*}
&\sum_{M\leqslant n\leqslant M+\Delta}t(n)\,e(n\alpha)\,w(n)\\
&=O(1)+\frac Ck\sum_{n=1}^\infty t(n)\,e\!\left(\frac{-n\overline h}k\right)\\
&\qquad\qquad\cdot
\int\limits_M^{M+\Delta}\frac{k^{1/2}}{n^{1/4}\,x^{1/4}}\sum_\pm (\pm1)e\!\left(\pm\frac{2\sqrt{nx}}k\right)e\!\left(\mp\frac18\right)g_\pm(x;n,k)\,e(\eta x)\,w(x)\,\mathrm dx,
\end{align*}
where
\[g_\pm(x;n,k)=1+\sum_{\ell=1}^Kc_\ell^\pm\,k^\ell\,n^{-\ell/2}\,x^{-\ell/2},\]
and $C$ is some real constant.

\subsection{The case $\eta\ll\Delta^{-1+\delta}$}
Write $X=k^2\,M\,\Delta^{3\delta-2}$. We shall handle separately the terms with $n>X$ and the terms with $n\leqslant X$.

The high-frequency terms with $n>X$ contribute
\begin{multline*}
\ll\frac1k\sum_{n>X}t(n)\,k^{1/2}\,n^{-1/4}\,e\!\left(\frac{-n\overline h}k\right)\\
\cdot\int\limits_M^{M+\Delta}x^{-1/4}\,g_\pm(x;n,k)\,e\!\left(\pm\frac{2\sqrt{nx}}k+\eta x\right)w(x)\,\mathrm dx.
\end{multline*}
Since we now have $X^{1/2}\,k^{-1}\,M^{-1/2}\gg\Delta^\delta\,\eta$, Lemma \ref{jutila-motohashi-lemma} says that the integral here is
\[\int\limits_M^{M+\Delta}\ldots\mathrm dx\ll_P
M^{-1/4}\left(\Delta\,n^{1/2}\,k^{-1}\,M^{-1/2}\right)^{-P}\Delta.\]
Provided that $P\geqslant2$, the contribution from these high-frequency terms is
\begin{align*}
&\ll\frac1k\sum_{n>X}\left|t(n)\right|\frac{k^{1/2}}{n^{1/4}}\,M^{-1/4}\left(\Delta\,n^{1/2}\,k^{-1}\,M^{-1/2}\right)^{-P}\Delta\\
&\ll k^{P-1/2}\,\Delta^{1-P}\,M^{P/2-1/4}\sum_{n>X}\left|t(n)\right|\,n^{-1/4-P/2}\\
&\ll k^{P-1/2}\,\Delta^{1-P}\,M^{P/2-1/4}\,X^{3/4-P/2},
\end{align*}
and for a fixed $P\gg_\delta1$, this is $\ll_\delta1$.

Let us consider next the low-frequency terms with $n\leqslant X$. These contribute
\begin{align*}
&\ll\frac1k\sum_{n\leqslant X}t(n)\,k^{1/2}\,n^{-1/4}\,e\!\left(\frac{-n\overline h}k\right)\\
&\qquad\cdot\int\limits_M^{M+\Delta}x^{-1/4}\,g_\pm(x;n,k)\,e\!\left(\pm\frac{2\sqrt{nx}}k+\eta x\right)w(x)\,\mathrm dx\\
&=k^{-1/2}\sum_{\substack{L\leqslant X/2\\\mathrm{dyadic}}}
\int\limits_M^{M+\Delta}x^{-1/4}\,w(x)\\
&\qquad\cdot\sum_{L<n\leqslant2L}t(n)\,n^{-1/4}\,g_\pm(x;n,k)\,e\!\left(\pm\frac{2\sqrt{nx}}k-\frac{n\overline h}k+\eta x\right)\mathrm dx.
\end{align*}
By Theorem \ref{long-nonlinear-estimate}, the conditions of which are met under the present circumstances, the sum $\sum_{L<n\leqslant2L}$ can be estimated by
\[\ll L^{-1/4}\,L^{1/2}\left(M^{1/2}\,k^{-1}\,L^{1/2}\right)^{1/3+\varepsilon}
\ll L^{5/12}\,M^{1/6+\varepsilon}\,k^{-1/3}.\]
Thus, the low-frequency terms contribute
\begin{align*}
&\ll k^{-1/2}\sum_{\substack{L\leqslant X/2\\\mathrm{dyadic}}}\Delta\,M^{-1/4}\,L^{5/12}\,M^{1/6+\varepsilon}\,k^{-1/3}\\
&\ll k^{-5/6}\,\Delta\,M^{-1/12+\varepsilon}\,X^{5/12}
\ll\Delta^{1/6}\,M^{1/3+\varepsilon},
\end{align*}
and we are finished with the case $\eta\ll\Delta^{-1+\delta}$.

\subsection{The case $\eta\gg\Delta^{-1+\delta}$}
This time we will choose $X=k^2\,\eta^2\,M$. The high-frequency terms with $n>2X$ are again handled in the same way as in the case $\eta\ll\Delta^{-1+\delta}$. For an integer $P\geqslant2$, we have
\begin{align*}
&\ll\frac1k\sum_{n>2X}t(n)\,k^{1/2}\,n^{-1/4}\,e\!\left(\frac{-n\overline h}k\right)\\
&\qquad\cdot\int\limits_M^{M+\Delta}x^{-1/4}\,g_\pm(x;n,k)\,e\!\left(\pm\frac{2\sqrt{nx}}k+\eta x\right)w(x)\,\mathrm dx\\
&\ll_P\frac1k\sum_{n>2X}\left|t(n)\right|\,k^{1/2}\,n^{-1/4}\,M^{-1/4}\left(\Delta\,n^{1/2}\,k^{-1}\,M^{-1/2}\right)^{-P}\Delta\\
&\ll_P k^{-1/2+P}\,\Delta^{1-P}\,M^{P/2-1/4}\,X^{3/4-P/2}.
\end{align*}
For $P\gg_\delta1$, this contribution is again $\ll_\delta1$.

If $X<1/2$, then the above already proves case 2 of the theorem, so let us assume that $X\gg1$. The remaining terms, the ones with $n\leqslant2X$, are then partitioned into two sets: those with $\left|n-X\right|\geqslant W$ and those with $\left|n-X\right|<W$, where $W=k^2\,M\,\eta\,\Delta^{-1+\delta}$.

So, let us consider the terms with $n\leqslant2X$ and $\left|n-X\right|\geqslant W$. The crucial observations here are that
\[\frac{\sqrt X}{k\sqrt x}-\frac{\sqrt n}{k\sqrt x}
\asymp\frac1{k\sqrt M}\int\limits_n^X\frac{\mathrm dt}{\sqrt t}
\gg\frac{\left|n-X\right|}{k\sqrt M\sqrt X}\gg\frac W{k\sqrt M\sqrt X}=\Delta^{-1+\delta},\]
and that, thanks to the assumption $\eta\ll M\,\Delta^{-2}$,
\[\left|\eta\right|-\frac{\sqrt X}{k\sqrt x}
=\frac{\sqrt X}{k\sqrt M}-\frac{\sqrt X}{k\sqrt x}
\asymp\frac{\sqrt X}k\int\limits_M^x\frac{\mathrm dt}{t^{3/2}}
\ll\frac{\sqrt X\,\Delta}{k\,M^{3/2}}=\frac{\eta\,\Delta}M\ll\Delta^{-1}.\]
Using these appropriately (depending on the sign of $\eta$), we conclude that
\[\frac{\mathrm d}{\mathrm dx}\left(\pm\frac{2\sqrt{nx}}k+\eta x\right)
=\pm\frac{\sqrt n}{k\sqrt x}+\eta\gg\Delta^{-1+\delta},\]
and so, by Lemma \ref{jutila-motohashi-lemma}, the terms under consideration contribute
\begin{align*}
&\ll_P\frac1k\sum_{\substack{n\leqslant2X,\\\left|n-X\right|\geqslant W}}\left|t(n)\right|\,k^{1/2}\,n^{-1/4}\,M^{-1/4}\,\Delta^{-\delta\,P}\,\Delta\\
&\ll_Pk^{-1/2}\,X^{3/4}\,M^{-1/4}\,\Delta^{-\delta\,P}\,\Delta,
\end{align*}
and for a fixed $P\gg_\delta1$ this is again $\ll_\delta1$.

Next, if $W\ll1$, then the remaining terms, the ones with $\left|n-X\right|<W$, contribute
\[\ll k^{-1/2}\,X^{\vartheta-1/4+\varepsilon}\,\Delta\,M^{-1/4}
\ll k^{-1/2}\,(k^2\eta^2M)^{\vartheta-1/4+\varepsilon}\,\Delta\,M^{-1/4},\]
and we have established case 3. Finally, only case 4 remains.

So, let us assume that $W\gg1$. The idea now is to exchange integration and summation, apply Theorem \ref{nonlinear-estimate} to the integrand with the parameters
\[M=X,\quad\Delta=W,\quad\text{and}\quad B=\frac{\sqrt x}k,\]
observing that the condition $\Delta\,F\gg M^2$ of Theorem \ref{nonlinear-estimate} holds, since it reduces to
\[W\cdot\frac{\sqrt X\sqrt M}k\gg X,\]
which follows from $k^2\,\eta\ll1\ll W$.
The remaining terms are then seen to contribute
\begin{align*}
&\ll\frac1k\sum_{X-W<n<X+W}t(n)\,k^{1/2}\,n^{-1/4}\,e\!\left(\frac{-n\overline h}k\right)\\
&\qquad\cdot\int\limits_M^{M+\Delta}x^{-1/4}\,g_\pm(x;n,k)\,e\!\left(\pm\frac{2\sqrt{nx}}k+\eta x\right)w(x)\,\mathrm dx\\
&\ll k^{-1/2}\int\limits_M^{M+\Delta}x^{-1/4}\,e(\eta x)\,w(x)\\
&\qquad\cdot\sum_{X-W<n<X+W}t(n)\,n^{-1/4}\,e\!\left(-\frac{n\overline h}k\pm\frac{2\sqrt{nx}}k\right)\mathrm dx\\
&\ll k^{-1/2}\int\limits_M^{M+\Delta}x^{-1/4}\,w(x)
\left(\frac WX\right)^{5/6}\,X^{-1/4}\,X^{1/2+\vartheta}\,\left(\frac{\sqrt X\sqrt M}k\right)^{1/3+\varepsilon}\,\mathrm dx\\
&\ll k^{-1/2}\,\Delta\,M^{-1/4}\left(\frac{k^2\,\eta\,M\,\Delta^{-1+\delta}}{k^2\,\eta^2\,M}\right)^{5/6}\left(k^2\,\eta^2\,M\right)^{1/4+\vartheta}\left(\frac{k\,\eta\,M}k\right)^{1/3+\varepsilon}\\
&\ll k^{-1/2}\,\Delta\,M^{-1/4}\left(\frac{\Delta^{-1+\delta}}\eta\right)^{5/6}k^{1/2+2\vartheta}\,\eta^{1/2+2\vartheta}\,M^{1/4+\vartheta}\,\eta^{1/3}\,M^{1/3+\varepsilon}\\
&\ll\left(k^2\,\eta^2\,M\right)^\vartheta\,\Delta^{1/6}\,M^{1/3+\varepsilon},
\end{align*}
and we are done.

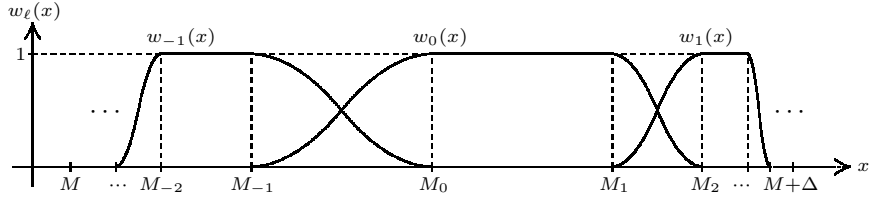
\begin{figure}[h]
\begin{center}
\setlength{\unitlength}{.5cm}
\begin{picture}(22,5.3)(-1,-.8)
\qbezier(-.5,0)(5.5,0)(21.5,0)
\qbezier(0,-.5)(0,1.75)(0,3.8)
\qbezier(1,-.1)(1,0)(1,.1)
\qbezier(2.1875,-.1)(2.1875,0)(2.1875,.1)
\qbezier(3.375,-.1)(3.375,0)(3.375,.1)
\qbezier(5.75,-.1)(5.75,0)(5.75,.1)
\qbezier(10.5,-.1)(10.5,0)(10.5,.1)
\qbezier(15.25,-.1)(15.25,0)(15.25,.1)
\qbezier(17.615,-.1)(17.615,0)(17.615,.1)
\qbezier(18.8075,-.1)(18.8075,0)(18.8075,.1)
\qbezier(19.40325,-.1)(19.40325,0)(19.40325,.1)
\qbezier(20,-.1)(20,0)(20,.1)
\put(3,3.3){$\scriptstyle{w_{-1}(x)}$}
\put(10,3.3){$\scriptstyle{w_0(x)}$}
\put(17,3.3){$\scriptstyle{w_1(x)}$}
\put(1.5,1.3){$\cdots$}
\put(19.5,1.3){$\cdots$}
\put(.7,-.6){$\scriptstyle M$}
\put(2,-.6){$\scriptstyle{\cdots}$}
\put(2.8,-.6){$\scriptstyle{M_{-2}}$}
\put(5.2,-.6){$\scriptstyle{M_{-1}}$}
\put(10.15,-.6){$\scriptstyle{M_0}$}
\put(14.9,-.6){$\scriptstyle{M_1}$}
\put(17.3,-.6){$\scriptstyle{M_2}$}
\put(18.4,-.6){$\scriptstyle{\cdots}$}
\put(19.2,-.6){$\scriptstyle{M+\Delta}$}
\qbezier(21.5,0)(21.3,.05)(21.1,.2)
\qbezier(21.5,0)(21.3,-.05)(21.1,-.2)
\qbezier(0,3.8)(.05,3.6)(.2,3.4)
\qbezier(0,3.8)(-.05,3.6)(-.2,3.4)
\qbezier(-.1,3)(0,3)(.1,3)
\put(-.45,2.85){$\scriptstyle1$}
\put(21.7,-.1){$\scriptstyle x$}
\put(-.65,4){$\scriptstyle{w_\ell(x)}$}
\multiput(0,0)(.25,0){70}{\qbezier(0,3)(.05,3)(.1,3)}
\multiput(3.375,0)(0,.25){12}{\qbezier(0,0)(0,.05)(0,.1)}
\multiput(5.75,0)(0,.25){12}{\qbezier(0,0)(0,.05)(0,.1)}
\multiput(10.5,0)(0,.25){12}{\qbezier(0,0)(0,.05)(0,.1)}
\multiput(15.25,0)(0,.25){12}{\qbezier(0,0)(0,.05)(0,.1)}
\multiput(17.615,0)(0,.25){12}{\qbezier(0,0)(0,.05)(0,.1)}
\multiput(18.8075,0)(0,.25){12}{\qbezier(0,0)(0,.05)(0,.1)}
\linethickness{.2mm}
\qbezier(2.1875,0)(2.484375,.05)(2.78125,1.5)
\qbezier(2.78125,1.5)(3.078125,2.95)(3.375,3)
\qbezier(5.75,0)(6.9375,.05)(8.125,1.5)
\qbezier(8.125,1.5)(9.3125,2.95)(10.5,3)
\qbezier(5.75,3)(6.9375,2.95)(8.125,1.5)
\qbezier(8.125,1.5)(9.3125,.05)(10.5,0)
\qbezier(15.25,3)(15.84125,2.95)(16.4325,1.5)
\qbezier(16.4325,1.5)(17.02375,.05)(17.615,0)
\qbezier(15.25,0)(15.84125,.05)(16.4325,1.5)
\qbezier(16.4325,1.5)(17.02375,2.95)(17.615,3)
\qbezier(18.8075,3)(18.9564375,2.95)(19.105375,1.5)
\qbezier(19.105375,1.5)(19.2543125,.05)(19.40325,0)
\linethickness{.25mm}
\qbezier(3.375,3)(4.5625,3)(5.75,3)
\qbezier(10.5,3)(12.615,3)(15.25,3)
\qbezier(17.615,3)(18.21125,3)(18.8075,3)
\end{picture}
\end{center}
\caption{\label{partition-of-unity}The weight functions $w_\ell(x)$ used in the proofs of Theorem \ref{short-estimate}, Proposition \ref{weight-function-removal-for-longer-sums} and Lemma \ref{weight-function-removal-from-last-voronoi-argument}.}
\end{figure}

\paragraph{Proof of Theorem \ref{short-estimate}.}
We can now remove the weight function $w$ from the estimates for short sums.
For this purpose we shall introduce a partition of unity of $\left]M,M+\Delta\right[$. Let us define a set of points $M_\ell$ for $\ell\in\mathbb Z$ by first setting
\[M_0=M+\frac\Delta2,\]
and then for each $\ell\in\mathbb Z_+$
\[M_{\pm\ell}=M+\frac\Delta2\pm\left(\frac\Delta4+\frac\Delta8+\ldots+\frac\Delta{2^{\ell+1}}\right).\]
We pick functions $w_\ell\in C_{\mathrm c}^\infty(\mathbb R)$ such that
each $w_\ell$ only takes values from $\left[0,1\right]$, $w_\ell$ is supported on $\left[M_{2\ell-1},M_{2\ell+2}\right]$, $w_\ell\equiv1$ on $\left[M_{2\ell},M_{2\ell+1}\right]$, and
\[w_\ell^{(\nu)}(x)\ll_\nu\left(\frac\Delta{4^{\left|\ell\right|}}\right)^{\!-\nu},\]
for $x\in\left[M_{2\ell-1},M_{2\ell+2}\right]$, uniformly in $\ell$. 
Furthermore, $w_\ell+w_{\ell+1}$ is to be $\equiv1$ on $\left[M_{2\ell+1},M_{2\ell+2}\right]$. Figure \ref{partition-of-unity} depicts the situation.

Let now $\Delta$ be such that $\Delta^{3/2+\delta}\ll M$ for some arbitrarily small $\delta\in\mathbb R_+$, and let $L\in\mathbb Z_+$ be such that
\[\Delta\,2^{-L}=M^{2/(5+6\vartheta)}.\]
Since $\Delta^{-1/2+\delta}\ll M\,\Delta^{-2}$, we have for any Farey approximation $\alpha=h/k+\eta$ of order $\Delta^{1/2-\delta}$ that
\[\left|\eta\right|\leqslant\frac1{k\,\Delta^{1/2-\delta}}\leqslant\frac1{\Delta^{1/2-\delta}}\ll\frac M{\Delta^2}.\]
Thus, we may apply Theorem \ref{estimate-for-smooth-short-sums} to get
\[\sum_{n\in\mathbb Z}t(n)\,e(n\alpha)\,w_\ell(n)\ll
\Delta^{1/6-\vartheta}\,M^{1/3+\vartheta+\varepsilon}
+\Delta\,M^{-1/4}
\ll\Delta^{1/6-\vartheta}\,M^{1/3+\vartheta+\varepsilon},\]
and so
\[\sum_{\ell=-L}^L\sum_{n\in\mathbb Z}t(n)\,e(n\alpha)\,w_\ell(n)
\ll\sum_{\ell=-L}^L\left(\frac{\Delta}{4^{\left|\ell\right|}}\right)^{\!1/6-\vartheta}M^{1/3+\vartheta+\varepsilon}\ll\Delta^{1/6-\vartheta}\,M^{1/3+\vartheta+\varepsilon},\]
and estimating by absolute values,
\begin{multline*}
\sum_{M\leqslant n\leqslant M+\Delta}t(n)\,e(n\alpha)\left(1-\sum_{\ell=-L}^Lw_\ell(n)\right)\\\ll M^{2/(5+6\vartheta)}\,M^{\vartheta+\varepsilon}\ll\Delta^{1/6-\vartheta}\,M^{1/3+\vartheta+\varepsilon}.
\end{multline*}

\begin{proposition}\label{weight-function-removal-for-longer-sums}
Let $M\in\left[1,\infty\right[$ and let $\Delta\in\left[1,M\right]$ satisfy $\Delta\gg M^{2/3}$. Also, let $\alpha\in\mathbb R$ have a rational approximation $\alpha=h/k+\eta$, where $h$ and $k$ are coprime integers with $1\leqslant k\ll M^{1/3-\varepsilon}$, and where $\eta\in\mathbb R$ satisfies $\eta\ll k^{-1}\,\Delta^{\varepsilon-1/2}$ and $\eta\ll M\,\Delta^{-2}$. Then
\[\sum_{M\leqslant n\leqslant M+\Delta}t(n)\,e(n\alpha)
\ll\Delta^{1/6-\vartheta}\,M^{1/3+\vartheta+\varepsilon}
+k^{-1/2}\,\Delta\,M^{-1/4}\left(k^2\,\eta^2\,M\right)^{\vartheta-1/4+\varepsilon}.\]
\end{proposition}

\paragraph{Proof.} This is very similar to the proof of Theorem \ref{short-estimate} above. In particular, we may use the same weight functions $w_\ell$, and we simply have an extra term on the right-hand side.

\bigbreak
We also need a slightly more complicated version:
\begin{proposition}\label{weight-removal-for-longer-sums-with-many-farey-fractions}
Let $M\in\left[1,\infty\right[$, let $\Delta\in\left[1,M\right]$ with $\Delta\gg M^{2/3}$, and let $\delta\in\left]0,1/2\right[$ be fixed. Also, let $\alpha\in\mathbb R$ have rational approximations
\[\alpha=\frac{h_1}{k_1}+\eta_1=\frac{h_2}{k_2}+\eta_2=\ldots=\frac{h_L}{k_L}+\eta_L,\]
where $L\in\mathbb Z_+$ is chosen so that $\Delta\,2^{-L}\asymp M^{2/3}$, and that $h_1,h_2,\ldots,h_L\in\mathbb Z$, $k_1,k_2,\ldots,k_L\in\mathbb Z_+$ with
\[k_1\leqslant\Delta^{1/2-\delta},\quad k_2\leqslant\left(\frac\Delta2\right)^{1/2-\delta},\quad\ldots,\quad k_L\leqslant\left(\frac\Delta{2^L}\right)^{1/2-\delta},\]
and $\eta_1,\eta_2,\ldots,\eta_L\in\mathbb R$ with
\[\left|\eta_1\right|\leqslant k_1^{-1}\,\Delta^{\delta-1/2},\quad
\left|\eta_2\right|\leqslant k_2^{-1}\left(\frac\Delta2\right)^{\delta-1/2},\quad
\ldots,\quad
\left|\eta_L\right|\leqslant k_L^{-1}\left(\frac\Delta{2^L}\right)^{\delta-1/2},\]
and assume that $\eta_\ell\ll M\,(\Delta\,4^{-\ell})^{-2}$ for each $\ell\in\left\{1,2,\ldots,L\right\}$.
Then
\begin{multline*}
\sum_{M\leqslant n\leqslant M+\Delta}t(n)\,e(n\alpha)
\ll_\delta\Delta^{1/6-\vartheta}\,M^{1/3+\vartheta+\varepsilon}\\
+\sum_{\ell=1}^Lk_\ell^{-1/2}\cdot\frac\Delta{2^L}\cdot M^{-1/4}\left(k_\ell^2\,\eta_\ell^2\,M\right)^{\vartheta-1/4+\varepsilon}.
\end{multline*}
Furthermore, here the term corresponding to a given value of $\ell$ can be deleted unless $k_\ell^2\,\eta_\ell\,M\,\Delta^{-1+\delta}\ll1\ll k_\ell^2\,\eta_\ell^2\,M$.
\end{proposition}
\paragraph{Proof.}
Again, the proof is very much similar to the proof of Theorem \ref{short-estimate}, but in this case each subsum $\sum t(n)\,e(n\alpha)\,w_\ell(n)$ is estimated with a Farey approximation appropriate for the length of the support of $w_\ell$, which is $\asymp\Delta\,4^{-\left|\ell\right|}$. 

\paragraph{Proof of Theorem \ref{holomorphic-short-estimate}.}
These estimates follow from the proofs of Theorems 5.5 and 5.7 in \cite{Ernvall-Hytonen--Karppinen2008}, except that Theorem 5.1 should be modified a little \cite{Ernvall-Hytonen--communication}. The term $k^{-1}\,\Delta\left|\eta\right|^{-1/2}M^{-1/2+\varepsilon}$ only appears in the case in which $k^2\,\eta^2\,M\gg1$, in which case the relevant estimate (on p.~27 of \cite{Ernvall-Hytonen--Karppinen2008}) is actually
\begin{multline*}
\widetilde A(M,\Delta,\alpha)\ll\Delta^{1/6}\,M^{1/3+\varepsilon}
+k^{-1/2}\,\Delta\,M^{-1/4}\left(k^2\,\eta^2\,M\right)^{\varepsilon-1/4}\\
\ll\Delta^{1/6}\,M^{1/3+\varepsilon}
+k^{-1/2}\,\Delta\,M^{-1/4}.
\end{multline*}
Thus, the second upper bound of Theorem 5.1 is actually
\[\sum_{M\leqslant n\leqslant M+\Delta}a(n)\,e(n\alpha)\,w(n)\ll\Delta^{1/6}\,M^{1/3+\varepsilon}+k^{-1/2}\,\Delta\,M^{-1/4+\varepsilon}.\]

\section{Proof of Theorem \ref{improved-estimate} and Corollary \ref{holomorphic-corollary}}

Let $U\in\mathbb R_+$. We shall pick a weight function $w\in C_{\mathrm c}^\infty(\mathbb R_+)$ taking only nonnegative real values, supported in $\left[M,M+\Delta\right]$, identically equal to $1$ in $\left[M+U,M+\Delta-U\right]$, for which
\[w^{(\nu)}(x)\ll U^{-\nu},\]
and whose derivatives are supported in $\left[M,M+U\right]\cup\left[M+\Delta-U,M+\Delta\right]$. Sums with this weight function can be estimated rather nicely:
\begin{lemma}\label{smooth-improved-estimates}
Let $X\in\mathbb R_+$ with $X\gg1$, let $M\in\left[1,\infty\right[$, and let $U$ and $w$ as above. Also, let $h$ and $k$ be coprime integers with $1\leqslant k\ll M^{1/2-\delta}$, where $\delta$ is a fixed positive real number. Then
\[\sum_{M\leqslant n\leqslant M+\Delta}t(n)\,e\!\left(\frac{nh}k\right)w(n)
\ll_\delta k^{1/2}\,X^{1/4}\,M^{1/4}+k^{3/2}\,X^{-1/4}\,M^{3/4}\,U^{-1}.\]
Furthermore, if we select $X=1/2$, we can forget the first term.
\end{lemma}

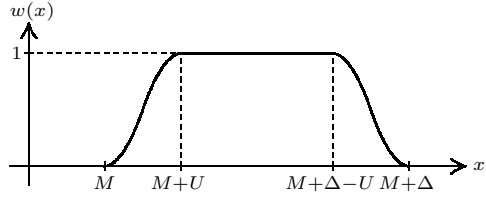
\begin{figure}[h]
\begin{center}
\setlength{\unitlength}{.5cm}
\begin{picture}(12,5.8)(-1,-.8)
\qbezier(-.5,0)(5.5,0)(11.5,0)
\qbezier(0,-.5)(0,1.75)(0,3.8)
\qbezier(2,-.1)(2,0)(2,.1)
\qbezier(4,-.1)(4,0)(4,.1)
\qbezier(8,-.1)(8,0)(8,.1)
\qbezier(10,-.1)(10,0)(10,.1)
\put(1.7,-.6){$\scriptstyle M$}
\put(3.2,-.6){$\scriptstyle{M+U}$}
\put(6.75,-.6){$\scriptstyle{M+\Delta-U}$}
\put(9.2,-.6){$\scriptstyle{M+\Delta}$}
\qbezier(11.5,0)(11.3,.05)(11.1,.2)
\qbezier(11.5,0)(11.3,-.05)(11.1,-.2)
\qbezier(0,3.8)(.05,3.6)(.2,3.4)
\qbezier(0,3.8)(-.05,3.6)(-.2,3.4)
\qbezier(-.1,3)(0,3)(.1,3)
\put(-.45,2.85){$\scriptstyle1$}
\put(11.7,-.1){$\scriptstyle x$}
\put(-.5,4){$\scriptstyle{w(x)}$}
\multiput(0,0)(.25,0){18}{\qbezier(0,3)(.05,3)(.1,3)}
\multiput(0,0)(0,.25){12}{\qbezier(4,0)(4,.05)(4,.1)}
\multiput(0,0)(0,.25){12}{\qbezier(8,0)(8,.05)(8,.1)}
\linethickness{.2mm}
\qbezier(2,0)(2.5,.05)(3,1.5)
\qbezier(3,1.5)(3.5,2.95)(4,3)
\qbezier(8,3)(8.5,2.95)(9,1.5)
\qbezier(9,1.5)(9.5,.05)(10,0)
\linethickness{.25mm}
\qbezier(4,3)(6,3)(8,3)
\end{picture}
\end{center}
\caption{\label{weight-function-for-smooth-long-sums}The weight function $w(x)$ of Lemma \ref{smooth-improved-estimates}.}
\end{figure}

\paragraph{Proof of Theorem \ref{improved-estimate}.}
Introducing the above weight function $w$ gives
\[\sum_{M\leqslant n\leqslant M+\Delta}t(n)\,e\!\left(\frac{nh}k\right)
\ll U\,M^{\vartheta+\varepsilon}+\sum_{M\leqslant n\leqslant M+\Delta}t(n)\,e\!\left(\frac{nh}k\right)w(n).\]
If we select $U=k^{2/3}\,M^{1/3-2\vartheta/3}$ and $X=k^{2/3}\,M^{1/3+4\vartheta/3}$ in Lemma \ref{smooth-improved-estimates}, we obtain
\begin{align*}
\sum_{M\leqslant n\leqslant M+\Delta}&t(n)\,e\!\left(\frac{nh}k\right)\\
&\ll U\,M^{\vartheta+\varepsilon}+k^{1/2}\,X^{1/4}\,M^{1/4}
+k^{3/2}\,X^{-1/4}\,M^{3/4}\,U^{-1}\\
&\ll k^{2/3}\,M^{1/3+\vartheta/3+\varepsilon},
\end{align*}
as required.

When $M^{3/(5+6\vartheta)-1/2+\vartheta}\ll k\ll M^{1/4+3\vartheta/8}$ we argue similarly, except that now the smoothing error is estimated by Theorem \ref{short-estimate} to be $\ll U^{1/6-\vartheta}\,M^{1/3+\vartheta+\varepsilon}$, and we choose $X=k\,M\,U^{-1}$ and $U=k^{3/(2-3\vartheta)}\,M^{(1-6\vartheta)/(4-6\vartheta)}$.
We observe that Theorem \ref{short-estimate} is applicable here since a little simplification shows that
\[U\ll\left(M^{5/18+\vartheta/3}\right)^{3/(2-3\vartheta)}\cdot M^{(1-6\vartheta)/(4-6\vartheta)}=M^{2/3}.\]
This choice of $X$ and $U$ leads to
\begin{align*}
\sum_{n\leqslant M}&t(n)\,e\!\left(\frac{nh}k\right)\\
&\ll U^{1/6-\vartheta}\,M^{1/3+\vartheta+\varepsilon}+k^{1/2}\,X^{1/4}\,M^{1/4}
+k^{3/2}\,X^{-1/4}\,M^{3/4}\,U^{-1}\\
&\ll k^{(1-6\vartheta)/(4-6\vartheta)}\,M^{3/(8-12\vartheta)+\varepsilon},
\end{align*}
as required.

When $M^{5/18+\vartheta/3}\ll k\ll M^{1/2-\varepsilon}$ we choose $X=k^2\,M\,U^{-2}$ and $U=k^{2/3}\,M^{13/27-2\vartheta/9}$, so that $U\gg M^{2/3}$, and get
\begin{align*}
\sum_{n\leqslant M}&t(n)\,e\!\left(\frac{nh}k\right)\\
&\ll U\,M^{\vartheta-2/9+\varepsilon}+k^{1/2}\,X^{1/4}\,M^{1/4}+k^{3/2}\,X^{-1}\,M^{3/4}\,U^{-1}\\
&=k^{2/3}\,M^{7/27+\vartheta/9+\varepsilon}.
\end{align*}

\paragraph{Proof of Lemma \ref{smooth-improved-estimates}.}
We shall feed the sum in question to the Voronoi type summation formula cited in Theorem \ref{full-voronoi-for-maass-waves} with the choice $f=w$. The series involving the $K$-Bessel function will be negligible: Pick any $A\in\left]2,\infty\right[$. Then the series involving the $K$-Bessel function can estimated as follows
\begin{align*}
&\ll_A\frac1k\sum_{n=1}^\infty\left|t(n)\right|\int\limits_M^{M+\Delta}w(x)\,k^A\,n^{-A/2}\,x^{-A/2}\,\mathrm dx\\
&\ll_A\frac1k\,\Delta\,(k\,M^{-1/2})^A\ll\frac1k\,\Delta\,M^{-\delta A}.
\end{align*}
For $A\gg_\delta1$, this is $\ll_\delta1$.

In the series involving the $J$-Bessel function, we apply \eqref{Voronoi-J-asymptotics} with $K=2$. The terms involving the error term contribute only
\begin{align*}
&\ll\frac1k\sum_{n=1}^\infty\left|t(n)\right|\int\limits_M^{M+\Delta}w(x)\,k^{5/2}\,n^{-5/4}\,x^{-5/4}\,\mathrm dx\\
&\ll k^{3/2}\,\Delta\,M^{-5/4}.
\end{align*}

We shall consider the series involving the $J$-function in two parts according to whether $n\leqslant X$ or $n>X$. The high-frequency terms $n>X$ are again treated by integrating by parts twice. However, here there will be a slight twist: the bound for the integral
\[\int\limits_M^{M+\Delta}w(x)\,k^{1/2}\,n^{-1/4}\,x^{-1/4}\left(1+C_\pm\,k\,n^{-1/2}\,x^{-1/2}\right)e\!\left(\pm\frac{2\sqrt{nx}}k\right)\mathrm dx\]
will be
\[\ll k^{5/2}\,n^{-5/4}\,M^{3/4}\left(M^{-2}\,\Delta+U^{-1}\right)
\ll k^{5/2}\,n^{-5/4}\,M^{3/4}\,U^{-1},\]
instead of $\ll k^{5/2}\,n^{-5/4}\,M^{3/4}\,U^{-2}\,\Delta$. The reason for this is that after having integrated by parts twice, the resulting integral is estimated by absolute values, and most of the terms in the integrands will be supported on $\mathrm{supp}\,w'$ which is a set of length $\ll U$. The only terms in which the integrand is supported in a larger set are those, which still feature $w(x)$ after differentiation, but here the other factors all give an extra $M^{-1}$ instead of mere $U^{-1}$ upon differentation.

Substituting the bound from integration by parts back into the series, we see that the contribution from the high-frequency terms is
\[\ll\frac1k\sum_{n>X}\left|t(n)\right|\,k^{5/2}\,n^{-5/4}\,M^{3/4}\,U^{-1}
\ll k^{3/2}\,X^{-1/4}\,M^{3/4}\,U^{-1}.\]

With the low-frequency terms, we estimate the integral in question by the first derivative test to get
\[\int\limits_M^{M+\Delta}\dots\mathrm dx
\ll k^{1/2}\,n^{-1/4}\,M^{-1/4}\,\frac{k\,\sqrt{M}}{\sqrt n},\]
and so the contribution from the low-frequency terms is
\begin{align*}
&\ll\frac1k\sum_{n\leqslant X}\left|t(n)\right|k^{1/2}\,n^{-1/4}\,M^{-1/4}\,k\,M^{1/2}\,n^{-1/2}
\ll k^{1/2}\,X^{1/4}\,M^{1/4}.
\end{align*}

\paragraph{Proof of Corollary \ref{holomorphic-corollary}.}
This is proved in exactly the same way as the corresponding result, Theorem 1, in \cite{Vesalainen2014}, except that Theorem \ref{holomorphic-short-estimate} above is used to estimate the smoothing error. The only change in the computations when $M^{1/4}\ll k\ll M^{5/18}$ is to observe that when $k\ll M^{5/18}$, we have $k^{3/2}\,M^{1/4}\ll M^{2/3}$. When $k\ll M^{5/18}$, we choose $U=k^{2/3}\,M^{13/27}$ instead of $U=k^{2/3}\,M^{11/24}$, and we observe that for these choices we have the required lower bound $U\gg M^{2/3}$.

\section{Proof of Theorem \ref{approximate-functional-equation}}

Let $J\in\mathbb Z_+$. In order to be able to apply the Voronoi summation formula, we shall consider the smoothed exponential sum
\begin{align}\label{smoothedsum}
\sum_{M_{-1}\leqslant n\leqslant M_2}t(n)\,w(n)\,e(\alpha n),
\end{align}
where $w$ is the weight function $\eta_J$ (see Section \ref{integral-section}) which corresponds to the interval $\left[M_{-1},M_2\right]$ with parameter $U\in\mathbb R_+$, which is defined as follows:
Let $d\in\mathbb R_+$ be a small constant depending on $\varepsilon$, and write $U=M^{1/2}\,\eta^{-1/2}\left(k^2\,\eta^2\,M\right)^d$. Let $M_{-1}=M-JU$, $M_1=M+\Delta$, and $M_2=M+\Delta+JU$. Also, we define $N_i=k^2\,\eta^2\,M_i$ for $i\in\left\{-1,1,2\right\}$ and $N=k^2\,\eta^2\,M$. The choice of $J$ depends on $d$, and a fortiori on $\varepsilon$.

\subsection{Estimating smoothing error}

First, we estimate the error caused by the introduction of the weight function $w$.
\begin{lemma}\label{painottamaton}
Let $M\in\left[1,\infty\right[$, and let $\alpha\in\mathbb R$, $h\in\mathbb Z$, $k\in\mathbb Z_+$ and $\eta\in\mathbb R$ be such that
\[\alpha=\frac hk+\eta,\quad k\leqslant M^{1/4},\quad(h,k)=1\quad\text{and}\quad
\left|\eta\right|\leqslant\frac1{k\,M^{1/4}}.\]
Furthermore, write $U=M^{1/2}\,\eta^{-1/2}\left(k^2\,\eta^2\,M\right)^d$, where $d\in\mathbb R_+$. Then, given $\varepsilon\in\mathbb R_+$, we have
\[\sum_{M\leqslant n\leqslant M+U}t(n)\,e(n\alpha)\ll M^{1/2}\left(k^2\,\eta^2\,M\right)^{\vartheta/2-1/12+\varepsilon},\]
for any fixed $d\ll_\varepsilon1$.
\end{lemma}

\noindent Now, by partial summation and Lemma $\ref{painottamaton}$ we have
\begin{align}\label{smoothingerror}
&\sum_{M_{-1}\leqslant n<M}t(n)\,e(\alpha n)\,w(n)+\sum_{M_1<n\leqslant M_2}t(n)\,e(\alpha n)\,w(n)\nonumber\\
&\qquad\ll  M^{1/2}\left(k^2\,\eta^2\,M\right)^{\vartheta/2-1/12+\varepsilon}.
\end{align}

\paragraph{Proof of Lemma \ref{painottamaton}.} Let us first dispose of the case $U\ll M^{2/3}$. In this case we have, by Theorem \ref{short-estimate},
\begin{align*}
\sum_{M\leqslant n\leqslant M+U}&t(n)\,e(n\alpha)
\ll U^{1/6-\vartheta}\,M^{1/3+\vartheta+\varepsilon}\\
&\ll\left(M^{1/2}\,\eta^{-1/2}\left(k^2\,\eta^2\,M\right)^d\right)^{1/6-\vartheta}M^{1/3+\vartheta+\varepsilon}\\
&\ll M^{1/2}\left(k^2\,\eta^2\,M\right)^{\vartheta/2-1/12+d/6-d\vartheta}\,M^\varepsilon\,k^{1/6-\vartheta}\,\eta^{1/12-\vartheta/2}.
\end{align*}
If $k^2\,\eta^2\,M\gg M^{1/4}$, then certainly
\[M^\varepsilon\,k^{1/6-\vartheta}\,\eta^{1/12-\vartheta/2}\leqslant M^\varepsilon\ll\left(k^2\,\eta^2\,M\right)^\varepsilon.\]
If $k^2\,\eta^2\,M\ll M^{1/4}$, then
\[k^{1/12-\vartheta/2}\,\eta^{1/12-\vartheta/2}\ll\left(M^{-3/8}\right)^{1/12-\vartheta/2}=M^{-1/32+3\vartheta/16},\]
so that
\begin{multline*}
M^\varepsilon\,k^{1/6-\vartheta}\,\eta^{1/12-\vartheta/2}\ll M^\varepsilon\,k^{1/12-\vartheta/2}\,M^{-1/32+3\vartheta/16}\\
\ll M^{\varepsilon+1/48-1/32-\vartheta/8+3\vartheta/16}\ll1.
\end{multline*}
Thus, in either case the sums of length $U\ll M^{2/3}$ are sufficiently small. The same argument also takes care of all the later terms which have the shape $U^{1/6-\vartheta}\,M^{1/3+\vartheta+\varepsilon}$.

Let us next focus on the case $U\gg M^{2/3}$. Let us first assume that $h/k$ is a Farey fraction of order $U^{1/2-\delta}$ for some small $\delta\in\left]0,1/2\right[$, sufficiently small depending on $\varepsilon$, i.e.\ that $\left|\eta\right|\leqslant k^{-1}\,U^{-1/2+\delta}$. Then the second error term from Proposition \ref{weight-function-removal-for-longer-sums} contributes
\begin{align*}
&\ll k^{-1/2}\,U\,M^{-1/4}\left(k^2\,\eta^2\,M\right)^{\varepsilon+\vartheta-1/4}\\
&\ll k^{-1/2}\,M^{1/2}\,\eta^{-1/2}\left(k^2\,\eta^2\,M\right)^dM^{-1/4}\left(k^2\,\eta^2\,M\right)^{\varepsilon+\vartheta-1/4}\\
&\ll M^{1/2}\left(k^2\,\eta^2\,M\right)^{d+\varepsilon+\vartheta/2-1/12}
\left(k^2\,\eta^2\,M\right)^{\vartheta/2-1/6-1/4}\\
&\ll M^{1/2}\left(k^2\,\eta^2\,M\right)^{d+\varepsilon+\vartheta/2-1/12},
\end{align*}
provided that $\eta\ll M\,U^{-2}$. But this condition holds since it reduces to
\[\eta\ll\frac1{\eta^{-1}\left(k^2\,\eta^2\,M\right)^{2d}},\]
and we have
\[\eta^2\ll\frac1{k^2\,M^{1/2}}\ll\left(k^2\,\eta^2\,M\right)^{-2d}\]
for sufficiently small $d$.

Let us observe next that if $U\gg M^{5/6}$, then $M^{5/4}\ll U^{3/2}$. Further,
\[k\,M\left(k^2\,\eta^2\,M\right)^{2d}\ll U^{3/2+\delta},\]
for sufficiently small $d$, so that, by the definition of $U$,
\[\eta=M\,U^{-2}\left(k^2\,\eta^2\,M\right)^{2d}\ll\frac1{k\,U^{1/2-\delta}}.\]
Thus, if $U\gg M^{5/6}$, then $h/k$ is indeed a Farey fraction of order $U^{1/2-\delta}$ and everything is fine.

The remaining length range is $M^{2/3}\ll U\ll M^{5/6}$, and the only problematic case is the one in which $\eta\gg k^{-1}\,U^{\delta-1/2}$. In this case we use Proposition \ref{weight-removal-for-longer-sums-with-many-farey-fractions} which involves many Farey approximations possibly different from $h/k$. Let us consider one such Farey approximation
\[\alpha=\frac{h_\ell}{k_\ell}+\eta_\ell,\]
where $\ell\in\left\{1,2,\ldots,L\right\}$, $L\in\mathbb Z_+$, $U\,2^{-L}\asymp M^{2/3}$, and $h_\ell\in\mathbb Z$, $k_\ell\in\mathbb Z_+$, $\eta_\ell\in\mathbb R_+$, $(h_\ell,k_\ell)=1$, $k_\ell\leqslant U^{1/2-\delta/2}$, and $\left|\eta\right|\leqslant k_\ell^{-1}\,(U\,4^{-\ell})^{\delta/2-1/2}$.

Let us observe that if we had $k_\ell\leqslant M^{1/4}/2$, then we would have
\[\frac1{k\,k_\ell}\leqslant\left|\frac{h_\ell}{k_\ell}-\frac hk\right|\leqslant\left|\eta\right|+\left|\eta_\ell\right|\leqslant\frac1{k\,M^{1/4}}+\frac1{k_\ell\,(U\,4^{-\ell})^{1/2-\delta/2}},\]
so that
\[1\leqslant k_\ell\,M^{-1/4}+k\left(\frac U{2^\ell}\right)^{\delta/2-1/2}\leqslant\frac12+M^{1/4+\varepsilon-1/3}=\frac12+M^{\varepsilon-1/12}=\frac12+o(1),\]
which is impossible. Thus, we must have $k_\ell\gg M^{1/4}$.

Let us now consider the term $k_\ell^{-1/2}\,U\,4^{-\ell}\,M^{-1/4}\left(k_\ell^2\,\eta_\ell^2\,M\right)^\vartheta$ which only arises when
\[k_\ell^2\,\eta_\ell\,M\,(U\,4^{-\ell})^{-1+\delta/2}\ll1\quad\text{and}\quad
k_\ell^2\,\eta_\ell^2\,M\gg1.\]
Let us check that the condition $\eta_\ell\ll M\,(U\,4^{-\ell})^{-2}$ required in this case holds. Namely, since
\[k_\ell^2\,\eta_\ell\,M\,(U\,4^{-\ell})^{-1+\delta/2}\ll1,\]
we have
\[\eta_\ell\ll\frac{(U\,4^{-\ell})^{1-\delta/2}}{k_\ell^2\,M}.\]
It is therefore enough that
\[\frac{(U\,4^{-\ell})^{1-\delta/2}}{k_\ell^2\,M}\ll\frac M{(U\,4^{-\ell})^2},\]
i.e. that $(U\,4^{-\ell})^{3-\delta/2}\ll k_\ell^2\,M^2$. Since $k_\ell\gg M^{1/4}$, it is enough that $U\,4^{-\ell}\ll M^{5/6+\delta/6}$, which is indeed true.

Next, we need to check that the term $k_\ell^{-1/2}\,U\,4^{-\ell}\,M^{-1/4}\left(k_\ell^2\,\eta_\ell^2\,M\right)^\vartheta$ is small enough. We have
\begin{align*}
&k_\ell^{-1/2}\,U\,4^{-\ell}\,M^{-1/4}\left(k_\ell^2\,\eta_\ell^2\,M\right)^\vartheta\\
&\quad\ll M^{-1/8}\,U\,4^{-\ell}\,M^{-1/4}\,(U\,4^{-\ell})^{-\vartheta+2\delta\vartheta}\,M^\vartheta\\
&\quad\ll (U\,4^{-\ell})^{1-\vartheta}\,M^{\vartheta-3/8+\varepsilon}.
\end{align*}
We have
\[(U\,4^{-\ell})^{1-\vartheta}\,M^{\vartheta-3/8}\ll (U\,4^{-\ell})^{1/6-\vartheta}\,M^{1/3+\vartheta}\]
if and only if
\[(U\,4^{-\ell})^{5/6}\ll M^{17/24},\]
or equivalently, $U\,4^{-\ell}\ll M^{17/20}$. But this holds, since $5/6\ll 17/20$.

Finally, we need to check that the condition $\eta_\ell\ll M\,(U\,4^{-\ell})^{-2}$ holds also in the case in which $k_\ell^2\,\eta_\ell^2\,M\gg1$ and $k_\ell^2\,\eta_\ell\,M\,(U\,4^{-\ell})^{-1+\delta/2}\gg1$. Since $U=M^{1/2}\,\eta^{-1/2}\left(k^2\,\eta^2\,M\right)^d$, the question is, whether
\[\eta_\ell\ll\eta\left(k^2\,\eta^2\,M\right)^{-2d}\,2^{2\ell}?\]
If this was not the case, then we could estimate
\begin{align*}
\eta_\ell&\gg\eta\left(k^2\,\eta^2\,M\right)^{-2d}2^{2\ell}
\gg k^{-1}\,U^{\delta-1/2}\left(k^2\,\eta^2\,M\right)^{-2d}2^{2\ell}\\
&\gg k_\ell^{-1}\left(\frac U{2^\ell}\right)^{\delta/2-1/2}
U^{\delta/2}\left(k^2\,\eta^2\,M\right)^{-2d}2^{3\ell/2+\ell\delta/2},\\
&\gg\eta_\ell\,M^{\delta/3-d},
\end{align*}
giving a contradiction when $d$ is sufficiently small, and we are done.

\subsection{Voronoi summation formula and saddle-points: the main terms}

Now we consider the smoothed sum (\ref{smoothedsum}). We start by estimating terms that arise when Voronoi summation formula is applied and Bessel functions are replaced by their asymptotic expressions. As always, the terms involving the $K$-Bessel function contribute a negligible amount. Asymptotics of the $J$-Bessel function lead to certain exponential integrals which are estimated by using standard tools. We will assume throughout the proof that $\eta\geqslant 0$, as the other case is similar.

The Voronoi summation formula says that
\begin{align*}
&\sum_{M_{-1}\leqslant n\leqslant M_2}t(n)\,e(n\alpha)\,w(n)\\
&=\frac{\pi i}{k\,\sinh\pi\kappa}\sum_{n=1}^\infty t(n)\,e\!\left(\frac{-n\overline h}k\right)\int\limits_{M_{-1}}^{M_2}\left(J_{2i\kappa}-J_{-2i\kappa}\right)\!\left(\frac{4\pi\sqrt{nx}}k\right)e(\eta x)\,w(x)\,\mathrm dx\\
&\qquad+\frac{4\cosh\pi\kappa}k\sum_{n=1}^\infty t(-n)\,e\!\left(\frac{n\overline h}k\right)\int\limits_{M_{-1}}^{M_2}K_{2i\kappa}\!\left(\frac{4\pi\sqrt{nx}}k\right)e(\eta x)\,w(x)\,\mathrm dx.
\end{align*}

The sum involving $K$-Bessel function contributes $\ll 1$ as before. Replacing the difference between $J$-Bessel functions by the asymptotic expression (\ref{Voronoi-J-asymptotics}) gives 

\begin{align*}
&\sum_{M_{-1}\leqslant n\leqslant M_2}t(n)\,e(n\alpha)\,w(n)\\
&=O(1)+\frac {C'}k\sum_{n=1}^\infty t(n)\,e\!\left(\frac{-n\overline h}k\right)\\
&\qquad\qquad\cdot
\int\limits_{M_{-1}}^{M_2}\frac{k^{1/2}}{n^{1/4}\,x^{1/4}}\sum_\pm (\pm1)e\!\left(\pm\frac{2\sqrt{nx}}k\right)e\!\left(\mp\frac18\right)g_\pm(x;n,k)\,e(\eta x)\,w(x)\,\mathrm dx,
\end{align*}
where
\[g_\pm(x;n,k)=1+\sum_{\ell=1}^Kc_\ell^\pm\,k^\ell\,n^{-\ell/2}\,x^{-\ell/2},\]
and $C'=i/\sqrt2$, just like in the proof of Theorem \ref{estimate-for-smooth-short-sums} as the error term from the $J$-Bessel asymptotics gives the contribution $\ll 1$ when $K$ is fixed and chosen to be large enough. Now by Lemma \ref{jutila-motohashi-lemma} we have

\[\sum_{n=1}^\infty t(n)\,n^{-1/4}\int\limits_{M_{-1}}^{M_2}\,k^{-1/2}\,x^{-1/4}\,e\!\left(x\eta+\frac{2\sqrt{nx}}k\right)g_+(x;n,k)\,w(x)\,\mathrm dx \ll 1,\]
and so the main terms come from the integrals involving $g_-$.

Let $c$ be a positive constant so that the term
\[x\eta-\frac{2\sqrt{nx}}{k}\gg\frac{\sqrt{nM}}{k}\]
when $n> cN$. A direct application of Lemma \ref{jutila-motohashi-lemma} gives

\[\sum_{n>cN} t(n)\,n^{-1/4}\int\limits_{M_{-1}}^{M_2}\,k^{-1/2}\,x^{-1/4}\,e\!\left(x\eta-\frac{2\sqrt{nx}}k\right)g_-(x;n,k)\,w(x)\,\mathrm dx \ll 1.\]

\noindent For the terms with $n\leqslant cN$ we split $g_-(x;n,k)$ into two parts $1$ and $g_-(x;n,k)-1$ and estimate corresponding terms differently. 

For the first term, using the second derivative test we get

\[ \int\limits_{M_{-1}}^{M_2}\,x^{-3/4}\,e\!\left(x\eta-\frac{2\sqrt{nx}}k\right)w(x)\,\mathrm dx\ll\frac{k^{1/2}}{n^{1/4}}.\]
 Therefore we have
\begin{align*}
&\sum_{n\leqslant cN} t(n)\,n^{-1/4}\int\limits_{M_{-1}}^{M_2}\,k^{-1/2}\,x^{-1/4}\,e\!\left(x\eta-\frac{2\sqrt{nx}}k\right)(g_-(x;n,k)-1)\,w(x)\,\mathrm dx \\
&\qquad\quad
\ll k\sum_{n\leqslant cN}\frac{\left|t(n)\right|}n\ll k\,N^{\varepsilon}=k\,(k^2\,\eta^2\,M)^{\varepsilon}.
\end{align*}
The remaining terms are treated using the first saddle point lemma, Theorem $\ref{saddle-point-lemma}$. 
For $1\leqslant n<cN$ we get
\begin{align*}
&\int\limits_{M_{-1}}^{M_2} e\!\left(x\eta-\frac{2\sqrt{nx}}k\right)w(x)\,x^{-1/4}\,\mathrm dx\\
&\qquad=\xi\!\left({n}\right)\cdot\frac{\sqrt 2\, n^{1/4}}{\sqrt k\, \eta}\,e\!\left(-\frac{n}{k^2\eta}+\frac18\right)\\
&\qquad\qquad+O\!\left((M_2-M_{-1})\left(1+M^{J}\,U^{-J}\right)M^{-1/4}\,e^{-A|\eta|M-A\sqrt{nM}/k}\right)\\
&\qquad\qquad+O\!\left(\frac{k^{3/2}}{n^{3/4}}+\chi(n)\,\frac{M^{1/4}k}{\sqrt n}\right)\\
&\qquad\qquad+O\!\left(M^{-1/4}\,U^{-J}\sum_{j=0}^J\left(\left|\eta-\frac{\sqrt n}{k\,\sqrt{M+jU}}\right|+\frac{n^{1/4}}{\sqrt k\, M^{3/4}}\right)^{-J-1}\right)\\
&\qquad\qquad+O\!\left(M^{-1/4}\,U^{-J}\sum_{j=0}^J\left(\left|\eta-\frac{\sqrt n}{k\,\sqrt{M+\Delta+jU}}\right|+\frac{n^{1/4}}{\sqrt k \,M^{3/4}}\right)^{-J-1}\right),
\end{align*} 
where we have written for simplicity $\xi(\cdot)=\xi_J(\cdot/(k^2\eta^2))$, and where
\begin{align*}
\left\{ \!\!
  \begin{array}{l l l}
    \xi(n)=0 \text{ and }\chi(n)=0 & \text{ if $n\leqslant N_{-1}$ or $n\geqslant N_2$,}\\
    \xi(n)=1 \text{ and }\chi(n)=0 & \text{ if $N\leqslant n\leqslant N_1$,}\\
    \xi(n)\ll 1
    \text{ and }\chi(n)=1 &  \text{ otherwise.}\\
  \end{array} \right.
\end{align*}
Furthermore, $\xi'$ is piecewise continuously differentiable and $\xi'(n)\ll(k^2\eta^2 U)^{-1}$ where the derivative exists.

\noindent The main term on the right-hand side produces the total contribution 
\begin{align*}
\frac{1}{k\eta}\sum_{N\leqslant n\leqslant N_1} t(n)\,e\!\left(-\frac{n\overline h}k-\frac{n}{k^2\eta}\right)&+\frac{1}{k\eta}\sum_{N_{-1}\leqslant n< N} t(n)\,\xi(n)\,e\!\left(-\frac{n\overline h}k-\frac{n}{k^2\eta}\right)\\
&\quad+\frac{1}{k\eta}\sum_{N_1<n\leqslant N_2} t(n)\,\xi(n)\,e\!\left(-\frac{n\overline h}k-\frac{n}{k^2\eta}\right).
\end{align*}
The first term is exactly what appears in the statement of the theorem. Let us first estimate the contribution of error terms arising from the saddle point lemma and after that estimate the contribution of other main terms.

\subsection{The error terms from the saddle point theorem}

The first error term contributes

\begin{align*}
&k^{-1/2}\sum_{1\leqslant n\leqslant cN}\frac{\left|t(n)\right|}{n^{3/4}}\,(\Delta+U)\left(1+M^{J}U^{-J}\right)M^{-1/4}\,\exp\!\left(-A|\eta|M^{1/4}-\frac{A\sqrt{nM}}k\right)\\
&\qquad\ll (\Delta+U)\,M^{J-1/4}\,N^{1/4}\cdot e^{-AM^{1/4}}\ll 1.
\end{align*}

The second error term is also easy to handle:

\begin{align*}
k^{-1/2}\sum_{1\leqslant n\leqslant cN}&\frac{\left|t(n)\right|}{n^{1/4}}\left(\frac{k^{3/2}}{n^{3/4}}+\chi(n)\frac{M^{1/4}k}{n^{1/2}}\right)\\
&\qquad\ll k\,(k^2\eta^2M)^{\varepsilon+d+\vartheta}\ll M^{1/2}\,(k^2\eta^2 M)^{d-1/12+\varepsilon}
\end{align*}
just by using partial summation.

The estimation of the other two error terms is covered by the following lemma.

\begin{lemma} Let $c$ be any given positive constant. Let $T$ be any number of type $M\pm jU$, where $0\leqslant j\leqslant J$. Then
\begin{multline*}
k^{-1/2}\sum_{1\leqslant n\leqslant cN}\frac{t(n)}{n^{1/4}}\,M^{-1/4}\,U^{-J}\left(\left|\eta-\frac{\sqrt n}{\sqrt T \,k}\right|+\frac{n^{1/4}}{\sqrt k\,M^{3/4}}\right)^{-J-1}\\
\ll\sqrt M\,(k^2\eta^2M)^{1/2-Jd}.
\end{multline*}
\end{lemma}

\paragraph{Proof.}
We estimate
the left-hand side as $\ll S_1+S_2+S_3$, where
\begin{align*}
S_1=k^{-1/2}\,M^{-1/4}\,U^{-J}\sum_{|n-k^2\eta^2T|\leqslant \sqrt N}n^{\vartheta+\varepsilon-1/4}\left(\frac{n^{1/4}}{\sqrt k M^{3/4}}\right)^{-J-1},
\end{align*} 
\begin{align*}
S_2=k^{-1/2}\,M^{-1/4}\,U^{-J}\sum_{1\leqslant n\leqslant k^2\eta^2T-\sqrt N}n^{\vartheta+\varepsilon-1/4}\left|\eta-\frac{\sqrt n}{\sqrt T k}\right|^{-J-1},
\end{align*}
and
\begin{align*}
S_3=k^{-1/2}\,M^{-1/4}\,U^{-J}\sum_{k^2\eta^2T+\sqrt N\leqslant n\leqslant cN}n^{\vartheta+\varepsilon-1/4}\left|\eta-\frac{\sqrt n}{\sqrt T k}\right|^{-J-1}.
\end{align*}
Next we compute the claimed upper bound for each of them. Observe that by partial summation
\begin{align*}
S_1&=M^{1/2}(k^2\eta^2M)^{-Jd}k^{J/2}M^{-J/4}\eta^{J/2}\sum_{|n-k^2\eta^2T|\leqslant \sqrt N}n^{\vartheta+\varepsilon-1/2-J/4}\\
&\ll M^{1/2}(k^2\eta^2M)^{-Jd}k^{J/2}M^{-J/4}\eta^{J/2}(k^2\eta^2M)^{\vartheta+\varepsilon+1/2-J/4}\\
&\ll M^{1/2}(k^2\,\eta^2M)^{\varepsilon-Jd+\vartheta}.
\end{align*}
The sum $S_2$ is estimated as follows:
\begin{align*}
&k^{-1/2}\,M^{1/4}\,M^{-J/2}\,\eta^{J/2}\,N^{-dJ}\sum_{n\leqslant k^2\eta^2T-\sqrt N}\frac{\left|t(n)\right|}{n^{1/4}}\left|\eta-\frac{\sqrt n}{k\sqrt T}\right|^{-J-1}\\
&\ll k^{-1/2}\,M^{1/4}\,M^{-J/2}\,\eta^{J/2}\,N^{-dJ}\,k^{J+1}\,M^{J/2+1/2}\\
&\qquad\cdot\sum_{n\leqslant k^2\eta^2T-\sqrt N}\frac{\left|t(n)\right|}{n^{1/4}}\cdot\frac{\left|k\eta\sqrt T+\sqrt n\right|^{J+1}}{\left|k^2\eta^2T-n\right|^{J+1}}\\
&\ll k^{-1/2}\,M^{-1/4}\,M^{-J/2}\,\eta^{J/2}\,N^{-dJ}\,k^{J+1}\,M^{J/2+1/2}\,N^{-J/2+1/2}\,N^{J/2+1/2}\,N^{3/4}\\
&\ll k^{1/2+J}\,\eta^{J/2}\,M^{1/4}\,N^{-dJ+3/4}
\ll M^{3/8}\,N^{3/4-dJ}\ll M^{1/2}\,N^{1/2-dJ}.
\end{align*}
Finally, the sum $S_3$ is estimated in the same manner as $S_2$.

\subsection{Removing the weight function $\xi$}

 By partial summation it is enough to deal with the sum without $\xi(n)$. Observe that $N-N_{-1}=N_2-N_1=k^2\,\eta^2\,J\,U$ and
\[k^2\,\eta^2\,U=(k^2\,\eta^2\,M)^{d+1/2}\,k\,\eta^{1/2}\ll(k^2\,\eta^2\,M)^{2/3}=N^{2/3}\]
for sufficiently small $d\in\mathbb R_+$.

Therefore, using Theorem $\ref{short-estimate}$ we get that other two main terms given by the main term of the saddle point lemma contribute
\begin{align*}
&\frac{1}{k\eta}\left(\sum_{N_{-1}\leqslant n<N}+\sum_{N_1<n\leqslant N_2} \right) t(n)\,\xi(n)\,e\!\left(-\frac{n\overline h}k-\frac{n}{k^2\eta}\right)\\
&\ll\frac1{k\eta}\left(k^2\,\eta^2\left(M^{1/2}\,\eta^{-1/2}\,k^2\,\eta^2\,M\right)^d\right)^{1/6-\vartheta}\cdot(k^2\,\eta^2\,M)^{1/3+\vartheta+\varepsilon}\\
&\ll(k^2\,\eta^2\,M)^{d(1/6-\vartheta)+\varepsilon} M^{1/2-\vartheta/2}\\
&\ll M^{1/2}(k^2\,\eta^2 \,M)^{d/6+\left(1/2-d\right)\vartheta-1/12} \\
&\ll M^{1/2}(k^2\,\eta^2\, M)^{d/6+\vartheta/2-1/12}
\end{align*}
for small enough $d\in\mathbb R_+$.

At this point we have proved that
\begin{align*}
\sum_{M_{-1}\leqslant n\leqslant M_2}t(n)\,w(n)\,e(\alpha n)=&\frac1{k\eta}\sum_{N\leqslant n\leqslant N_1}t(n)\,e(-\beta n)+O(M^{1/2}\,(k^2\eta^2M)^{1/2-Jd})\\
&\quad\quad+O(M^{1/2}\,(k^2\,\eta^2\,M)^{\varepsilon+d/6-1/12+\vartheta/2})\\
&\qquad\quad+O\left(k\,(k^2\,\eta^2\,M)^{\varepsilon}\right).
\end{align*}
Furthermore, using \eqref{smoothingerror}, this tells that

\begin{align*}
\sum_{M\leqslant n\leqslant M+\Delta}t(n)\,e(\alpha n)=&\frac1{k\eta}\sum_{N\leqslant n\leqslant N_1}t(n)\,e(-\beta n)+O(M^{1/2}\,(k^2\eta^2M)^{1/2-Jd})\\
&\quad\quad+O\left(M^{1/2}\,(k^2\,\eta^2\,M)^{\varepsilon+d/6-1/12+\vartheta/2}\right)\\
&\qquad\qquad\quad+O\left(k\,(k^2\eta^2M)^{\varepsilon}\right).
\end{align*}
Choosing $J$ sufficiently large depending on $d$, and letting $d\in\mathbb R_+$ to be arbitrarily small finishes the proof.

\section{Proof of Theorem \ref{logarithm-removal}}

We shall prove theorem \ref{logarithm-removal} first near rational points, and then iterate approximate functional equation in the remaining cases until either we end up near a rational point or the sum in question has become shorter than some given constant.

\subsection{Logarithm removal near rational points}

The following lemma, which we will soon prove, covers the logarithm removal near rational points.

\begin{lemma}
\label{logarithm-removal-near-rational-points}
Let $M\in\left[1,\infty\right[$, let $\alpha\in\mathbb R$, and $h\in\mathbb Z$ and $k\in\mathbb Z$ be coprime with $1\leqslant k\leqslant M^{1/4}$, and $\alpha=h/k+\eta$ with $\left|\eta\right|\leqslant k^{-1}\,M^{-1/4}$. If $k^2\,\eta^2\,M<1/2$, then
\[\sum_{M\leqslant n\leqslant2M}t(n)\,e(n\alpha)\ll
k^{(1-6\vartheta)/(4-6\vartheta)}\,M^{3/(8-12\vartheta)+\varepsilon}.\]
In particular, with the exponent $\vartheta=7/64$ the upper bound is $\ll M^{203/428+\varepsilon}\ll M^{1/2}$.
\end{lemma}

The following Voronoi type identity for Maass forms can be found in Section 12 of Meurman's paper \cite{Meurman1988}.
\begin{theorem}\label{voronoi-identity-of-meurman}
Let $x\in\left[1,\infty\right[$, and let $h$ and $k$ be coprime integers with $k\geqslant1$. Then
\begin{align*}
\sum_{n\leqslant x}'t(n)\,e\!\left(\frac{nh}k\right)
&=\frac{2\pi}{k\,\sinh\pi\kappa}
\sum_{n=1}^\infty
t(n)\,e\!\left(\frac{-n\overline h}k\right)
\int\limits_0^x
\Re\!\left(i\,J_{2i\kappa}\!\left(\frac{4\pi\sqrt{nv}}k\right)\right)\mathrm dv\\
&\qquad+\frac{4\cosh\pi\kappa}k\sum_{n=1}^\infty
t(-n)\,e\!\left(\frac{n\overline h}k\right)
\int\limits_0^x
K_{2i\kappa}\!\left(\frac{4\pi\sqrt{nv}}k\right)
\mathrm dv,
\end{align*}
and where the series are boundedly convergent for $x$ restricted in any bounded subinterval of $\left[1,\infty\right[$.
\end{theorem}
\noindent The integrals involving the $J$-Bessel function will have an asymptotic expansion reminiscent of those for the $J$-Bessel function itself. The following asymptotics for the $J$-Bessel function integral are obtained from Section 6 of \cite{Meurman1988}, and the asymptotics for the integral involving the $K$-Bessel function is easily obtained from the asymptotic properties of $K_{2i\kappa}$.
\begin{lemma}\label{Bessel-integral-asymptotics}
Let $n\in\mathbb Z_+$, $x\in\left[1,\infty\right[$, and $k\in\mathbb Z_+$.
If $n\,x\gg k^2$, then we have an asymptotic expansion
\begin{align*}
&\int\limits_0^x\Re\left(i\,J_{2i\kappa}\!\left(\frac{4\pi\sqrt{nv}}k\right)\right)\mathrm dv\\
&\qquad=k^{3/2}\,n^{-3/4}\,x^{1/4}\sum_{\pm}A_{1,\pm}\,e\!\left(\pm\frac{2\sqrt{nx}}k\right)
+A_2\,k^2\,n^{-1}\\
&\qquad\quad+k^{5/2}\,n^{-5/4}\,x^{-1/4}\sum_{\pm}A_{3,\pm}\,e\!\left(\pm\frac{2\sqrt{nx}}k\right)
+O_\kappa(k^{7/2}\,n^{-7/4}\,x^{-3/4}),
\end{align*}
where $A_{1,+}$, $A_{1,-}$, $A_2$, $A_{3,+}$ and $A_{3,-}$ are some constants only depending on $\kappa$, and the implicit constant in the lower bound $n\,x\gg k^2$. Similarly, we have the asymptotic expansion
\begin{align*}
\int\limits_0^xK_{2i\kappa}\!\left(\frac{4\pi\sqrt{nv}}k\right)\mathrm dv
=B_2\,k^2\,n^{-1}+O_{\kappa,C}(k^{2+C}\,n^{-1-C/2}\,x^{-C/2}),
\end{align*}
where $C\in\mathbb R_+$ is arbitrary and $B_2$ is a constant only depending on $\kappa$ and the implicit constant in $n\,x\gg k^2$.
\end{lemma}

We need one more lemma before the proof as the special value $L(1,h/k)$ will appear there.
\begin{lemma}\label{twisted-L-function-at-1}
Let $h$ and $k$ be coprime integers with $k\geqslant1$. Then
\[\sum_{n=1}^\infty\frac{t(n)}n\,e\!\left(\frac{nh}k\right)\ll k^\varepsilon.\]
\end{lemma}

\paragraph{Proof.} It is proved in \cite{Meurman1988} that the rationally additively twisted $L$-function attached to our fixed Maass form,
\[L\!\left(s,\frac hk\right)=\sum_{n=1}^\infty\frac{t(n)}{n^s}\,e\!\left(\frac{nh}k\right),\]
at first defined only for complex numbers $s$ with $\Re s>1$, has an entire analytic extension to $\mathbb C$. Furthermore, this $L$-function, in a sense, satisfies a functional equation with $\Gamma$-factors. Using the fact that the twisted $L$-functions are $\ll_\delta1$ on the vertical line $\Re s=1+\delta$ for any fixed $\delta\in\mathbb R_+$, the functional equations combined with Stirling's formula easily give the bound
\[L\!\left(s,\frac hk\right)\ll_\delta k^{1+2\delta}\left(1+\left|t\right|\right)^{1+2\delta}\quad\text{on the vertical line}\quad\Re s=-\delta.\]
Phragm\'en--Lindel\"of principle then tells us that
\[L\!\left(s,\frac hk\right)\ll_\delta k^{1+\delta-\sigma}\left(1+\left|t\right|\right)^{1+\delta-\sigma}\]
in the vertical strip $-\delta\leqslant\Re s\leqslant1+\delta$. Applying this with $s=1$ gives the Lemma. For more details about the functional equations used here, we refer to Section 2 of \cite{Meurman1988}.

\paragraph{Proof of Lemma \ref{logarithm-removal-near-rational-points}.} We begin by integrating by parts:
\begin{multline*}
\sum_{M\leqslant n\leqslant2M}t(n)\,e(n\alpha)
=\left.e(\eta x)\sum_{n\leqslant x}t(n)\,e\!\left(\frac{nh}k\right)\right]_M^{x=2M}\\
-2\pi i\eta\int\limits_M^{2M}e(\eta x)\sum_{n\leqslant x}t(n)\,e\!\left(\frac{nh}k\right)\,\mathrm dx.
\end{multline*}
As $k\leqslant M^{1/4}$, Theorem \ref{improved-estimate} immediately tells us that the substitution terms are $\ll k^{(1-6\vartheta)/(4-\vartheta)}\,M^{3/(8-12\vartheta)+\varepsilon}$. We will prove that the term involving the integral is actually $\ll k^{1/2}\,M^{1/4}$. The full Voronoi identity for Maass forms tells us that
\begin{align*}
&\eta\int\limits_M^{2M}e(\eta x)\sum_{n\leqslant x}t(n)\,e\!\left(\frac{nh}k\right)\mathrm dx\\
&\qquad=\frac{2\pi\,\eta}{k\,\sinh\pi\kappa}\sum_{n=1}^\infty t(n)\,e\!\left(\frac{-n\overline h}k\right)\int\limits_M^{2M}e(\eta x)\int\limits_0^x\Re\!\left(i\,J_{2i\kappa}\!\left(\frac{4\pi\sqrt{nv}}k\right)\right)\mathrm dv\,\mathrm dx\\
&\qquad\quad+\frac{4\eta\,\cosh\pi\kappa}k\sum_{n=1}^\infty t(-n)\,e\!\left(\frac{n\overline h}k\right)\int\limits_M^{2M}e(\eta x)\int\limits_0^xK_{2i\kappa}\!\left(\frac{4\pi\sqrt{nv}}k\right)\mathrm dv\,\mathrm dx.
\end{align*}
We emphasize that termwise integration of the series is allowed since the series converge boundedly. We note that the integral $\int_0^x$ involving the $K$-Bessel function has better asymptotic behaviour than the similar integral involving the $J$-Bessel function, and since the two series otherwise have largely the same shape, it is enough to consider the series involving $J_{2i\kappa}$.

Next we simply replace the $\int_0^x\Re(i\dots)\mathrm dv$ by the asymptotics given by Lemma \ref{Bessel-integral-asymptotics}. We start with the contribution from either of the first main terms. Since $k^2\,\eta^2\,M<1/2$, we have
\[\frac{\mathrm d}{\mathrm dx}\left(\pm\frac{2\sqrt{nx}}k+\eta\right)
=\pm\frac{\sqrt n}{k\,\sqrt x}+\eta\asymp\frac{\sqrt n}{k\,\sqrt x}\asymp n^{1/2}\,k^{-1}\,M^{-1/2}.\]
Thus, using the first derivative test, the contribution from these terms is
\begin{align*}
&\ll\frac\eta k\sum_{n=1}^\infty\left|t(n)\right|\frac{k^2}n\int\limits_M^{2M}e(\eta x)\,\frac{n^{1/4}\,x^{1/4}}{k^{1/2}}\,e\!\left(\pm\frac{2\sqrt{nx}}k\right)\mathrm dx\\
&\ll\frac1{k^2\,M^{1/2}}\sum_{n=1}^\infty\left|t(n)\right|\frac{k^2}n\cdot\frac{n^{1/4}\,M^{1/4}}{k^{1/2}}\cdot\frac{k\,M^{1/2}}{n^{1/2}}
\ll k^{1/2}\,M^{1/4}.
\end{align*}
The contribution from the constant term of the asymptotics is
\begin{align*}
&\ll\frac\eta k\sum_{n=1}^\infty t(n)\,e\!\left(\frac{-n\overline h}k\right)
\frac{k^2}n\int\limits_M^{2M}e(\eta t)\,\mathrm dt\\
&\ll k\sum_{n=1}^\infty\frac{t(n)}n\,e\!\left(\frac{-n\overline h}k\right)
\ll k^{1+\varepsilon}\ll k^{1/2}\,M^{1/8+\varepsilon}.
\end{align*}
The contribution from the third main terms is clearly smaller than that from the first main terms since $k\,n^{-1/2}\,x^{-1/4}\ll1$. Finally, the contribution from the $O$-term of the asymptotics contributes
\begin{align*}
&\ll\frac\eta k\sum_{n=1}^\infty\left|t(n)\right|k^{7/2}\,n^{-7/4}\int\limits_M^{2M}x^{-3/4}\,\mathrm dx\\
&\ll k^{-1}\,M^{-1/2}\,k^{5/2}\,M^{1/4}=k^{3/2}\,M^{-1/4}\ll k^{1/2},
\end{align*}
and we are done.

\subsection{Away from rational points; applying the approximate functional equation}

When $k^2\,\eta^2\,M\gg1$, the logarithm removal is implemented quite easily using the approximate functional equation. The result will be as follows:
\begin{lemma}\label{logarithm-removal-away-from-rational-points}
Let $M\in\left[1,\infty\right[$, let $\alpha\in\mathbb R$, let $h$ and $k$ be coprime integers with $1\leqslant k\leqslant M^{1/4}$, and let $\alpha=h/k+\eta$ with $\left|\eta\right|\leqslant k^{-1}\,M^{-1/4}$. If $k^2\,\eta^2\,M\gg1$, then
\[\sum_{M\leqslant n\leqslant2M}t(n)\,e(n\alpha)\ll M^{1/2}.\]
\end{lemma}

We start by applying the approximate functional equation, obtaining:
\begin{multline*}
\frac1{M^{1/2}}\sum_{M\leqslant n\leqslant2M}t(n)\,e(n\alpha)
=\frac1{(k^2\,\eta^2\,M)^{1/2}}\sum_{k^2\eta^2M\leqslant n\leqslant2k^2\eta^2M}t(n)\,e(n\beta)\\
+O\big((k^2\,\eta^2\,M)^{\vartheta/2-1/12+\varepsilon}\big),
\end{multline*}
where $\beta=-\overline h/k-(k^2\,\eta)^{-1}$. Write for $\beta$ a rational approximation $\beta=h_1/k_1+\eta_1$ with $h_1$ and $k_1$ coprime and $1\leqslant k_1\leqslant(k^2\,\eta^2\,M)^{1/4}$ and with remainder satisfying $\left|\eta_1\right|\leqslant k_1^{-1}\,(k^2\,\eta^2\,M)^{-1/4}$. If $k_1^2\,\eta_1^2\,(k^2\,\eta^2\,M)<1/2$, then the first term on the right-hand side is $\ll1$ by Lemma \ref{logarithm-removal-near-rational-points}, and the error term is clearly $\ll1$, and we are done.

If instead $k_1^2\,\eta_1^2\,(k^2\,\eta^2\,M)\gg1$, then we apply the approximate functional equation again to the right-hand side, and iterate the above argument as many times as necessary. Since the length of the new sum from the approximate functional equation is at most the square root of the length of the previous sum, the exponential sum term will eventually be covered by Lemma \ref{logarithm-removal-near-rational-points} or become shorter than some constant length. In either case, the sum will ultimately be $\ll1$, and the error terms will form a nice geometric progression which sums up to $\ll1$.

\section{Proof of Theorem \ref{longer-short-estimates}}

Let us begin with a simple corollary to Theorems \ref{approximate-functional-equation} and \ref{short-estimate}.

\begin{lemma}\label{corollary-to-approximate-functional-equation}
Let $M\in\left[1,\infty\right[$ and $\Delta\in\left[1,M\right]$ with $M\gg1$ and $M^{2/3}\ll\Delta\ll M^{3/4}$. Furthermore, let $\alpha\in\mathbb R$, $h\in\mathbb Z$, $k\in\mathbb Z_+$ and $\eta\in\mathbb R$ with $h$ and $k$ coprime, $\alpha=h/k+\eta$, $k\leqslant M^{1/4}$ and $\left|\eta\right|\leqslant k^{-1}\,M^{-1/4}$. Also, let $k^2\,\eta^2\,M=M^\gamma$ with $\gamma\in\mathbb R_+$. Then
\[\sum_{M\leqslant n\leqslant M+\Delta}t(n)\,e(n\alpha)
\ll\Delta^{1/6-\vartheta}\,M^{1/3+\vartheta+\varepsilon}
+M^{1/2-(1/12-\vartheta/2)\gamma+\varepsilon}.\]
\end{lemma}

\paragraph{Proof.} The approximate functional equation of Theorem \ref{approximate-functional-equation} immediately gives
\begin{multline*}
\sum_{M\leqslant n\leqslant M+\Delta}t(n)\,e(n\alpha)\\
=\frac1{k\eta}\sum_{k^2\eta^2M\leqslant n\leqslant k^2\eta^2(M+\Delta)}t(n)\,e(n\beta)+O(M^{1/2-(1/12-\vartheta/2)\gamma+\varepsilon}),
\end{multline*}
where $\beta=-\overline h/k-1/(k^2\,\eta)$.
We also have $k^2\,\eta^2\,\Delta\ll\left(k^2\,\eta^2\,M\right)^{2/3}$ since this is equivalent to
\[\left(k^2\,\eta^2\right)^{1/3}\Delta\ll M^{2/3},\]
which is true in view of
\[\left(k^2\,\eta^2\right)^{1/3}\Delta\ll M^{-1/6}\,M^{3/4}=M^{7/12}\ll M^{2/3}.\]
Thus, we may use Theorem \ref{short-estimate} to estimate
\begin{multline*}
\sum_{k^2\eta^2M\ll n\ll k^2\eta^2(M+\Delta)}t(n)\,e(n\beta)\\\ll
\left(k^2\,\eta^2\,\Delta\right)^{1/6-\vartheta}\,\left(k^2\,\eta^2\,M\right)^{1/3+\vartheta+\varepsilon}\ll k\,\eta\,\Delta^{1/6-\vartheta}\,M^{1/3+\vartheta+\varepsilon},
\end{multline*}
and the lemma has been proved.
\bigbreak
Next, we shall prove from the Voronoi summation formula an another estimate for the same sum which works nicely for a different range of $k^2\,\eta^2\,M$.

\begin{lemma}\label{last-application-of-voronoi}
Let $M,\Delta\in\mathbb R_+$ with $M\gg1$ and $M^{\beta}\ll\Delta\ll M^{3/4}$, where $\beta\in\left[2/3,3/4\right]$. Furthermore, let $\alpha\in\mathbb R$, $h,k\in\mathbb Z$ and $\eta\in\mathbb R$ with $h$ and $k$ coprime, $1\leqslant k\leqslant M^{1/4}$, $\alpha=h/k+\eta$ and $\left|\eta\right|\leqslant k^{-1}\,M^{-1/4}$. Also, let $w\in C_{\mathrm c}^\infty(\mathbb R_+)$ be supported on $\left[M,M+\Delta\right]$, take values only from $\left[0,1\right]$ and satisfy $w^{(\nu)}(x)\ll_\nu\Delta^{-\nu}$ for all $\nu\in\mathbb R_+$ and $\nu\in\left\{0\right\}\cup\mathbb Z_+$.

If $k^2\,\eta^2\,M\leqslant1/2$, then
\[\sum_{M\leqslant n\leqslant M+\Delta}t(n)\,e(n\alpha)\,w(n)\ll M^{9/8-\beta}.\]

If $k^2\,\eta^2\,M\gg1$, then, for any $S\in\left[1,\infty\right[$,
\begin{multline*}
\sum_{M\leqslant n\leqslant M+\Delta}t(n)\,e(n\alpha)\,w(n)
\ll k^3\left|\eta\right|^{3/2}M^{3/2}\,\Delta^{-1}\left(k^2\,\eta^2\,M\right)^{\vartheta+\varepsilon}\,S^{-1}\\
+k^{-1/2}\,\Delta\,M^{-1/4}\left(S+k^2\,\eta^2\,\Delta\right)\left(k^2\,\eta^2\,M\right)^{-1/4+\vartheta+\varepsilon}+M^{9/8-\beta}.
\end{multline*}
\end{lemma}

\paragraph{Proof.} We shall apply the full Voronoi summation formula for one last time and write as before
\begin{align*}
&\sum_{M\leqslant n\leqslant M+\Delta}t(n)\,e(n\alpha)\,w(n)\\
&=\frac{\pi i}{k\sinh\pi\kappa}\sum_{n=1}^\infty t(n)\,e\!\left(\frac{-n\overline h}k\right)\int\limits_M^{M+\Delta}\left(J_{2i\kappa}-J_{-2i\kappa}\right)\!\left(\frac{4\pi\sqrt{nx}}k\right)e\!\left(\eta x\right)w\!\left(x\right)\mathrm dx\\
&\qquad+\frac{4\cosh\pi\kappa}k\sum_{n=1}^\infty t(-n)\,e\!\left(\frac{n\overline h}k\right)\int\limits_M^{M+\Delta}K_{2i\kappa}\!\left(\frac{4\pi\sqrt{nx}}k\right)e(\eta x)\,w(x)\,\mathrm dx.
\end{align*}
As in the proof of Theorem \ref{estimate-for-smooth-short-sums}, the $K$-terms give only the small contribution $\ll1$. Also, applying the asymptotics \eqref{Voronoi-J-asymptotics} with a sufficiently large  fixed $K\in\mathbb Z_+$, the error terms also give $\ll1$. Thus, we have
\begin{align*}
&\sum_{M\leqslant n\leqslant M+\Delta}t(n)\,e(n\alpha)\,w(n)
=k^{-1/2}\sum_{n=1}^\infty t(n)\,n^{-1/4}\,e\!\left(\frac{-n\overline h}k\right)\\
&\qquad\cdot\int\limits_M^{M+\Delta}x^{-1/4}\sum_\pm A_\pm\,e\!\left(\pm\frac{2\sqrt{nx}}k\right)e\!\left(\pm\frac18\right)g_\pm(x;n,k)\,e(\eta x)\,w(x)\,\mathrm dx+O(1),
\end{align*}
where $A_\pm$ are constants, and
\[g_\pm(x;n,k)=1+\sum_{\ell=1}^Kc_\ell^\pm\,k^\ell\,n^{-\ell/2}\,x^{-\ell/2}.\]

If $k^2\,\eta^2\,M\leqslant1/2$, then we may simply estimate using Lemma \ref{jutila-motohashi-lemma} with $P=2$ that the infinite series is
\begin{align*}
&\ll k^{-1/2}\sum_{n=1}^\infty\left|t(n)\right|n^{-1/4}\,M^{-1/4}\left(\Delta\,n^{1/2}\,k^{-1}\,M^{-1/2}\right)^{-2}\Delta\\
&\ll k^{3/2}\,M^{3/4}\,\Delta^{-1}
\ll M^{3/8}\,M^{-\beta}\,M^{3/4}\ll M^{9/8-\beta}.
\end{align*}
Thus, we may focus on the case $k^2\,\eta^2\,M\gg1$.

Next, writing $X=k^2\,\eta^2\,M$, the high-frequence terms $n>2X$ contribute, again using Lemma \ref{jutila-motohashi-lemma} with $P=2$,
\begin{align*}
&\ll k^{-1/2}\sum_{n>2X}\left|t(n)\right|n^{-1/4}\,M^{-1/4}\left(\Delta\,n^{1/2}\,k^{-1}\,M^{-1/2}\right)^{-2}\Delta\\
&\ll k^{3/2}\,M^{3/4}\,\Delta^{-1}\,X^{-1/4}\ll M^{3/8}\,M^{3/4}\,M^{-\beta}\,M^{-\gamma/4}\ll M^{9/8-\beta-\gamma/4}.
\end{align*}
And so we are left with
\begin{align*}
&\sum_{M\leqslant n\leqslant M+\Delta}t(n)\,e(n\alpha)\,w(n)
=k^{-1/2}\sum_\pm A_\pm\sum_{n\leqslant2X}t(n)\,n^{-1/4}\,e\!\left(\frac{-n\overline h}k\right)\\
&\quad\cdot\int\limits_M^{M+\Delta}x^{-1/4}\,g_\pm(x;n,k)\,e\!\left(\pm\frac{2\sqrt{nx}}k+\eta x\right)w(x)\,\mathrm dx+O(M^{9/8-\beta-\gamma/4})+O(1).
\end{align*}

We shall split the sum $\sum_{n\leqslant2X}$ into three parts
\[\sum_{n\leqslant2X}=\sum_{n<X-S}+\sum_{X-S\leqslant n\leqslant X'+S}+\sum_{X'+S<n\leqslant2X},\]
where for simplicity $X'=k^2\,\eta^2\left(M+\Delta\right)$ and $S$ is a parameter to be chosen later but which satisfies $S\gg1$. The first and third sums might be empty; this happens when $S\gg k^2\,\eta^2\,M$. Also, large values of $S$ pose no problems in the middle terms as they are estimated via absolute values. Let us first consider the case $S\ll k^2\,\eta^2\,M$.

The third sum is estimated using Lemma \ref{jutila-motohashi-lemma} with $P=2$ to get
\begin{align*}
&k^{-1/2}\sum_{X'+S<n\leqslant2X}\dots\\
&\ll k^{-1/2}\sum_{X'+S<n\leqslant2X}\left|t(n)\right|n^{-1/4}\cdot M^{-1/4}\,\Delta^{-P}\left(\frac{\sqrt n}{k\sqrt{M+\Delta}}-\left|\eta\right|\right)^{\!-P}\Delta\\
&\ll k^{-1/2}\,M^{-1/4}\,\Delta^{-1}\,X^{\vartheta-1/4+\varepsilon}\,k^2\,M\sum_{X'+S<n\leqslant2X}\frac1{(\sqrt n-\sqrt{X'})^2}\\
&\ll k^{3/2}\,M^{3/4}\,\Delta^{-1}\,X^{\vartheta-1/4+\varepsilon}\sum_{X'+S<n\leqslant2X}\frac X{\left(n-X'\right)^2}\\
&\ll k^{3/2}\,M^{3/4}\,\Delta^{-1}\,X^{\vartheta+3/4+\varepsilon}\,S^{-1}
\ll k^3\,\left|\eta\right|^{3/2}\,M^{3/2}\,\Delta^{-1}\,X^{\vartheta+\varepsilon}\,S^{-1}.
\end{align*}
The first terms are estimated similarly, but with a dyadic split over the range of $n$:
\begin{align*}
&k^{-1/2}\sum_{n<X-S}\dots\\
&\ll k^{-1/2}\sum_{n<X-S}\left|t(n)\right|n^{-1/4}\cdot M^{-1/4}\,\Delta^{-P}\left(\frac{\sqrt n}{k\sqrt{M+\Delta}}-\left|\eta\right|\right)^{\!-P}\Delta\\
&\ll k^{-1/2}\,M^{-1/4}\,\Delta^{-1}\,k^2\,M\sum_{\substack{L\leqslant X-S\\\mathrm{dyadic}}}\sum_{L\leqslant n<2L}\frac{\left|t(n)\right|}{n^{1/4}\left|\sqrt n-\sqrt{X}\right|^2}\\
&\ll k^{3/2}\,M^{3/4}\,\Delta^{-1}\,X\sum_{\substack{L\leqslant X-S\\\mathrm{dyadic}}}L^{\vartheta-1/4+\varepsilon}\sum_{L\leqslant n<2L}\frac1{\left|n-X\right|^2}\\
&\ll k^{3/2}\,M^{3/4}\,\Delta^{-1}\,X^{\vartheta+3/4+\varepsilon}\,S^{-1}
\ll k^3\,\left|\eta\right|^{3/2}\,M^{3/2}\,\Delta^{-1}\,X^{\vartheta+\varepsilon}\,S^{-1}.
\end{align*}
The middle terms are estimated by absolute values to get
\begin{align*}
k^{-1/2}\sum_{X-S\leqslant n\leqslant X'+S}
&\ll k^{-1/2}\sum_{X-S\leqslant n\leqslant X'+S}\left|t(n)\right|n^{-1/4}\cdot\Delta\,M^{-1/4}\\
&\ll k^{-1/2}\,\Delta\,M^{-1/4}\left(S+k^2\,\eta^2\,\Delta\right)X^{-1/4+\vartheta+\varepsilon}.
\end{align*}
When $S\leqslant X/2$, say, then the last estimate follows immediately from estimate by absolute values. If $S>X/2$, we may use a dyadic split to estimate:
\begin{align*}
\sum_{X-S\leqslant n\leqslant X'+S}\frac{\left|t(n)\right|}{n^{1/4}}
&\ll\sum_{\substack{X-S\ll L\ll X'+S\\\mathrm{dyadic}}}\sum_{L<n\leqslant2L}\frac{\left|t(n)\right|}{n^{1/4}}\\
&\ll\sum_{\substack{X-S\ll L\ll X'+S\\\mathrm{dyadic}}}L^{\vartheta+3/4+\varepsilon}
\ll X^{\vartheta+3/4+\varepsilon}\ll S\,X^{\vartheta-1/4+\varepsilon}.
\end{align*}
Finally, if $S\gg k^2\,\eta^2\,M$, then the first terms and third terms do not exist, and the middle terms satisfy the same upper bound as before.

We have obtained
\begin{multline*}
\sum_{M\leqslant n\leqslant M+\Delta}t(n)\,e(n\alpha)\,w(n)
\ll k^3\left|\eta\right|^{3/2}M^{3/2}\,\Delta^{-1}\,X^{\vartheta+\varepsilon}\,S^{-1}\\
+k^{-1/2}\,\Delta\,M^{-1/4}\left(S+k^2\,\eta^2\,\Delta\right)X^{-1/4+\vartheta+\varepsilon}+M^{9/8-\beta-\gamma/4}+1.
\end{multline*}

\begin{lemma}\label{weight-function-removal-from-last-voronoi-argument}
Let $M\in\left[1,\infty\right[$ and $\Delta\in\left[1,M\right]$ with $M\gg1$ and $M^{2/3}\ll\Delta\ll M^{3/4}$, and let $\alpha\in\mathbb R$, $h\in\mathbb Z$, $k\in\mathbb Z_+$ and $\eta\in\mathbb R$ with $\alpha=h/k+\eta$, $(h,k)=1$, $k\leqslant M^{1/4}$ and $\left|\eta\right|\leqslant k^{-1}\,M^{-1/4}$.

If $k^2\,\eta^2\,M\leqslant1/2$, then
\[
\sum_{M\leqslant n\leqslant M+\Delta}t(n)\,e(n\alpha)
\ll M^{11/24}+M^{4/9+\vartheta/3+\varepsilon}.
\]

If $k^2\,\eta^2\,M\gg1$, then
\[
\sum_{M\leqslant n\leqslant M+\Delta}t(n)\,e(n\alpha)
\ll M^{11/24}+M^{4/9+\vartheta/3+\varepsilon}+M^{3/8+\gamma/4+\gamma\vartheta+\varepsilon}+\Delta\,M^{-1/4+\mu+\varepsilon},
\]
provided that $\Delta\ll M^{1+\mu-3\gamma/4-\gamma\vartheta}$, where $\gamma\in\mathbb R$ is such that $k^2\,\eta^2\,M=M^\gamma$, and $\mu\in\left[0,\infty\right[$.
\end{lemma}

\paragraph{Proof.} Let us use the same weight functions $w_\ell$ as in the proof of Theorem \ref{short-estimate} (cf. also Fig. \ref{partition-of-unity}), and let us pick $L\in\mathbb Z_+$ so that $\Delta\,4^{-L}\asymp M^{2/3}$. Theorem \ref{short-estimate} says that
\begin{multline*}
\sum_{M\leqslant n\leqslant M+\Delta}t(n)\,e(n\alpha)\left(1-\sum_{\ell=-L}^Lw_\ell(n)\right)\\
\ll\left(M^{2/3}\right)^{1/6-\vartheta}\,M^{1/3+\vartheta+\varepsilon}\ll M^{4/9+\vartheta/3+\varepsilon}.
\end{multline*}

Next, let us consider a single value $\ell\in\left\{-L,-L+1,\ldots,L\right\}$. If $k^2\,\eta^2\,M\leqslant1/2$, then the previous lemma gives the upper bound
\[\sum_{n\in\mathbb Z}t(n)\,e(n\alpha)\,w_\ell(n)\ll M^{9/8-2/3}\ll M^{11/24},\]
so assume that $k^2\,\eta^2\,M\gg1$.

If $k^2\left|\eta\right|M\,\Delta^{-1}\gg k^2\,\eta^2\,\Delta$, then we apply the previous lemma with $S=k^2\left|\eta\right|M\,\Delta^{-1}$ which has been optimized so that the first two terms in the upper bound coincide. We get
\begin{multline*}\sum_{n\in\mathbb Z}t(n)\,e(n\alpha)\,w_\ell(n)\\
\ll k^{-1/2}\,\Delta\,M^{-1/4}\,k^2\left|\eta\right|M\,\Delta^{-1}\left(k^2\,\eta^2\,M\right)^{-1/4+\vartheta+\varepsilon}+M^{9/8-2/3}.
\end{multline*}
Writing $\left|\eta\right|=M^{(\gamma-1)/2}\,k^{-1}$, this is
\begin{align*}
\sum_{n\in\mathbb Z}t(n)\,e(n\alpha)\,w_\ell(n)
&\ll k^{1/2}\,M^{1/4}\left(k^2\,\eta^2\,M\right)^{1/4+\vartheta+\varepsilon}+M^{11/24}\\
&\ll M^{3/8+\gamma/4+\gamma\vartheta+\varepsilon}+M^{11/24}.
\end{align*}
If $k^2\left|\eta\right|M\,\Delta^{-1}\ll k^2\,\eta^2\,\Delta$, then we apply the previous lemma with $S=k^2\,\eta^2\,\Delta$ and get
\begin{multline*}
\sum_{n\in\mathbb Z}t(n)\,e(n\alpha)\,w_\ell(n)
\ll k\left|\eta\right|^{-1/2}\Delta^{-2}\,M^{3/2}\left(k^2\,\eta^2\,M\right)^{\vartheta+\varepsilon}\\
+k^{3/2}\,\eta^2\,M^{-1/4}\,\Delta^2\left(k^2\,\eta^2\,M\right)^{\vartheta-1/4+\varepsilon}+M^{9/8-2/3}.
\end{multline*}
Now $M\ll\left|\eta\right|\Delta^2$, so that the first term is
\[\ll k\left|\eta\right|^{3/2}\Delta^2\,M^{-1/2}\left(k^2\,\eta^2\,M\right)^{\vartheta+\varepsilon},\]
and the second term is, writing again $\left|\eta\right|=M^{(\gamma-1)/2}\,k^{-1}$,
\begin{align*}
&=k\,\left|\eta\right|^{3/2}\Delta^2\,M^{-1/2}\left(k^2\,\eta^2\,M\right)^{\vartheta+\varepsilon}
\ll k^{-1/2}\,\Delta^2\,M^{3\gamma/4-5/4+\gamma\vartheta+\varepsilon}.
\end{align*}
This is $\ll\Delta\,M^{-1/4+\mu+\varepsilon}$, provided that $\Delta\ll M^{1+\mu-3\gamma/4-\gamma\vartheta}$,
and we are done.

\paragraph{Proof of Theorem \ref{longer-short-estimates}.}
We shall get the result by combining Lemmas \ref{corollary-to-approximate-functional-equation} and \ref{weight-function-removal-from-last-voronoi-argument}. To optimize the terms involving $\gamma$, we choose $\gamma_0$ so that
\[\frac12-\left(\frac1{12}-\frac\vartheta2\right)\gamma_0=\frac38+\frac{\gamma_0}4+\gamma_0\vartheta,\quad\text{i.e.}\quad \gamma_0=\frac3{12\vartheta+8}.\]

So, let $\alpha=h/k+\eta$ with $h$ and $k$ coprime integers, $1\leqslant k\leqslant M^{1/4}$, and $\eta\in\mathbb R$ with $\left|\eta\right|\leqslant k^{-1}\,M^{-1/4}$. When $k^2\,\eta^2\,M\gg M^{\gamma_0}$, Lemma \ref{corollary-to-approximate-functional-equation} gives
\[\sum_{M\leqslant n\leqslant M+\Delta}t(n)\,e(n\alpha)
\ll\Delta^{1/6-\vartheta}\,M^{1/3+\vartheta+\varepsilon}
+M^{3/8+(3+12\vartheta)/(32+48\vartheta)+\varepsilon}.\]
When $\Delta\ll M^{1+\mu-(3/4+\vartheta)\gamma_0}$ and $k^2\,\eta^2\,M\ll M^{\gamma_0}$, Lemma \ref{weight-function-removal-from-last-voronoi-argument} gives
\begin{multline*}
\sum_{M\leqslant n\leqslant M+\Delta}t(n)\,e(n\alpha)
\ll M^{11/24}+M^{4/9+\vartheta/3+\varepsilon}\\+M^{3/8+(3+12\vartheta)/(32+48\vartheta)+\varepsilon}
+\Delta\,M^{-1/4+\mu+\varepsilon},
\end{multline*}
where the exponent $\mu\in\left[0,\infty\right[$ is to be chosen later.
It is easy to check that
\[\max\left\{\frac{11}{24},\frac49+\frac\vartheta3\right\}\leqslant\frac38+\frac{3+12\vartheta}{32+48\vartheta}\quad\text{for}\quad\vartheta\in\left[0,\frac7{64}\right].\]
Also, it is easy to check that for $\Delta\ll M^{1+\mu-(3/4+\vartheta)\gamma_0}$, we have
\[\Delta^{1/6-\vartheta}\,M^{1/3+\vartheta+\varepsilon}
\ll M^{3/8+(3+12\vartheta)/(32+48\vartheta)+\varepsilon}\quad\text{for}\quad
\vartheta\in\left[0,\frac7{64}\right],\]
provided that $\mu\leqslant(3+12\vartheta)/(32+48\vartheta)$.
Combining the above facts gives the estimate
\[\sum_{M\leqslant n\leqslant M+\Delta}t(n)\,e(n\alpha)
\ll M^{3/8+(3+12\vartheta)/(32+48\vartheta)+\varepsilon}+\Delta\,M^{-1/4+\mu+\varepsilon},\]
for $\Delta\ll M^{1+\mu-(3/4+\vartheta)\gamma_0}$.
However, when $\Delta\asymp M^{1+\mu-(3/4+\vartheta)\gamma_0}$, it is easy to check that
\[M^{3/8+(3+12\vartheta)/(32+48\vartheta)+\varepsilon}\asymp\Delta\,M^{-1/4+\mu+\varepsilon},\quad\text{for any fixed}\quad\vartheta\in\left[0,\frac7{64}\right],\]
with equal exponents $\varepsilon$, of course, if we choose
\[\mu=\frac{3\vartheta}{32+48\vartheta}.\]
This choice trivially satisfies the required upper bound $\leqslant(3+12\vartheta)/(32+48\vartheta)$.
Thus, by splitting longer sums into sums of length $M^{1+\mu-(3/4+\vartheta)\gamma_0}$, and estimating these subsums separately, we have
\[\sum_{M\leqslant n\leqslant M+\Delta}t(n)\,e(n\alpha)
\ll M^{3/8+(3+12\vartheta)/(32+48\vartheta)+\varepsilon}+\Delta\,M^{-1/4+\mu+\varepsilon},\]
with the weaker condition $\Delta\ll M^{3/4}$.

\section{$\Omega$-results from second moments}

We end with some details about the $\Omega$-results. Two important results are actually mean square results for which the key is the following truncated Voronoi identity due to Meurman \cite{Meurman1988}.
\begin{theorem}\label{truncated-Voronoi-for-Maass-waves}
Let $x\in\left[1,\infty\right[$, let $N\in\mathbb R_+$ be such that $N\ll x$, let $k$ be a positive integer with $k\ll\sqrt x$, and let $h$ be an integer coprime to $k$. Then
\begin{multline*}
\sum_{n\leqslant x}t(n)\,e\!\left(\frac{nh}k\right)
=\frac1{\pi\,\sqrt2}\,k^{1/2}\,x^{1/4}\sum_{n\leqslant N}t(n)\,e\!\left(\frac{-n\overline h}k\right)n^{-3/4}\\
\cdot\cos\!\left(\frac{4\pi\sqrt{nx}}k-\frac\pi4\right)+O\!\left(k\,x^{1/2+\vartheta+\varepsilon}\,N^{-1/2}\right).
\end{multline*}
\end{theorem}
\noindent
This is simpler than the formulation in \cite{Meurman1988} where care was taken to retain an explicit dependence on~$\psi$. For a fixed $\psi$, the above formulation follows easily. This simplified formulation can be used in the same manner as the truncated Voronoi identities for the error terms of the Dirichlet divisor problem or circle problem, or for holomorphic cusp forms.

For long exponential sums with rational additive twists, we have the following result, which is a Maass form analogue of Theorem 1.2 in \cite{Jutila1987a}.
\begin{theorem}
For $M\in\left[1,\infty\right[$, and for coprime integers $h$ and $k$ for which $1\leqslant k\leqslant M$, we have
\[\int\limits_1^M\left|\sum_{n\leqslant x}t(n)\,e\!\left(\frac{nh}k\right)\right|^2\mathrm dx=A\,k\,M^{3/2}+O(k^2\,M^{1+2\vartheta+\varepsilon})+O(k^{3/2}\,M^{5/4+\vartheta+\varepsilon}),\]
where
\[A=\frac1{6\pi^2}\sum_{n=1}^\infty\frac{\left|t(n)\right|^2}{n^{3/2}}.\]
\end{theorem}
For sufficiently small $k$, the main term dominates:
\begin{corollary}
For $M\in\left[1,\infty\right[$ and for coprime integers $h$ and $k$ for which $1\leqslant k\ll M^{1/2-2\vartheta-\varepsilon}$, we have
\[\int\limits_1^M\left|\sum_{n\leqslant x}t(n)\,e\!\left(\frac{nh}k\right)\right|^2\mathrm dx\asymp k\,M^{3/2}.\]
\end{corollary}
\noindent
In particular, these sums are $\Omega(k^{1/2}\,M^{1/4})$.

For shorter sums we have the following theorem, which is similar to Theorem 6 in \cite{Vesalainen2014}.
\begin{theorem}\label{short-mean-square}
Let $M\in\left[1,\infty\right[$, let $\Delta\in\mathbb R_+$ with $M^\varepsilon\ll\Delta\ll M^{1/2-\varepsilon}$, and let $h$ and $k$ be coprime integers with $1\leqslant k\ll\Delta^{1/2-\varepsilon}\,M^{-\vartheta}$. Then
\[\int\limits_M^{2M}\left|\sum_{x\leqslant n\leqslant x+\Delta}t(n)\,e\!\left(\frac{nh}k\right)\right|^2\mathrm dx\asymp M\,\Delta.\]
\end{theorem}
\noindent
In particular, the sum in question is $\Omega(\Delta^{1/2})$ when $\Delta\gg M^{2\vartheta+\varepsilon}$.

This follows for intance by following the proof of Theorem 6 in \cite{Vesalainen2014}: the main difference is that the error terms $k^2\,\Xi\,M^{\varepsilon}+k\,\Xi\,\Delta^{1/2}\,M^{\varepsilon}$ are to be replaced by $k^2\,\Xi\,M^{2\vartheta+\varepsilon}+k\,\Xi\,\Delta^{1/2}\,M^{\vartheta+\varepsilon}$. Except for the error term of the truncated Voronoi identity, the exponent $\vartheta$ never appears as the Fourier coefficients are estimated by the Rankin--Selberg estimate.

\section*{Acknowledgements}

The authors would like to express their deep gratitude for the encouragement and support of Dr. A.-M. Ernvall-Hyt\"onen.

During this research, the first author was funded by the Academy of Finland project Number Theory Finland and the Doctoral Programme for Mathematics and Statistics of the University of Helsinki. The second author was funded by Finland's Ministry of Education through the Doctoral School of Inverse Problems, Academy of Finland's Centre of Excellence in Inverse Problems Research and the Foundation of Vilho, Yrj\"o and Kalle V\"ais\"al\"a.

\end{document}